\documentclass[10pt]{article}

\usepackage[utf8]{inputenc}

\usepackage{amsmath}
\usepackage{amsthm}
\usepackage{amssymb}  

\usepackage{varioref}       
\usepackage[hidelinks]{hyperref} 
\usepackage{cleveref}   

\usepackage{caption}
\usepackage{subcaption}
\usepackage{graphicx}  
\usepackage{placeins}  
\usepackage{grffile}
\usepackage[usenames,dvipsnames]{xcolor}  
\usepackage{listings} 
\usepackage{xcolor}

\usepackage{tikz}
\usetikzlibrary{shapes.geometric, arrows}
\usetikzlibrary{calc}

\usepackage{cancel}

\usepackage{xargs} 

\usepackage{datatool, filecontents}
\DTLsetseparator{;}
\DTLloadrawdb[noheader, keys={thekey,thevalue}]{runsinfo}{runsinfo.csv}
\newcommand{\perplexityinsert}[1]{\DTLgetvalueforkey{\datavalue}{thevalue}{runsinfo}{thekey}{#1}\datavalue}

\usepackage{csvsimple}
\makeatletter
\csvset{
  autotabularcenter/.style={
    file=#1,
    after head=\csv@pretable\begin{tabular}{|*{\csv@columncount}{c|}}\csv@tablehead,
    table head=\hline\csvlinetotablerow\\\hline,
    late after line=\\,
    table foot=\\\hline,
    late after last line=\csv@tablefoot\end{tabular}\csv@posttable,
    command=\csvlinetotablerow},
}
\makeatother

\definecolor{codegreen}{rgb}{0,0.6,0}
\definecolor{codegray}{rgb}{0.5,0.5,0.5}
\definecolor{codepurple}{rgb}{0.58,0,0.82}
\definecolor{backcolour}{rgb}{0.95,0.95,0.92}

\lstdefinestyle{mystyle}{
    backgroundcolor=\color{backcolour},   
    commentstyle=\color{codegreen},
    keywordstyle=\color{magenta},
    numberstyle=\tiny\color{codegray},
    stringstyle=\color{codepurple},
    basicstyle=\ttfamily\footnotesize,
    breakatwhitespace=false,         
    breaklines=true,                 
    captionpos=b,                    
    keepspaces=true,                 
    numbers=left,                    
    numbersep=5pt,                  
    showspaces=false,                
    showstringspaces=false,
    showtabs=false,                  
    tabsize=2
}

\newcommand\e{\varepsilon}

\renewcommand\t{\tilde}

\renewcommand\>{\rangle}
\renewcommand\leq{\leqslant}
\renewcommand\geq{\geqslant}
\def\gsim{\mathrel{\raisebox{-4pt}{$\stackrel{\textstyle>}{\sim}$}}}

 \def\cM{{\cal M}}

 \def\cT{{\cal T}}
 
  \def\cS{{\cal S}}
 \def\cF{{\cal F}}
 \def\cC{{\cal C}}
 \def\cL{{\cal L}}
 \def\cK{{\cal K}}
  \def\cO{{\cal O}}

 \def\o{\overline}
 \def\nl{\newline}
 \def\vp{\varphi}
 \def\({\left (}
 \def\){\right )}

\newcommand{\R}{\mathbb R}

\newcommand{\pP}{\mathbb P}

\def\Chi{\raise .3ex\hbox{\large $\chi$}}
\newcommand{\be}{\begin{equation}}
\newcommand{\ee}{\end{equation}}
\newcommand{\iref}[1]{{\rm (\ref{#1})}}

\newtheorem{theorem}{Theorem}[section]

\newtheorem{remark}[theorem]{Remark}
\newtheorem{definition}[theorem]{Definition}
\newtheorem{corollary}[theorem]{Corollary}

\DeclareMathOperator*{\argmin}{arg\,min}

\newcommand{\corr}[1]{\textcolor{black}{#1}}

\definecolor{gridcolor}{RGB}{192, 192, 192} 
\definecolor{domaincolor}{RGB}{34, 139, 34} 

\def\edgeXcoord{2.3}

\newcommand{\functionA}[1]{((#1)^2/4+(#1)/3+2)}
\newcommand{\functionB}[1]{(-((#1)-3.2)^2/3+1.1)}
\newcommand{\functionAbracket}[1]{{\functionA{#1}}}
\newcommand{\functionBbracket}[1]{{\functionB{#1}}}

\newcommand{\PlotSubCellDiagramD}{
            
    \begin{tikzpicture}[node distance=\whalf cm, scale=1, every node/.style={transform shape}]
    
        \draw[ultra thick] plot [domain=-1:\edgeXcoord] (\x,\functionAbracket{\x});
        \draw[ultra thick] plot [domain=\edgeXcoord:4] (\x,\functionBbracket{\x});
                    
        \fill[fill=yellow!50] (-1, 0) -- plot [domain=-1:\edgeXcoord] (\x,\functionAbracket{\x}) -- (\edgeXcoord, 0) -- cycle;
        \fill[fill=yellow!50] (\edgeXcoord, 0) -- plot [domain=\edgeXcoord:4] (\x,\functionBbracket{\x}) -- (4, 0) -- cycle;
    
        \draw[ultra thick, dotted] (-1,0) -- (-1, \functionAbracket{-1});
        \draw[ultra thick, dotted] (0,0) -- (0, \functionAbracket{0});
        \draw[ultra thick, dotted] (1,0) -- (1, \functionAbracket{1});
        \draw[ultra thick, dotted] (2,0) -- (2, \functionAbracket{2});
        \draw[ultra thick, dotted] (3,0) -- (3, \functionBbracket{3});
        \draw[ultra thick, dotted] (4,0) -- (4, \functionBbracket{4});
        
        \draw[ultra thick, MidnightBlue, opacity=1] (-1,\functionAbracket{-0.5}) -- (0,\functionAbracket{-0.5});
        \draw[ultra thick, MidnightBlue, opacity=1] (0,\functionAbracket{0.5}) -- (1,\functionAbracket{0.5});
        \draw[ultra thick, MidnightBlue, opacity=1] (1,\functionAbracket{1.5}) -- (2,\functionAbracket{1.5});
        \newcommand{\average}{{\functionA{(2.0+\edgeXcoord)/2}*(\edgeXcoord-2)+\functionB{(\edgeXcoord+3.0)/2}*(3-\edgeXcoord)}}
        \draw[ultra thick, MidnightBlue, opacity=1] (2, \average) -- (3, \average);
        \draw[ultra thick, MidnightBlue, opacity=1] (3,\functionBbracket{3.5}) -- (4,\functionBbracket{3.5});
        
        \draw[thick, dashed, red, opacity=1] (\edgeXcoord,0) -- (\edgeXcoord, \functionAbracket{\edgeXcoord});
        \node[red, opacity=1] at (\edgeXcoord, -0.5) {$X_0$};
        
        \node[MidnightBlue, opacity=1] at (-0.5, -0.5) {$\mathbf{i-1}$};
        \node[ForestGreen, opacity=1] at (0.5, -0.5) {$\mathbf{i}$};
        \node[MidnightBlue, opacity=1] at (1.5, -0.5) {$\mathbf{i+1}$};
        \node[MidnightBlue, opacity=1] at (3.5, -0.5) {$\mathbf{i+3}$};
    
    \end{tikzpicture}
}

\newcommand*\samethanks[1][\value{footnote}]{\footnotemark[#1]}

\title{High order recovery of geometric interfaces from cell-average data}
\author{Albert Cohen\thanks{Laboratoire Jacques-Louis Lions, Sorbonne Universit\'e, 
4 place Jussieu, 75005 Paris, France (albert.cohen@sorbonne-universite.fr and agustin.somacal@sorbonne-universite.fr)}, 
Olga Mula\thanks{TU Eindhoven, Department of Mathematics and Computer Science, 5612 AZ Eindhoven, Netherlands (o.mula@tue.nl)}, Agustin Somacal\samethanks[1]}
\date{\today}

\begin{document}
\begin{filecontents*}{runsinfo.csv}
piecewise_constant;Piecewise Constant
nn_linear;NN Linear
elvira;ELVIRA
elvira_w;ELVIRA-W
elvira_w_oriented;ELVIRA-W Oriented
linear_obera;OBERA Linear
linear_aero;AERO Linear
linear_aero_w;AERO-W Linear
quadratic_obera_non_adaptive;OBERA Quadratic 3x3
quadratic_obera;OBERA Quadratic
quadratic_aero;AERO Quadratic
piecewise-constant;Piecewise constant
nn-linear;NN Linear
elvira-w;ELVIRA-W
elvira-w-oriented;ELVIRA-WO
linear-obera;OBERA Linear
linear-aero;AEROS Linear
linear-aero-w;AEROS-W Linear
quadratic-obera-non-adaptive;OBERA Quadratic
quadratic-obera;OBERA Quadratic l2
quadratic-aero;AEROS Quadratic
linear-aero-consistent;AEROS Linear Column Consistent
linear-obera-w;LVIRA
obera-circle;OBERA Circle
obera-circle-vander;OBERA Circle ReParam
nn-linear-torch;NN Linear Torch
nn-linear-keras;NN Linear Keras
cceweight;100
circler;0.232
circlex;0.511
circley;0.486
upwind;UpWind
elvira-oriented;ELVIRA Oriented
aero-linear;AERO-Linear
aero-linear-oriented;AERO Linear Oriented
quadratic;AERO-Quadratic
aero-lq;AERO-Linear-Quadratic
aero-lq-vertex;AERO-Quadratic Vertex
obera-aero-lq-vertex;OBERA-AERO Linear Quadratic Vertex
quadratic-oriented;AERO Quadratic Oriented
qelvira;ELVIRA AERO-Quadratic
elvira-time;0.04
elvira-w-oriented-time;0.02
linear-aero-time;0.0077
linear-aero-consistent-time;0.0027
linear-aero-w-time;0.0075
linear-obera-time;0.9
linear-obera-w-time;0.8
nn-linear-time;0.0011
obera-circle-time;0.5075
obera-circle-vander-time;0.2472
piecewise-constant-time;nan
quadratic-aero-time;0.005
quadratic-obera-time;0.3002
quadratic-obera-non-adaptive-time;0.5
aero-qelvira-vertex;ELVIRA-WO + AEROS Quadratic + TEM + AEROS Vertex
aero-qelvira-vertex45;AERO-Quadratic Orient Vertex
poly02half;Quadratic
poly0-elvira;Constant + ELVIRA
poly02h-elvira;Constant + ELVIRA + Quadratic
poly02h-qelvira;Constant + QELVIRA + Quadratic
poly02h-quadratic;Constant + AERO-QUADRATIC + Quadratic
poly2h-elvira;Quadratic + ELVIRA
poly2h-qelvira;Quadratic + QELVIRA
poly2h-quadratic;Quadratic + AERO-QUADRATIC
std-elvira-time;0.009
std-elvira-w-oriented-time;0.002
std-linear-aero-time;0.0083
std-linear-aero-consistent-time;0.01
std-linear-aero-w-time;0.0086
std-linear-obera-time;0.4
std-linear-obera-w-time;0.3
std-nn-linear-time;0.0001
std-obera-circle-time;0.3736
std-obera-circle-vander-time;0.2469
std-piecewise-constant-time;nan
std-quadratic-aero-time;0.006
std-quadratic-obera-time;0.1001
std-quadratic-obera-non-adaptive-time;0.2
qlow-elvira-time;0.03
qlow-elvira-w-oriented-time;0.01
qlow-linear-aero-time;0.0024
qlow-linear-aero-consistent-time;0.0015
qlow-linear-aero-w-time;0.0024
qlow-linear-obera-time;0.5
qlow-linear-obera-w-time;0.4
qlow-nn-linear-time;0.001
qlow-obera-circle-time;0.2589
qlow-obera-circle-vander-time;0.1217
qlow-piecewise-constant-time;nan
qlow-quadratic-aero-time;0.002
qlow-quadratic-obera-time;0.1786
qlow-quadratic-obera-non-adaptive-time;0.2
qhigh-elvira-time;0.05
qhigh-elvira-w-oriented-time;0.02
qhigh-linear-aero-time;0.0268
qhigh-linear-aero-consistent-time;0.0021
qhigh-linear-aero-w-time;0.0269
qhigh-linear-obera-time;2
qhigh-linear-obera-w-time;1
qhigh-nn-linear-time;0.0014
qhigh-obera-circle-time;0.818
qhigh-obera-circle-vander-time;0.7597
qhigh-piecewise-constant-time;nan
qhigh-quadratic-aero-time;0.02
qhigh-quadratic-obera-time;0.4859
qhigh-quadratic-obera-non-adaptive-time;0.8
median-elvira-time;0.04
median-elvira-w-oriented-time;0.02
median-linear-aero-time;0.0028
median-linear-aero-consistent-time;0.0017
median-linear-aero-w-time;0.0028
median-linear-obera-time;0.9
median-linear-obera-w-time;0.7
median-nn-linear-time;0.0011
median-obera-circle-time;0.4479
median-obera-circle-vander-time;0.1826
median-piecewise-constant-time;nan
median-quadratic-aero-time;0.003
median-quadratic-obera-time;0.2768
median-quadratic-obera-non-adaptive-time;0.5
nn-linear-keras-time;0.0016
std-nn-linear-keras-time;0.0003
qlow-nn-linear-keras-time;0.001
qhigh-nn-linear-keras-time;0.002
median-nn-linear-keras-time;0.0016
nn-linear-skkeras-time;0.0326
std-nn-linear-skkeras-time;0.0079
qlow-nn-linear-skkeras-time;0.0284
qhigh-nn-linear-skkeras-time;0.041
median-nn-linear-skkeras-time;0.0307
nn-linear-skkeras;NN Linear SKKeras
sknn-fluxlines;skNN Flux Lines
fluxlines;NN Flux Lines
ml-qelvira;ML ELVIRA AERO-Quadratic
nn-fluxlines;NN Flux Lines
nn-fluxquadratics;NN Flux Quadratics
nn-fluxlinesquadratics;NN Flux Lines-Quadratics
ml-vql;ML ELVIRA AERO-Quadratic Vertex
obera-circle-vander-ml;OBERA Circle ReParam ML-ker
quadratic-obera-ml-time;0.1279
std-quadratic-obera-ml-time;0.0642
qlow-quadratic-obera-ml-time;0.0639
qhigh-quadratic-obera-ml-time;0.2578
median-quadratic-obera-ml-time;0.1122
quadratic-obera-ml;OBERA Quadratic ML-ker
nn-quadratic-3x3;NN Quadratic 3x3
nn-quadratic-3x7;NN Quadratic 3x7
nn-quadratic-3x7-params;NN Quadratic 3x7 params
nn-quadratic-3x3-time;0.0013
nn-quadratic-3x7-time;0.0014
nn-quadratic-3x7-params-time;0.0012
std-nn-quadratic-3x3-time;0.0001
std-nn-quadratic-3x7-time;0.0002
std-nn-quadratic-3x7-params-time;0.0002
qlow-nn-quadratic-3x3-time;0.0011
qlow-nn-quadratic-3x7-time;0.0012
qlow-nn-quadratic-3x7-params-time;0.0011
qhigh-nn-quadratic-3x3-time;0.0015
qhigh-nn-quadratic-3x7-time;0.0018
qhigh-nn-quadratic-3x7-params-time;0.0016
median-nn-quadratic-3x3-time;0.0012
median-nn-quadratic-3x7-time;0.0014
median-nn-quadratic-3x7-params-time;0.0012
elvira-w-oriented-ml;ELVIRA-W Oriented ML-ker l2
elvira-w-oriented-ml-time;0.0028
std-elvira-w-oriented-ml-time;0.0003
qlow-elvira-w-oriented-ml-time;0.0024
qhigh-elvira-w-oriented-ml-time;0.0033
median-elvira-w-oriented-ml-time;0.0028
obera-circle-vander-ml-time;0.2745
std-obera-circle-vander-ml-time;0.2601
qlow-obera-circle-vander-ml-time;0.1241
qhigh-obera-circle-vander-ml-time;0.9165
median-obera-circle-vander-ml-time;0.1979
nn-quadratic-3x7-adapt;NN Quadratic 3x7 Adapt S
nn-quadratic-3x7-params-adapt;NN Quadratic 3x7 params Adapt S
elvira-w-oriented-l1;ELVIRA-W Oriented l1
linear-obera-w-l2;OBERA-W Linear l2
linear-obera-w-ml;OBERA-W Linear ML-ker l1
quadratic-obera-non-adaptive-l1;OBERA Quadratic 3x3 l1
quadratic-obera-l1;OBERA Quadratic l1
quadratic-obera-ml-params;OBERA Quadratic ML-ker 3x3
quadratic-obera-ml-points-adapt;OBERA Quadratic ML-ker 3x3 ReParam Adapt S
quadratic-obera-ml-points;OBERA Quadratic ML-ker 3x3 ReParam
cubic-aero;AEROS Cubic
cubic-aero-stencil4;AEROS Cubic
cubic-obera;OBERA Cubic
quartic-aero;AEROS Quartic
quartic-obera;OBERA Quartic
color-piecewise-constant;0.89, 0.47, 0.76
color-linear-obera;0.12, 0.47, 0.71
color-elvira;0.12, 0.47, 0.71
color-elvira-w-oriented;0.58, 0.4, 0.74
color-linear-obera-w;0.58, 0.4, 0.74
color-quadratic-aero;0.17, 0.63, 0.17
color-quadratic-obera-non-adaptive;0.17, 0.63, 0.17
color-quartic-aero;0.84, 0.15, 0.16
quartic-aero-time;0.007
std-quartic-aero-time;0.005
qlow-quartic-aero-time;0.003
qhigh-quartic-aero-time;0.02
median-quartic-aero-time;0.004
aero_qelvira_vertex;ELVIRA + TEM + AEROS Quadratic + AVROS
aero_qelvira_tem;ELVIRA + TEM + AEROS Quadratic
aero_qvertex;AEROS Quadratic + TEM + AVROS
aero_qtem;AEROS Quadratic + TEM
aero-qelvira-tem;ELVIRA-WO + AEROS Quadratic + TEM
aero-qvertex;AEROS Quadratic + TEM + AEROS Vertex
aero-qtem;AEROS Quadratic + TEM
obera-qtem;OBERA Quadratic + TEM
qelvira-kml;ML-Kernel ELVIRA AERO-Quadratic
corners-aero-qelvira-vertex-time;0.05
corners-aero-qtem-time;0.004
corners-aero-qvertex-time;0.04
corners-quadratic-aero-time;0.001
scheme-aero-qelvira-vertex-time;3
scheme-elvira-w-oriented-time;0.5
scheme-quadratic-aero-time;0.2
scheme-upwind-time;0.09
\end{filecontents*}

\maketitle

\begin{abstract}
We consider the problem of recovering characteristic functions $u:=\Chi_\Omega$ from cell-average data on a coarse grid, and where $\Omega$ is a compact set of $\R^d$. This task arises in very different contexts such as image processing, inverse problems, and the accurate treatment of interfaces in finite volume schemes.
While linear recovery methods are known to perform poorly, nonlinear strategies based on local reconstructions of the jump interface 
$\Gamma:=\partial\Omega$ by geometrically simpler interfaces may offer significant improvements.
We study two main families of local reconstruction schemes,
the first one based on nonlinear least-squares fitting,
the second one based on the explicit computation of
a polynomial-shaped curve fitting the data, which yields
simpler numerical computations and high order geometric fitting. For each of them,
we derive a general theoretical framework which allows us to control the recovery error by the error of best approximation up to a fixed multiplicative constant. Numerical tests in 2d illustrate the expected approximation order of these strategies. Several extensions are discussed, in particular
the treatment of piecewise smooth interfaces with corners.
\end{abstract}


\section{Introduction}

\subsection{Reconstruction from cell-averages}
We consider the problem of reconstructing a function $u:D\to \R$ defined
on a multivariate domain $D\subset \R^d$ from 
cell averages
\begin{equation}
\label{eq:local-average}
a_T(u):=\frac 1 {|T|} \int_T u(x)dx, \quad T\in \cT,
\end{equation}
over a partition $\cT$ of $D$. This task occurs in various contexts, the most notable ones being:
\begin{enumerate}
\item 
\textbf{Image processing:} here $u$ is the light intensity of an image and $\cT$ represents a grid of pixels in dimension $d=2$ or voxels in dimension $d=3$. Various processing tasks are facilitated by the reconstruction of the image at the continuous level,  for example when applying operations that are not naturally compatible with the pixel grid such as rotations, or when changing the 
format of the pixel grid such as in super-resolution.
\item
\textbf{Hyperbolic transport PDE's:} here $u$ is a solution to such
an equation and $\cT$ is a computational grid, typically in dimension $d=1,2$ or $3$. 
Finite volume schemes 
evolve the cell average data by computing at each time step 
the numerical fluxes at the interfaces between each cell. Several such
schemes are based on an intermediate step that reconstructs
simple approximations to $u$ on each cell and compute the numerical fluxes by 
applying the transport operator to these approximations.
\item \textbf{Inverse Problems:} numerous inversion tasks can be formulated as the recovery of a function from observational data, and this data could typically come in the form of local averages of the type \eqref{eq:local-average}. 
\end{enumerate}
In this paper, we consider functions defined on the unit cube
$$
D:=[0,1]^d,
$$
and partitions $\cT_h$ of $D$ based on uniform cartesian meshes,
that is, consisting of cells of the form
$$
T=h(k+D),  \quad k:=(k_1,\dots,k_d)\in \{0,\dots,l-1\}^d,
$$
where $h:=\frac 1 l >0$ is the side-length 
of each cell in $\cT_h$, for some $l>1$. The cardinality of the partition
is therefore
$$
\corr{N}:=\#(\cT_h)=l^d=h^{-d}.
$$

We are thus interested in reconstruction operators $R$
that return an approximation $\t u=R(a)$
to $u$ from the $\corr{N}$-dimensional vector 
$a=a(u)=(a_T(u))_{T\in \cT_h}\in \R^{\corr{N}}$.
The most trivial one is the piecewise constant function
\be
\t u:=\sum_{T\in \cT_h} a_T(u) \Chi_T,
\label{eq:pwc}
\ee
which is for example used in the Godunov finite volume scheme.
Elementary arguments show that this
reconstruction is first order accurate: if $u\in W^{1,p}(D)$, one has
$$
\|u-\t u\|_{L^p(D)} \leq Ch \|\nabla u\|_{L^p(D)} \sim \corr{N}^{-\frac 1 d},
$$
where $C$ is a fixed constant, and the exponent 
in this estimate cannot be improved for smoother functions. 

A simple way to raise the order of accuracy 
is by reconstructing on each cell polynomials of higher degree
using neighbouring cell averages.
For example, in the univariate case $d=1$ and for some fixed $m\geq 1$,
we associate to each interval $T_k:=[kh,(k+1)h]$ the centered stencil 
consisting of the cells $T_l$ for $l=k-m,\dots,k+m$. Then, there exists
a unique polynomial $p_k\in \pP_{2m}$ such that
$$
a_{T_l}(p_k)=a_{T_l}(u), \quad l=k-m,\dots,k+m.
$$
We then define a piecewise polynomial reconstruction by
$$
\t u:=\sum_{k=0,\dots,\corr{N}-1} p_k \Chi_{T_k}.
$$
This strategy can be generalized to higher dimension $d>1$ in a straightforward
manner: for each cell $T$ we consider the stencil $S=S_T$
of $(2m+1)^d$ cells centered around $T$ and define the piecewise polynomial reconstruction
$$
\t u:=\sum_{T\in \cT_h} p_T \Chi_T,
$$
where $p_T$ is the unique polynomial of degree $2m$ in each variable
such that 
$$
a_{\t T}(p_T)=a_{\t T}(u), \quad \t T\in S_T.
$$
For example in the bivariate case $d=2$, we may use $3\times 3$ stencils 
to reconstruct bi-quadratic polynomials on each cell. Standard approximation
theory arguments show that these local reconstruction operators now satisfy 
accuracy estimates of the form
$$
\|u-\t u\|_{L^p(D)} \leq Ch^{r} |u|_{W^{r,p}(D)} \sim \corr{N}^{-\frac r d},
$$
for $r\leq 2m+1$. 

\begin{remark} Polynomials of odd degree can also be constructed by 
using non-centered stencils. Also note that non-centered stencils 
need to be used when approaching the boundary of $D$, but 
this does not affect the above estimate.
\end{remark}

These classical methods are therefore efficient to reconstruct
smooth functions with a rate of accuracy that optimally reflects their amount of smoothness. Unfortunately 
they are doomed to perform poorly in the case of functions $u$ that 
are piecewise smooth with jump discontinuities accross hypersurfaces. 
For example, a piecewise constant reconstruction will 
have ${\cal O}(1)$ error on each cell that is crossed by the interface.
Since the amount of such cells is of order \corr{$\corr{N}^{d-1}=h^{1-d}$}, we cannot
expect a reconstruction error better than 
\be
\|u-\t u\|_{L^p} \gsim (h^d h^{1-d})^{\frac 1 p}=h^{\frac 1 p}=\corr{N}^{-\frac 1 {dp}}.
\label{firstorder}
\ee
In particular, this reconstruction has first order accuracy $\cO(h)$
in the $L^1$ norm.

The use of higher order polynomial cannot improve this rate.
In fact, a fundamental obstruction is the fact that 
all the above methods produce approximations $\t u$
that depend linearly on $a(u)$, and therefore belong to a linear
space of dimension $\corr{N}$. The bottleneck of such methods for a given 
class of functions $\cK$ is therefore given by the so-called Kolmogorov 
$\corr{N}$-width defined by 
$$
d_{\corr{N}}(\cK)_{L^p} :=\inf_{{\rm dim}(V_{\corr{N}})=\corr{N}} \sup_{u\in \cK}\;\; \min_{v\in V_{\corr{N}}} \|u-v\|_{L^p}.
$$
Then, it can be shown that for very simple classes $\cK$ of discontinuous functions
such as those of the form $u=\Chi_H$, where $H$ is any half-space passing
through $D$,
the $\corr{N}$-width in $L^p$ precisely behaves like $\corr{N}^{-\frac 1 {dp}}$, see \cite{NonLinearReduced}. 
In summary, any linear method cannot do much better than the low order
piecewise constant method, even for interfaces that are infinitely smooth.

Improving the accuracy in the reconstruction of piecewise smooth functions
from cell averages therefore motivates the development and study of  
nonlinear reconstruction strategies, which is at the heart of this work. 
We first recall the main existing approaches.

\subsection{Reconstruction of discontinuous interfaces}

\begin{figure}
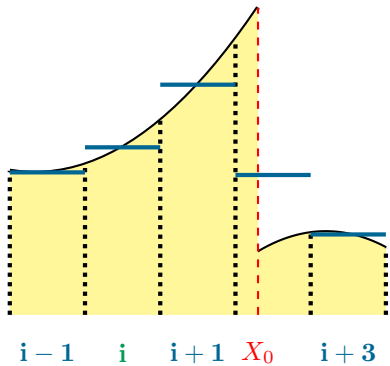

    \centering
    \PlotSubCellDiagramD
    \caption{ENO-SR in one dimension: the jump point $X_0$ is identified
    by matching the average on the singular cell with the piecewise polynomial reconstruction.}
    \label{subcell_schema}
\end{figure}

One first approach aiming to tackle jump discontinuities 
while maintaining high order approximation in the smooth regions
was proposed for the univariate case $d=1$ by Ami Harten in terms of ENO
(Essentially Non Oscillatory) and ENO-SR (Subcell Resolution).
The ENO strategy \cite{ENO1987} is based on selecting
for each cell $T$ a stencil $S_T$ that should not include the cell 
$T^*$ which contains the point of jump discontinuity. This is achieved
by choosing among the stencils that contains $T$ the one where the 
cell average values have the least numerical variation. For $T\neq T^*$
such adaptively selected stencils will tend to avoid $T^*$. 

As a consequence high order reconstruction can be preserved in all cells where $u$ is smooth, and  in addition this yields for free a singularity detection mechanism which identifies the singular cell $T^*$ that is avoided from both side by the stencil selection.  
The ENO-SR strategy \cite{ENO-SR1989} then consists in reconstructing 
in this singular cell by extending the polynomials fitted on both sides until a point 
for which the resulting average match the observed average. The position of this 
point can therefore be identified by solving a simple algebraic equation. This strategy
is very effective in the univariate case as illustrated in 
Figure \ref{subcell_schema}.

While the ENO stencil selection can be generalized to higher dimension \cite{Arandiga2008},
the ENO-SR strategy does not have a straightforward multivariate version.
This is due to the fact that, instead of a single point, the jump discontinuity to be identified
is now a hypersurface (curve in 2d, surface in 3d, ...) that cannot be described
by finitely many parameters. The approximate recovery of geometric interfaces
from cell-average has been the object of continuous investigation, with the 
particular focus on functions of the form
\begin{equation}
\label{eq:cartoon-function}
u=\Chi_\Omega,
\end{equation}
where $\Omega\subset D$ is a set with boundary $\Gamma=\partial \Omega$ having
a certain smoothness. For such characteristic functions, that could for example
represent a two phase flow without any mixing or the evolution of a front, the whole
difficulty is concentrated in the recovery of $\Gamma$ since the smooth parts
are the trivial constant functions $0$ or $1$.

In this paper, we denote by 
$$
S_h:=\{T\in \cT_h \; : \: |T\cap \Omega | \neq 0 \; {\rm and} \;  |T\cap \Omega^c | \neq 0\},
$$
the set of cells that are crossed by the interface $\Gamma$ in a non trivial way. Such cells are termed as {\it singular}, the other as {\it regular}. Let us note that, \corr{for functions $u$ of the form \eqref{eq:cartoon-function},} singular cells are characterized by the property 
$$
0<a_T(u)<1,
$$
and can thus be identified from the cell-average data.

Practical computational strategies, often termed as {\it volume of fluid} methods, 
consist in using the cell-average data to reconstruct a local 
approximation of the interface
that can be described by finitely many parameters, such as lines in 2d or planes in 3d.
This idea was introduced in \cite{SLIC}, and was significantly improved
in \cite{LVIRApuckett1991volume} with the LVIRA algorithm
that consists in reconstructing in each singular cell 
a linear interface whose parameters are found by least-square minimization of  
the difference between exact and reconstructed cell averages on 
centered $3\times 3$ stencil, in the 2d case.

\corr{Note} that this continuous least-square minimization 
is performed over a nonlinear set. This induces
a substantial computational burden that could 
be avoided through a variant, the ELVIRA algorithm, in which the
line selection is made between $6$ possible configurations explicitly computed
from the cell averages. One main result is the fact that this reconstruction returns
precisely the true interface if this one is indeed a straight line. We refer to 
\cite{Pilliod2004} for a comparative survey on these reconstruction algorithms \corr{ 
and to \cite{VOFSLIC, VtrackingSurvey, Vtracking, Zaleski, ZaleskiVOF, WY2010, DS2008, RV2008} for improvements
and applications in the domain of 2d and 3d fluid mechanics.}

 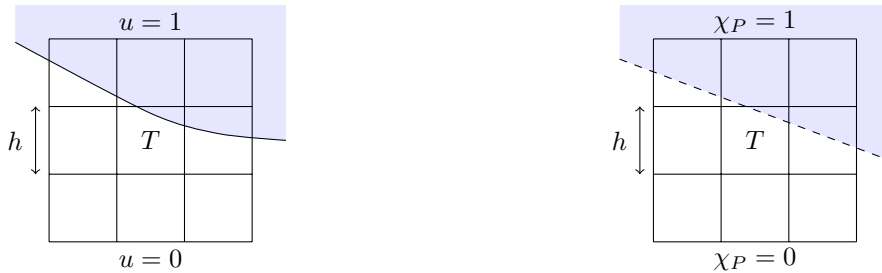
\begin{figure}[ht]
\begin{center}
\begin{minipage}{.4\textwidth}
\centering
\begin{tikzpicture}[scale=0.9]
  \draw[step=1cm, xshift=-1.5cm, yshift=-1.5cm] (0,0) grid (3,3);
  \draw [draw=none, fill=blue, fill opacity=0.1] plot [smooth] coordinates {(-2, 1.45) (0,0.4) (1,0.1) (2,-0)}--(2,2)--(-2,2)--cycle;
  \draw plot [smooth] coordinates {(-2, 1.45) (0,0.4) (1,0.1) (2,-0)};
  \node at (0,1.75) {$u=1$};
  \node at (0, -1.75) {$u=0$};
   \node at (0, 0) {$T$};
   \draw[<->] (-1.7,-0.5)--(-1.7,0.5);
   \node at (-2,0) {$h$};
\end{tikzpicture}
\end{minipage}
\hspace{1cm}
\begin{minipage}{.4\textwidth}
\centering
\begin{tikzpicture}[scale=0.9]
  \draw[step=1cm, xshift=-1.5cm, yshift=-1.5cm] (0,0) grid (3,3);
  \draw[dashed](-2,1.2)--(2,-0.3);
  \draw[draw=none, fill=blue, fill opacity=0.1](-2,1.2)--(2,-0.3)--(2,2)--(-2,2)--cycle;
  \node at (0,1.75) {$\chi_P=1$};
  \node at (0, -1.75) {$\chi_P=0$};
  \node at (0, 0) {$T$};
   \draw[<->] (-1.7,-0.5)--(-1.7,0.5);
   \node at (-2,0) {$h$};
\end{tikzpicture}
\end{minipage}
\end{center}
 \caption{Local approximation of a smooth interface by a line interface.}
 \label{strip}
\end{figure}

Intuitively the advantage of locally fitting a line or plane interface is that it has
the ability to better approximate the interface $\Gamma$ if it is smooth,
so that it is expected to improve the low order of accuracy \iref{firstorder}
of linear methods. More precisely, if $\Gamma$ has $\cC^2$ regularity, on each cell
of side-length $h$, it can be approximated with Hausdorff distance $\cO(h^2)$
by a line in 2d, a plane in 3d, etc. Therefore, if the locally reconstructed linear interface 
is optimally fitted, the $\cO(1)$ error is observed on a strip of volume $\cO(h^{2+d-1})$,
see Figure \ref{strip}.  Since the amount of singular cell is of order $\#(\cS_h)\sim h^{1-d}$, we may hope
for a reconstruction error with improved order,
\be
\|u-\t u\|_{L^p} \sim (h^{d+1} h^{1-d})^{\frac 1 p}=h^{\frac 2 p}=\corr{N}^{-\frac 2 {dp}},
\label{highorder}
\ee
that is, we double the rate compared to \iref{firstorder}. In particular, we obtain second order
accuracy $\cO(h^2)$ in the $L^1$ norm.

However, to our knowledge such estimates
have never been rigourously proved for the aforementionned method. In addition the 
approximation power of linear interface is also limited and we cannot hope to improve the
above rate for interfaces that are smoother than $\cC^2$. In this context, the objective of this paper is twofold:
\begin{enumerate}
\item introduce a theoretical framework for the
rigourous convergence analysis of local interface reconstructions from cell averages,
\item within this framework develop reconstruction methods going beyond linear interfaces
and provably achieving higher order of accuracy.
\end{enumerate}

\subsection{Outline}

In this paper, we will essentially work in the bivariate case $d=2$
which makes the exposition simpler while most of our discussion 
can be carried over to higher dimension.

The recovery methods that we study are local: on each cell $T$ identified as
singular, the unknown function $u=\Chi_\Omega$ is approximated by a simpler 
$\t u=\Chi_{\Omega_T}$ picked from a family $V_n$ that can be described by $n$ parameters 
and that enjoy certain approximation properties for interfaces having prescribed smoothness.
This approximation is computed from the cell averages 
of a rectangular stencil $S_T$ of $m\geq n$ cells centered around $T$.
We begin \S 2 by giving examples of such families and 
discussing their approximation properties for prior classes 
$\cK_s$
of $\Chi_\Omega$ associated to sets $\Omega$ with
$\cC^s$ boundaries for $s>1$. Our ultimate goal is
to develop recovery schemes such that the error 
between $u$ and $\t u$ is near optimal
in the sense of being bounded (up to a multiplicative
constant) by the 
error of best approximation of $u$ by elements 
from $V_n$. This, in particular means that 
the recovery should be exact if the true 
function $u$ belongs to $V_n$.

We introduce in \S 3 a first class of recovery strategies which we call Optimization Based Edge Reconstruction Algorithms (OBERA). Similar to LVIRA it is based on least-square fitting
of simpler interfaces such as lines, but also circles or polynomials that
allow to raise the order of accuracy.
We show that near optimality of this recovery
is ensured by an inverse stability inequality 
which can in particular
be established for line edges and $3\times 3$ stencils.
Unfortunately, this property seems more difficult to 
prove when raising the order of accuracy, which also
leads to more difficult nonlinear optimization procedures.

As a more manageable alternative, we consider 
recovery methods that are based on the identification of a certain preferred orientation 
- vertical or horizontal - for describing the interface by a function in the vicinity 
of each such cell. This leads us to the second class of recovery methods, 
termed as Algorithms for Edge Reconstruction 
using Oriented Stencils (AEROS) which is discussed in \S 4. It avoids continuous 
optimization by finding an edge described by a polynomial function
$y=p_T(x)$ or $x=p_T(y)$ after having used the previously discussed 
orientation mechanism. The polynomial $p_T$ is identified by 
simple linear equations. We show that this process satisfies 
an exactness and stability property that leads to the near
optimal recovery bound. In addition, the rate of convergence can be raised 
to an arbitrarily high order by raising the degree of $p_T$, however at 
the price of using larger stencils. The analysis of the orientation
selection mechanism, based on the Sobel filter, is 
postponed to the Appendix, where it is shown that it correctly classifies the cells
when $h$ is sufficiently small. 

All these methods are numerically tested and compared in \S 5. The differences in terms of convergence rates are confirmed, and 
we also compare the computational costs, which are by far less for 
AEROS than OBERA. We also discuss and test the 
specific recovery of corner singularities in the interface.
Finally, in the case of linear transport, we illustrate the propagation of error 
when using finite volume schemes 
based on these local recovery strategies.

An open-source python framework\footnote{\url{https://github.com/agussomacal/SubCellResolution}} 
is made available to show the methods presented here but specially to allow an easy way of creating, 
testing, and comparing new subcell resolution and interface reconstruction 
methods without the need to re-implement everything from scratch.

\section{Numerical analysis of local recovery methods}
\label{sec:analysis}

\subsection{Local approximation by nonlinear families}

The methods that we study in this paper for the recovery of 
the unknown $u= \Chi_\Omega$ are based on local approximations
of $u$ on each cell $T\in \cT_h$ that is identified as singular
by simpler characteristic functions picked from an $n$ parameter family $V_n$.

Let us give three examples that will be used further. We stress
that in all such examples $V_n$ is not an $n$-dimensional linear 
space, but instead should be thought as an $n$-dimensional nonlinear manifold.
\nl
\nl
{\bf Example 1: linear interfaces.} These are functions of the form $v=\Chi_H$ where
$H$ is a half-plane with a line interface $L=\partial H$. Such functions are described 
by $n=2$ parameters. One convenient description
is by the pair $(r,\theta)$, where 
$r\geq 0$ is the offset distance between the center $z_T$ of the cell $T$ of interest and the linear interface
and $\theta\in [0,2\pi[$ is the angle between this line and the horizontal axis.
In other words, in this case we use
$$
V_2:=\{v_{r,\theta}:=\Chi_{\{\<z-z_T,e_\theta\>\leq r\} } \: ,\; \theta\in [0,2\pi[; \; r>0\},
$$
where $e_\theta=(-\sin(\theta),\cos(\theta))$.
Of course, the $d$-dimensional generalization by half-spaces is a $d$-parameter family
where the unit normal vector $e_\theta$ lives on the $d-1$-dimensional unit sphere.
\nl
\nl
{\bf Example 2: circular interfaces.} These are functions of the form $v=\Chi_D$
or $v=\Chi_{D^c}$ where $D$ is a disc with circular interface, and $D^c$ its complement.
It is easily seen that the corresponding space $V_3$ is 
now a $3$ parameter family, and its $d$-dimensional generalization 
by characteristic function of balls and their complements is a $d+1$ parameter family. The
idea of using circles instead of lines is to increase the approximation capability.
However, as we will see, the next family turns out to be more effective both from
the point of view of theoretical analysis and computational simplicity.
\nl
\nl
{\bf Example 3: oriented graphs.} These are functions of the form $v=\Chi_P$ where 
$P$ is the subgraph or the epigraph of a function $p\in W_n$, either
applied to the coordinate $x$ or $y$, where $W_n$ is a linear space of univariate 
functions. In other words, $P$ is given by one among the four 
equations
\be
y\leq p(x), \quad y\geq p(x), \quad x\leq p(y), \quad x\geq p(y).
\label{subgraphs}
\ee
Of course, the corresponding space $V_n$ is an $n$-parameter family that is not a linear space, while $W_n$ is. Note that the linear interfaces of Example 1 are a particular case  
where $W_2$ is the set of affine functions. Raising the order of accuracy
will be achieved by taking for $W_n$ the space of polynomials of degree $n-1$.
The $d$ dimensional generalisation is obtained by taking for $W_n$ a linear space of 
functions of $d-1$ variables and considering $P$ to be defined by one of the 2d
equations
$$
x_i\leq p(x_1,\dots,x_{i-1},x_{i+1},\dots,x_d), \quad x_i\geq p(x_1,\dots,x_{i-1},x_{i+1},\dots,x_d), \quad i=1,\dots d,
$$
for some $p\in W_n$. In particular $W_n$ could be a space of multivariate polynomials. 
\nl
\nl
{\bf Example 4: piecewise linear interface.} These are functions of the form
$u=\Chi_{H_1\cap H_2}$ where $H_1$ and $H_2$ are two half planes. Therefore
the interface consists of two half lines that touch at a corner point $x_0$.
The corresponding space $V_4$ is a $4$ parameter family, for example
by considering the coordinates $x_0=(x_1,x_2)$ and the angles $\theta_1$ and $\theta_2$
of the normal vectors to the two lines. The goal of this family is to better approximate
piecewise smooth interfaces that have corner singularities.
Note that Example 1 may be viewed as a particular case where the two half-planes
coincide and there is no corner point.
\nl

Given such a family $V_n$ and a set $S\subset D$, we denote by
$$
e_n(u)_S=\min_{v\in V_n} \|u-v\|_{L^1(S)},
$$
the error of best approximation on $S$ measured in the $L^1$-norm.

\begin{remark}
Throughout this paper, we shall systematically measure error in $L^1$ norm
which is the most natural since in the case of $u=\Chi_{\Omega}$ and $\t u=\Chi_{\t \Omega}$,
this error is simply the area of the symmetric difference between domains, that is,
$$
\|u-\t u\|_{L^1(D)}=|\Omega\, \Delta\, \t \Omega|=| \Omega\cup\t \Omega - \Omega \cap \t \Omega |.
$$
Note however that estimates in $L^p$ norms can be derived 
from $L^1$ estimates in a straightforward manner since
$$
\|u-\t u\|_{L^p(D)}=|\Omega\, \Delta \, \t \Omega |^{1/p}=\|u-\t u\|_{L^1(D)}^{1/p}.
$$
\end{remark}

In our analysis of the reconstruction error on a singular cell $T$,
we will need to estimate the local approximation error on a stencil 
$S=S_T$ that consists of finitely many cells surrounding $T$.
We shall systematically consider rectangular stencils of symmetric
shape centered around $T$, so in the 2d case they are of size 
$$
m=(2k+1) \times (2l+1),
$$
for some fixed $k,l\geq 1$.

The order of magnitude $e_n(u)_S$ both depends on the type of family $V_n$
that one uses and on the smoothness property of the boundary $\Gamma=\partial \Omega$.
We describe these properties by introducing prior classes of 
characteristic functions $\Chi_\Omega$ of sets $\Omega\subset D$ 
with boundary of a prescribed H\"older smoothness. There exists
several equivalent definitions of a $\cC^s$ domain. We follow the
approach from \cite{Necas}, that expresses the fact that the boundary
can locally be described by graphs of $\cC^s$ functions (see also Chapter 4 of \cite{Adams}).

\begin{definition}
Let $s>0$. A domain $\Omega\subset \R^2$ is of class $\cC^s$ if and only if
there exists an $R>0$, $P>1$ and $M>0$,
such that for any point $z_0\in \partial\Omega$, the following holds:
there exists an orthonormal system $(e_1,e_2)$ 
and a function $\psi\in \cC^s([-R,R])$
with $\|\psi\|_{\cC^s}\leq M$, taking its value in $[-PR,PR]$
and such that
$$
z\in \Omega \iff z_2\leq \psi(z_1),
$$
for any $z=z_0+z_1e_1+z_2e_2$ with 
$|z_1|\leq R$ and $|z_2| \leq PR$.
\label{defcs1}
\end{definition}

Here, we have used the usual definition 
\[
\|\psi\|_{\cC^s}=\sup_{0\leq k\leq \lfloor s\rfloor} \|\psi^{(k)}\|_{L^\infty([-R,R])}
+\sup_{s,t\in [-R,R]} |s-t|^{ \lfloor s\rfloor-s} \Big |\psi^{(\lfloor s\rfloor)}(s)-\psi^{(\lfloor s\rfloor)}(t)\Big |,
\]
for the H\"older norm. In the case of integer smoothness, we use the convention that $\cC^s$
denotes functions with Lipschitz derivatives up to order $s-1$, so that in particular
the case $s=1$ corresponds to domains with Lipschitz boundaries. This definition
naturally extends to domains of $\R^d$ with $d>2$ with $\psi$ now being a $\cC^s$ function 
of $d-1$ variables.

We can immediately derive a first local approximation error estimate for the 
the space $V_2$ of linear interfaces from the above Example 1: let $u=\Chi_\Omega$
with $\Omega$ a domain of class $\cC^s$. Then, if $S$ is 
a $2k+1\times 2l+1$ stencil centered around a cell $T$ that is
crossed by the interface, we apply the above Definition \ref{defcs1} taking $z_0=z_T$
the center of $T$. We assume that the sidelength $h$ is small enough so that
the stencil $S$ is contained in the rectangle
$\{z=z_0+z_1e_1+z_2e_2 \, : \;, |z_1|\leq R, \,|z_2| \leq PR\}$, 
in which $\partial \Omega$ is described by the graph of the $\cC^s$ function $\psi$.
For $z\in S$, we have in addition that $|z_1|\leq C_0 h \leq R$ where $C_0$ depends only on $l$ and $k$.
Using a Taylor expansion and the smoothness of $\psi$ we find that there exists
an affine function $a$ such that
\be
\|\psi-a\|_{L^\infty([-C_0h,C_0h])} \leq C_1 h^{r}, \quad r:=\min\{s,2\},
\label{psia}
\ee
where $C_1$ depends on $C_0$ and the bound $M$ on the $\cC^s$ norm of $\psi$.
For example, we can take for $a$ the Taylor polynomials of order $1$ at $z_1=0$
when $s>1$, which corresponds to match the tangent of the interface at the point $z_0+\psi(0)e_2$,
or of order $0$ when $s\leq 1$. Then, taking $v=\Chi_H\in V_2$, where
$$
H:=\{z_2\leq a(z_1)\},
$$
is the corresponding half-space, it follows that
$$
\|u-v\|_{L^1(S)}=|S \cap (\Omega \, \Delta \, H)| \leq C h^{r+1},
$$
where $C$ depends on $(M,l,k)$. In summary, for the local approximation error of 
$\cC^s$ domains by a linear interface we have 
\be
e_n(u)_S \leq C h^{r+1}, \quad r:=\min\{s,2\}.
\label{linlocerror}
\ee
The same reasoning in $d$ dimensions delivers a local approximation 
estimate of order $h^{r+d-1}$.

One way to raise this order of accuracy for smoother domains
is to use approximation by circular interfaces from Example 2 since this allows
us to locally match the curvature in addition to the tangent. In turn we 
reach a similar estimate of order $h^{r+1}$ however with $r:=\min\{s,3\}$. 
One more systematic way of raising the order arbitrarily high is to use
approximation by oriented subgraphs from Example 3. This
approach is central to the AEROS strategies discussed in \S 5 and we
thus discuss it below in more detail.

We begin with the observation that when $s>1$, 
the unit tangent vector varies continuously on $\partial \Omega$ if $\Omega$
is of class $\cC^s$. It follows that
locally around any point $z_0$, this vector remains away either from the 
horizontal vector $(1,0)$ or from the vertical vector $(1,0)$. This allows us
to locally describe the boundary by graphs of functions of the standard cartesian
coordinates, as expressed by the following alternate definition of $\cC^s$ domains.

\begin{definition}
\label{defcart}
Let $s>1$. A domain $\Omega\subset \R^2$ is of class $\cC^s$ if and only if
there exists an $R>0$, $P>1$ and $M>0$,
such that for any point $z_0=(x_0,y_0)\in \partial\Omega$, the following holds:
there exists a function $\psi\in \cC^s([-R,R])$
with $\|\psi\|_{\cC^s}\leq M$, taking its value in $[-PR,PR]$
and such that the membership in $\Omega$ of a point $z=(x,y)$ is equivalent
to one of the two equations
\be
\label{yx}
y\leq \psi(x), \quad y\geq \psi (x), 
\ee
when $|x-x_0|\leq R$ and $|y-y_0| \leq PR$, or one of the two equations
\be
x\leq \psi(y), \quad x\geq \psi (y), 
\label{xy}
\ee
when $|y-y_0|\leq R$ and $|x-x_0| \leq PR$.
\label{defcs2}
\end{definition}

The generalization of this alternate definition to higher dimension is straighforward
by considering equations one of the form
$$
x_i \leq \psi(x_1,\dots,x_{i-1},x_{i+1},\dots, x_d) \quad {\rm or}\quad 
x_i \geq \psi(x_1,\dots,x_{i-1},x_{i+1},\dots, x_d),\quad i=1,\dots,d,
$$
with $\psi$ a $\cC^s$ function of $d-1$ variables defined on $[-R,R]^{d-1}$.

\begin{remark}
We stress that this definition is only valid for $s>1$ and not for less
smooth domains such as Lipschitz domains. For example if $\Omega$ is
a rectangle with side oriented along principal axes, then no such
local parametrization can be derived if $z_0$ is a corner point.
\end{remark}

Consider now the family $V_n$ of oriented subgraphs from Example 3,
associated with the linear space $W_n=\pP_{n-1}$ of univariate polynomials
of degree $n-1$. Let $u=\Chi_\Omega$
with $\Omega$ a domain of class $\cC^s$ for some $s>1$. Then, if $S$ is 
a $(2k+1)\times (2l+1)$ stencil centered around a cell $T$ that is
crossed by the interface, we apply the above Definition \ref{defcs2} taking $z_0=z_T$
the center of $T$. Without loss of generality, assume for example that 
the description of $\Omega$ near $z_0$ is by the equation
$$
y\leq \psi(x), 
$$
for $z=(x,y)$ in the rectangle $\{|x-x_0|\leq R, \;|y-y_0| \leq PR\}$.
We assume that the sidelength $h$ is small enough so that
the stencil $S$ is contained in this rectangle. For $z=(x,y)\in S$ we thus have 
$$
|x-x_0| \leq  C_0h \leq R, \quad C_0=k+\frac 1 2.
$$
Using Taylor formula and the smoothness of $\psi$ we find that there exists
a polynomial $p\in \pP_{n-1}$ such that
\be
\|\psi-p\|_{L^\infty([-C_0h,C_0h])} \leq C_1 h^{r}, \quad r:=\min\{s,n\},
\label{polapp}
\ee
where $C_1$ depends on $C_0$ and the bound $M$ on the $\cC^s$ norm of $\psi$.
For example, we can take for $p$ the Taylor polynomials of order $\t n=\min\{\lceil s\rceil-1, n-1\}$
at $x=x_0$. Therefore, taking $v$, where
$$
v:=\Chi_{\{y\leq p(x)\}}\in V_n,
$$
the corresponding subgraph, it follows that
$$
\|u-v\|_{L^1(S)}=2 C_0 hC_1 h^{r}=Ch^{r+1}.
$$
where $C$ depends on $(M,k)$. We treat the other cases $y\geq \psi(x)$, $x\leq \psi(y)$ and $x\geq \psi(y)$
in a similar manner. In summary, for the local approximation error of 
$\cC^s$ domains by polynomial oriented subgraphs, we have 
\be
e_n(u)_S \leq C h^{r+1}, \quad r:=\min\{s,n\}.
\label{pollocerror}
\ee
The same reasoning in $d$ dimensions delivers a local approximation 
estimate of order $h^{r+d-1}$.

\subsection{Near optimal recovery from cell averages}

The recovery methods that we develop in \S \ref{sec:obera} and \S \ref{sec:aeros} are based on
recovering on each singular cell $T$ an element $\t u_T \in V_n$ where $V_n$
is a given nonlinear family, based on the data of the cell-averages 
$$
a_S(u)=(a_{\t T}(u))_{\t T\in S},
$$
where $S=S_T$ is a rectangular stencil centered around $T$.
It can therefore be summarized by a local nonlinear recovery operator
$$
R_T: \R^m \to V_n
$$
where $m=(2k+1) \times (2l+1)$ is the size of the stencil, such that
$$
\t u_T=R_T(a_S(u)).
$$
We are interested in deriving a favorable comparison between the local recovery error
$\|u-\t u_T\|_{L^1(T)}$ and the error of best approximation by $V_n$
whose magnitude can be estimated depending on the amount of smoothness of the boundary, as 
previously discussed. 

\begin{definition}
The local recovery procedure is said to be
near-optimal over a class of function $\cK$ if there exists a fixed 
constant $C$ so that one has
\be
\|u-R_T(a_S(u))\|_{L^1(T)} \leq Ce_n(u)_S,
\label{localrecbound}
\ee
for all $u$ in this class. In particular $C$ should be independent of
the considered singular cell $T$ and mesh size $h$.
\end{definition}

\begin{remark}
In the above definition, the recovery error on $T$ is bounded by the approximation
error on the larger stencil $S$. This is due to the fact that the
recovery operator $R_T$ acts on the cell averages $a_{\t T}(u)$ for all cells $\t T\subset S$. 
\end{remark}

On regular cells $T\in \cT_h\setminus \cS_h$, that is, such that $a_T(u)=0$ or $a_T(u)=1$, we 
simply define the reconstruction by the constant value
$$
\t u_T=a_T(u),
$$
which is then the exact value of $u$. The global reconstruction of $u$ from the cell-averages 
$(a_T(u))_{T\in \cT_h}$ is given by the function 
$$
\t u=\sum_{T\in \cT_h} \t u_T \Chi_T.
$$
The global $L^1$ error can be estimated by aggregating all local error estimates, which
thus gives
\be
\|u-\t u\|_{L^1(D)} = \sum_{T\in \cS_h} \|u-\t u_T\|_{L^1(T)}
\leq  C\sum_{T\in \cS_h} e_n(u)_{S_T}.
\label{agreg}
\ee
where $C$ is the stability constant in \iref{localrecbound} and where $S_T$ denotes the
stencil centered at $T$ which is used in the reconstruction of $\t u_T$.

If $\Omega$ is a $\cC^s$ domain with $s\geq 1$, that is, at least a Lipschitz domain,
one has the cardinality estimate
\be
\#(\cS_h)\leq Ch^{-1}
\label{card}
\ee
where $C$ depends on the length of $\Gamma$. 
We may thus derive a global error estimate by combining \iref{agreg}, \iref{card}
and the local approximation estimates \iref{linlocerror} and \iref{pollocerror}:
is $\Omega$ is a $\cC^s$ domain with $s\geq 1$, we obtain
\be
\|u-\t u\|_{L^1(D)} \leq C h^{r}, 
\label{globerror}
\ee
with $r:=\min\{s,2\}$ when using local recovery by linear interfaces
and $r:=\min\{s,n\}$ when using local recovery by polynomial 
subgraphs of degree $n-1$. This estimate generalizes to the higher dimensional case,
combining local approximation estimates with the cardinality estimate
$\#(\cS_h)\leq Ch^{1-d}$.

Our central objective is now to propose local recovery methods that
provably satisfy the near optimal recovery bound \iref{localrecbound}.
For this, we start by remarking that this bound implies  
the property 
\be
R_T(a_S(v))=v, \quad \forall v\in V_n,
\label{const}
\ee
that is, the recovery is exact
for elements from $V_n$.
This property itself implies that any element from $V_n$ 
should be exactly characterized
by its cell-average on the stencil $S$. In other words, the averaging operator
$$
v\mapsto a_S(v)=(a_{\t T}(v))_{\t T\in S},
$$
should be injective from $V_n$ to $\R^m$. For this to hold, the classes $V_n$ 
from Examples 1, 2, 3, 4 need to be restricted. 

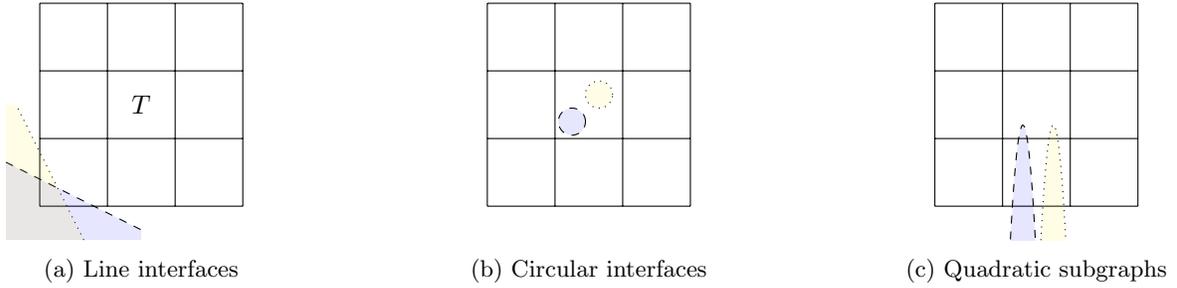
\begin{figure}
	\centering
		\begin{minipage}{.3\textwidth}
		\centering
		\begin{tikzpicture}[scale=0.9]
		  \draw[step=1cm] (0,0) grid (3,3);
		  \clip (-0.5, -0.5) rectangle (3.5, 3.5);
		  \draw[dashed](-0.5,{0.4+0.5*0.5})--(1.5,{0.4-1.5*0.5});
		  \draw[draw=none, fill=blue, fill opacity=0.1](-0.5,{0.4+0.5*0.5})--(1.5,{0.4-1.5*0.5})--(1.5, -0.5)--(-0.5,-0.5)--cycle;
		  \draw[dotted]({0.4+0.5*0.5}, -0.5)--({0.4-1.5*0.5}, 1.5);
		  \draw[draw=none, fill=yellow, fill opacity=0.1]({0.4+0.5*0.5},-0.5)--({0.4-1.5*0.5}, 1.5)--(-0.5, 1.5)--(-0.5,-0.5)--cycle;
		  \node at (1.5, 1.5) {$T$};
		\end{tikzpicture}		
	\subcaption{Line interfaces}
	\label{fig:noninjlines}
	\end{minipage}
    \hfill
		\begin{minipage}{.3\textwidth}
		\centering
		\begin{tikzpicture}[scale=0.9]
		  \draw[step=1cm] (0,0) grid (3,3);
		  \clip (-0.5, -0.5) rectangle (3.5, 3.5);
		  \filldraw[dashed, fill=blue, fill opacity=0.1] (1.25, 1.25) circle (0.2);
		  \filldraw[dotted, fill=yellow, fill opacity=0.1] (1.65, 1.65) circle (0.2);
		\end{tikzpicture}
		\subcaption{Circular interfaces}
		\label{fig:noninjcircles}
		\end{minipage}
	\hfill
		\begin{minipage}{.3\textwidth}
		\centering
		\begin{tikzpicture}[scale=0.9]
		  \draw[step=1cm] (0,0) grid (3,3);
		  \clip (-0.5, -0.5) rectangle (3.5, 3.5);
		  \filldraw[dashed, fill=blue, fill opacity=0.1] plot[smooth,domain=1:2] (\x, {-25*2*(\x-1.3)^2+1.2});
		  \filldraw[dotted, fill=yellow, fill opacity=0.1] plot[smooth,domain=1:2] (\x, {-25*2*(\x-1.75)^2+1.2});
		\end{tikzpicture}
		\subcaption{Quadratic subgraphs}
		\label{fig:noninjparab}
		\end{minipage}
\caption{Cases of non-injectivity: two elements of $V_n$ having same averages on a $3\times 3$ stencil.}
\end{figure}

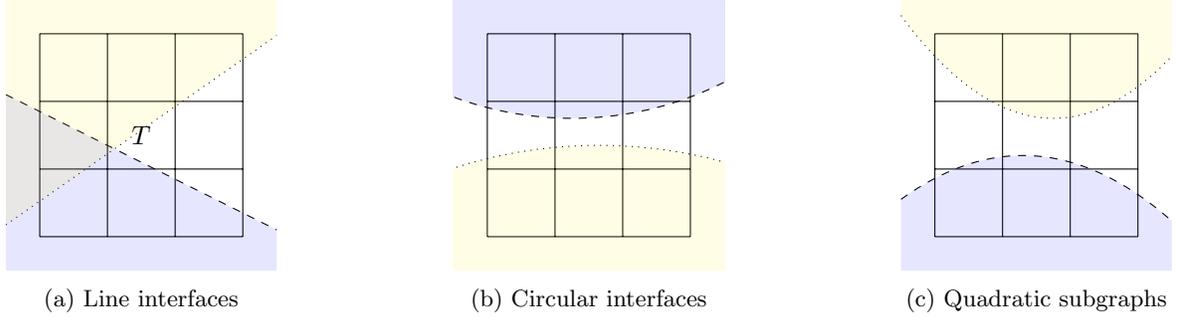
\begin{figure}
\centering
		\begin{minipage}[b]{.3\textwidth}
		\centering
		\begin{tikzpicture}[scale=0.9]
		  \draw[step=1cm] (0,0) grid (3,3);
		  \clip (-0.5, -0.5) rectangle (3.5, 3.5);
		  \draw[dashed] plot[smooth,domain=-0.5:3.5] (\x, {-0.5*(\x-1.3)+1.2});
		  \draw[draw=none, fill=blue, fill opacity=0.1](-0.5,-0.5)--plot[smooth,domain=-0.5:3.5] (\x, {-0.5*(\x-1.3)+1.2})--(3.5,-0.5)--cycle;
		  \draw[dotted] plot[smooth,domain=-0.5:3.5] (\x, {0.7*(\x-1.75)+1.75});
		  \draw[draw=none, fill=yellow, fill opacity=0.1](-0.5,3.5)--plot[smooth,domain=-0.5:3.5] (\x, {0.7*(\x-1.75)+1.75})--(3.5,3.5)--cycle;
		  \node at (1.5, 1.5) {$T$};
		\end{tikzpicture}
	\subcaption{Line interfaces}
	\label{fig:injlines}
	\end{minipage}
    \hfill
		\begin{minipage}[b]{.3\textwidth}
		\centering
		\begin{tikzpicture}[scale=0.9]
		  \draw[step=1cm] (0,0) grid (3,3);		  
		  \clip (-0.5, -0.5) rectangle (3.5, 3.5);
		  \filldraw[dashed, fill=blue, fill opacity=0.1] (1.25, 6.75) circle (5);
		  \filldraw[dotted, fill=yellow, fill opacity=0.1] (1.65, -5.65) circle (7);
		\end{tikzpicture}
		\subcaption{Circular interfaces}
		\label{fig:injcircles}
		\end{minipage}
    \hfill
		\begin{minipage}[b]{.3\textwidth}
		\centering
		\begin{tikzpicture}[scale=0.9]
		  \draw[step=1cm] (0,0) grid (3,3);
		  \clip (-0.5, -0.5) rectangle (3.5, 3.5);
		  \draw[dashed] plot[smooth,domain=-0.5:3.5] (\x, {-0.2*(\x-1.3)^2+1.2});
		  \draw[draw=none, fill=blue, fill opacity=0.1](-0.5,-0.5)--plot[smooth,domain=-0.5:3.5] (\x, {-0.2*(\x-1.3)^2+1.2})--(3.5,-0.5)--cycle;
		  \draw[dotted] plot[smooth,domain=-0.5:3.5] (\x, {0.3*(\x-1.75)^2+1.75});
		  \draw[draw=none, fill=yellow, fill opacity=0.1](-0.5,3.5)--plot[smooth,domain=-0.5:3.5] (\x, {0.3*(\x-1.75)^2+1.75})--(3.5,3.5)--cycle;
		\end{tikzpicture}
		\subcaption{Quadratic subgraphs}
		\label{fig:injparab}
		\end{minipage}
\caption{Illustration of the restrictions ensuring injectivity}
\end{figure}

In the case of linear interfaces, if $H$ and $\t H$
are two half-spaces that contain the stencil $S$, the functions $v=\Chi_H$ and $\t v=\Chi_{\t H}$
obviously have the same cell-average vector $a_S$ with component identically equal to $1$.
Asking that the linear interface passes through the stencil is thus necessary
but not sufficient as illustrated on Figure \iref{fig:noninjlines}:  two lines passing only through 
a corner cell may result in identical cell averages. The correct restriction on $V_n$ in this case
is obtained by imposing that the linear interface passes through the central cell, which in the case of 
a $3\times 3$ stencil suffices to ensure injectivity as recalled in \S 4
and illustrated on Figure \ref{fig:injlines}.
An important observation is that this type of restriction does not affect
the local error estimates \iref{linlocerror} since we are
precisely considering a stencil $S$ centered at a cell $T$ that contains the 
interface. Therefore, the tangent of the interface at the point $z_0+\psi(0)e_2$
delivering this estimate satisfies this restriction.

In the case of circular interface, asking that the disk intercepts the central cell is
not sufficient as  illustrated on Figure \iref{fig:noninjcircles}: two discs $D$ and $\t D$ of equal size
and contained in the central cell will result in identical cell averages for 
$v=\Chi_D$ and $\t v=\Chi_{\t D}$. In this case, an additional restriction should be
that the radius of the disc is sufficiently large compared to the size of the stencil, for example
by imposing that the disc center is not contained in $S$, as illustrated on Figure \ref{fig:injcircles}.

In the case of subgraphs of polynomial functions, again two subgraphs passing through
several cells of the stencil might have the same cell averages. This typically occurs
when the polynomials are too peaky, as illustrated on Figure \iref{fig:noninjparab}. 
In this case,  the additional restriction should 
be that the range of $p$ remains inside the stencil. By this we mean, say for a 
subgraph of the type $y\leq p(x)$ and a rectangular stencil $S=[a,b]\times [c,d]$,
one has that $p([a,b])\subset [c,d]$, as illustrated in Figure \iref{fig:injparab}. Similarly to linear interfaces, the local
error estimate \iref{pollocerror} is not affected by such a restriction, as we 
discuss in \S 5.

A similar type of observation shows that there is no hope to uniquely characterize a piecewise
linear interface from the cell averages on a given stencil if it is too peaky, that is, 
the opening angle of the cone embraced by the two lines 
cannot be arbitrarily close to $0$ or $2\pi$. In other words, corners cannot
be arbitraritly acute or obtuse.

\section{Reconstruction by optimization (OBERA)}
\label{sec:obera}

\subsection{Presentation of the method}

Optimization-Based Edge Reconstruction Algorithms (OBERA) consists in recovering on each
singular cell $T\in \cS_h$ a recovery $\t u_T\in V_n$ by a best fit of the available cell-average data 
on the stencil $S=S_T$. For this purpose, we solve a minimization problem of the form
$$
\t u_T =R_T(a_S(u))\in \argmin_{v\in V_n} \|a_{S}(v)-a_{S}(u)\|
$$
where $\|\cdot\|$ is a given norm on $\R^m$ where $m:=\#(S)$ is
the size of the stencil. In a practical implementation, one first simple choice is to use
 the Euclidean $\ell^2$ norm, that is, minimize the loss function
$$
\cL(u,v):=\sum_{\t T\in S} |a_{\t T}(u)-a_{\t T}(v)|^2,
$$
over all $v\in V_n$. The case of linear interface corresponds to the
LVIRA method \cite{LVIRApuckett1991volume}. Note that
$v\in V_n$ is defined through an appropriate parametrization as in Example 1, 2, 3 and 4, 
$$
\mu \in \cM \subset \R^n \mapsto v_\mu\in V_n,
$$
where $\cM$ is the restricted range of the parameter 
$\mu$ that defines the family $V_n$.
Therefore the optimization is done in practice by searching for
$$
\mu^* \in \argmin_{\mu \in  \cM} \cL(u,v_\mu)
$$
and taking $\t u_T=v_{\mu^*}$.

It is interesting to note that this recovery method is not consistent in
the sense that it does not guarantee that
\be
a_T(\t u_T)=a_T(u),
\label{consistency}
\ee
a property that is required, typically in finite volume methods since
it reflects the conservation of mass. In order to restore this property, one
possibility is to define the recovery by solving the constrained optimization problem
\be
\t u_T \in \argmin_{v\in V_n} \{\cL(u,v) \, : \, a_T(v)=a_T(u) \}
\label{constraint}
\ee
In practice, this can be emulated by modifying the loss function into
\be
\cL(u,v):=K|a_{T}(u)-a_{T}(v)|^2+\sum_{\t T\in S,\t T\neq T} |a_{\t T}(u)-a_{\t T}(v)|^2,
\label{eq:kloss}
\ee
and taking $K\gg 1$ (in our numerical tests we took $K=100$). 
As explained below, this constrained recovery satisfies similar
error bounds as its unconstrained counterpart.

The practical difficulty of the OBERA lies in the quick and accurate computation 
of the cell averages $a_{S_T}(v)$ for any given cell $T$ and $v\in V_n$, that is, have a fast evaluation
procedure for the parameter to average map
$$
\mu \in \cM\subset \R^n \mapsto a_{S_T}(v_\mu) \in \R^m,
$$
as it is needed to calculate $\cL(u,v_\mu)$ 
at each iteration of the optimization algorithm.
While analytic expressions are easily available for linear interfaces, they become more difficult
to derive for more general curves like polynomials or circle interfaces. In such cases,
the exact computation of $a_S(v_\mu)$ is possible if for each cell $\t T \in S_T$
one is able identify the points where the curved interface crosses its boundary.
We followed this approach in our numerical implementation.
Another option relies on quadrature methods as used in \cite{Despres2020}, however
at the expense of potentially many evaluations of $v$. Finally, another perspective is
to use machine learning methods in order to derive a cheaply computable surrogate of the
parameter to average map.

Even with such tools in hand, the computation of $a_S(v_\mu)$ for the many parameter values $\mu$ that
are explored through the optimization process
results in a time overhead that one would like to avoid. For linear interfaces, this
was achieved by the ELVIRA method, as it only requires $6$ calculations of $a_S(v_\mu)$ 
in order to decide which $\mu^*$ should be retained \corr{(see Figure \ref{fig:elviracalculations})}. This can also be avoided
for more general interfaces
having higher order geometric approximations by the AEROS approach
that we present in \S \ref{sec:aeros}.

\subsection{Analysis of the recovery error}

In order to prove that the recovery error is near optimal in the
$L^1(S)$ norm, we follow a general strategy introduced in \cite{NonLinearReduced}
which is based on comparing the continuous $L^1(S)$ norm of functions
and the discrete $\ell^1(\R^m)$ norm of their cell-averages. 

In the 2d case, one obviously has on the one hand the inequality
\be
h^2\|a_S(v)\|_{\ell^1}\leq  \|v\|_{L^1(S)}, \quad v\in L^1(S), 
\label{stab}
\ee
which is obtained
by summing up 
$$
h^2|a_{\t T}(v)|=\Big |\int_{\t T} v \Big | \leq \|v\|_{L^1(\t T)},
$$
over all $\t T\in S$. This property reflects the stability of the averaging operator,
between the continuous and the conveniently normalized discrete norm. 
Note that in more general dimension $d$ the normalizing factor is $h^d$.

Conversely, we say that the family $V_n$ satisfies an {\it inverse stability property} if
there exists a constant $C$ independent of $h$ such that
\be
  \|v-\t v\|_{L^1(S)}\leq Ch^2\|a_S(v)-a_S(\t v)\|_{\ell^1}, \quad \forall \, v,\t v\in V_n.
 \label{invstab}
\ee
 We stress that such a property cannot hold for general pairs of integrable functions,
 their membership in  $V_n$ is critical. Note that this property is a more
 quantitative version of the injectivity of the map $v\mapsto a_S(v)$ from
 $V_n$ to $\R^m$. Its validity is thus conditioned to a proper restriction of the
 classes $V_n$ from the various Examples 1, 2, 3, and 4, as already explained in \S 2.2.
 In the particular case of linear edges the following result was proved in \cite{NonLinearReduced}.
 
 \begin{theorem}
 \label{theoinvlin}
 Let  $S$ be the $3\times 3$ stencil centered at $T$
 and let $V_2$ be the family of linear interfaces from Example 1, with the
 restriction that the linear interfaces passes through $T$. Then
 \iref{invstab} holds and the best constant is $C=\frac 3 2$.
 \end{theorem}
 
The above stability and inverse stability property allow us to assess
the recovery error in the following way. We first write that for any $v\in V_n$
$$
\|u-\t u_T\|_{L^1(T)}
\leq \|u-v\|_{L^1(T)}+ \|v-\t u_T\|_{L^1(T)}
\leq   \|u-v\|_{L^1(T)}+ Ch^2\|a_S(v)-a_S(\t u_T)\|_{\ell^1},
$$
where we have used \iref{invstab}. We then have
$$
\|a_S(v)-a_S(\t u_T)\|_{\ell^1}
\leq \sqrt m \|a_S(v)-a_S(\t u_T)\|_{\ell^2}
\leq  2\sqrt m \|a_S(v)-a_S(u)\|_{\ell^2},
$$
where the first inequality is Cauchy-Schwartz and the second 
comes by triangle inequality and the $\ell^2$ minimization property of $\t u_T$.
Finally, we have
$$
 \|a_S(v)-a_S(u)\|_{\ell^2},
\leq  \|a_S(v)-a_S(u)\|_{\ell^1}\leq h^{-1}\|u-v\|_{L^1(S)},
$$
by using \iref{stab}. Combining all these,  and using that $v\in V_n$ is arbitrary, we obtain that
$$
\|u-\t u_T\|_{L^1(T)}\leq (1+2C\sqrt m) e_n(u)_S,
$$
which is summarized in the following.

\begin{theorem}
Under \iref{stab} and \iref{invstab}, the recovery by $\ell^2$ minimization is near optimal 
in $L^1$ norm with multiplicative constant $1+2C\sqrt m$. In the case of linear interfaces using $3\times 3$ stencils (LVIRA), this constant is $1+2\frac 3 2 \sqrt 9 =10$. 
\end{theorem}

As observed in \S 2.2, near optimal local recovery allows us to 
derive convergence rates for smooth domains in terms of the global 
error estimate \iref{globerror}. This gives the following.

\begin{corollary}
If $\Omega$ is a $\cC^s$ domain with $s\geq 1$, 
the LVIRA method which is OBERA recovery based on linear interfaces converges
in $L^1$ norm at rate $\cO(h^r)$ with $r=\min\{s,2\}$.
\label{col:lvira-convergence}
\end{corollary}

\begin{remark}
Note that if we were using a more general $\ell^p$ norm for the data fitting, we would
obtain a similar result with constant $1+C2 m^{1-\frac 1 p}$. The fact that we do not need
to restrict ourselves just to $p=2$ had been mentioned in \cite{Pilliod2004} concerning
the LVIRA method.
\end{remark}

As previously remarked, we can modify the OBERA approach in order to impose the consistency condition \iref{consistency}, by solving the constrained optimization problem \iref{constraint}, that is optimizing inside the subfamily
$$
\t V_n:=\{v \in V_n \; : \; a_T(v)=a_T(u)\} \subset V_n
$$
It is obvious that if the inverse stability property \iref{invstab}
is valid for $V_n$, it is also valid for the smaller set $\t V_n$.
We thus reach a similar estimate
$$
\|u-\t u_T\|_{L^1(T)}\leq (1+2C\sqrt m) \t e_n(u)_S,
$$
where $\t e_n(u)_S$ is the error of best approximation of $u$ in 
$L^1(S)$ norm from the element of $\t V_n$, that is best approximation
from $V_n$ under the consistency constraint. 

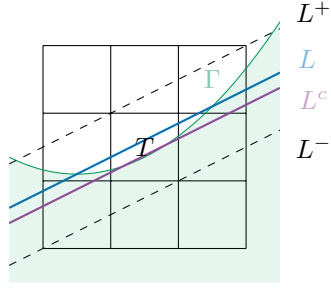
\begin{figure}
\centering
	\begin{tikzpicture}[scale=0.9]
	  \draw[step=1cm] (0,0) grid (3,3);

	  \node[ForestGreen, opacity=0.5] at (2.5, 2.5) {$\Gamma$};
	  \node at (4, 1.5) {$L^-$};
	  \node[NavyBlue, opacity=0.5] at (3.9, 2.8) {$L$};
	  \node[Purple, opacity=0.5] at (4, 2.2) {$L^c$};
	  \node at (4, 3.5) {$L^+$};

	  \clip (-0.5, -0.5) rectangle (3.5, 3.5);

	  \draw[solid, ForestGreen] plot[smooth,domain=-0.5:3.5] (\x, {0.5*(\x-1.5)+1.35+0.25*((\x-1.5)^2)});
	  \draw[dashed] plot[smooth,domain=-0.5:3.5] (\x, {0.5*(\x-1)+2});
	  \draw[solid, Purple, thick] plot[smooth,domain=-0.5:3.5] (\x, {0.5*(\x-1.5)+1.375});
	  \draw[solid, NavyBlue, thick] plot[smooth,domain=-0.5:3.5] (\x, {0.5*(\x-1.5)+1.6});
	  \draw[dashed] plot[smooth,domain=-0.5:3.5] (\x, {0.5*(\x-2)+1});

	  \draw[draw=none, fill=ForestGreen, fill opacity=0.1](-0.5,-0.5)--plot[smooth,domain=-0.5:3.5] (\x, {0.5*(\x-1.5)+1.35+0.25*((\x-1.5)^2)})--(3.5,-0.5)--cycle;
	  \node at (1.5, 1.5) {$T$};
	\end{tikzpicture}
\caption{By shifting the linear interface $L$ one achieves
average consistency on $T$ by a linear interface $L^c$ having the
same order of accuracy.}
\label{fig:ccconsistency}
\end{figure}

We thus need to 
understand if $\t e_n(u)_S$ satisfies similar size estimates as
$e_n(u)_S$.  While this cannot be ensured in full generality (for example 
the set $\t V_n$ could be empty), the following simple argument shows that this
indeed holds in the particular case of linear interface: the estimate 
\iref{psia} shows that the function $\psi$ describing the interface 
in the stencil $S$ satisfies
$$
a^-\leq \psi \leq a^+,
$$
where $a^-=a-C_1h^r$ and $a^+=a+C_1h^r$ are two affine functions that
parametrize two linear interfaces $L^-$ and $L^+$ that circumscribe $\Gamma$
in $S$, as illustrated on Figure \ref{fig:ccconsistency}. For the corresponding
halfplanes $H^-$ and $H^+$ thus satisfy
$$
a_T(\Chi_{H^-}) \leq a_T(u) \quad{\rm and} \quad a_T(\Chi_{H^+}) \geq a_T(u).
$$
Therefore, by sliding continuously a linear interface
between $L^-$ and $L^+$, which corresponds to 
the affine function $a+t$ when $t$ varies in $[-C_1h^r,C_1h^r]$, 
there exists an intermediate interface $L^c$ corresponding
to a particular $t^c$ and halfplane $H^c$ for which one has
$$
a_T(\Chi_{H^c})=a_T(u).
$$
Therefore $v=\Chi_{H^c}\in \t V_n$, and since one also has
\be
\|\psi-a^c\|_{L^\infty([-C_0h,C_0h])} \leq C_1 h^{r}, \quad r:=\min\{s,2\},
\label{psiac}
\ee
we reach the same estimate for $\t e_n(u)_S$ as the one obtained for $e_n(u)_S$.

From the theoretical perspective, an open problem is to establish
the inverse stability bound \iref{invstab} for nonlinear families $V_n$ offering higher order
approximation properties than the linear interface, for which the proof of Theorem \ref{theoinvlin}
is already quite involved. This, together with the already mentioned computational 
complexity of the optimization process, leads us to give up on OBERA for higher order
geometrical approximation in favor of the AEROS approach that we next discuss.

\section{Reconstruction on oriented stencils (AEROS)}
\label{sec:aeros}

\subsection{Presentation of the method}
\label{subsec:aeros}

Algorithms for Edge Reconstruction using Oriented Stencils (AEROS) are based 
on recovering an element from the family $V_n$ presented in Example 3. We thus
recover on each singular cell $T\in \cS_h$ a domain having one of the four subgraph or epigraph forms 
\iref{subgraphs}, that is, an interface having one of the two Cartesian forms
$y=p(x)$ or $x=p(y)$, where $p\in W_n$ a linear space of dimension $n$. 

More specifically we consider
for some fixed $k\geq 1$ the space
$$
W_{2k+1}=\pP_{2k},
$$
of polynomials of even degree $2k$. We then use stencils $S_T$
containing $T$ of the form $(2k+1)\times L$
or $L\times (2k+1)$, when reconstructing an interface of the form $y=p(x)$
or $x=p(y)$ respectively, for some $L>0$.
As we further explain, the value of $L$ and the exact positioning 
of $T$ inside $S_T$ may depend on the considered cell $T$.

As a first step we need to identify for each $T\in \cS_h$ the exact 
orientation of the subgraph or epigraph
that we decide to use. The decision must be based on the available data of the cell averages .
We have already noticed that if $\Omega$ is a $\cC^s$ 
domain for $s>1$, then it can itself be locally described by (at least) one of the four forms
with a function $\psi\in \cC^s$ describing the interface, as expressed by \iref{yx} and \iref{xy} in Definition \ref{defcart}.
Our objective is that our choice of form in the recovery is consistent with the form of the 
exact interface over the stencil $S_T$ for each cell $T$. 

\begin{figure}
    \centering    \includegraphics[width=0.45\linewidth]{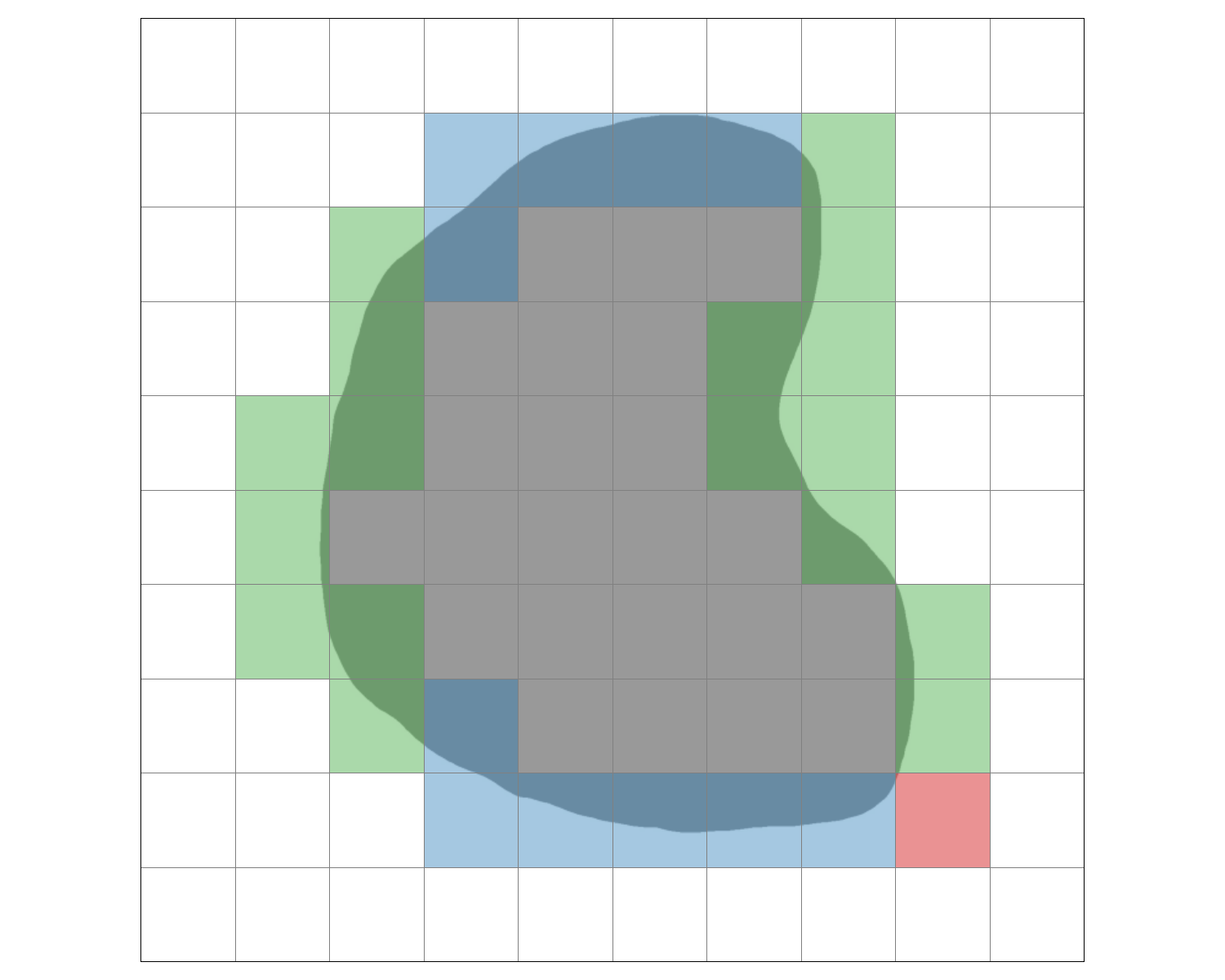}
    \includegraphics[width=0.45\linewidth]{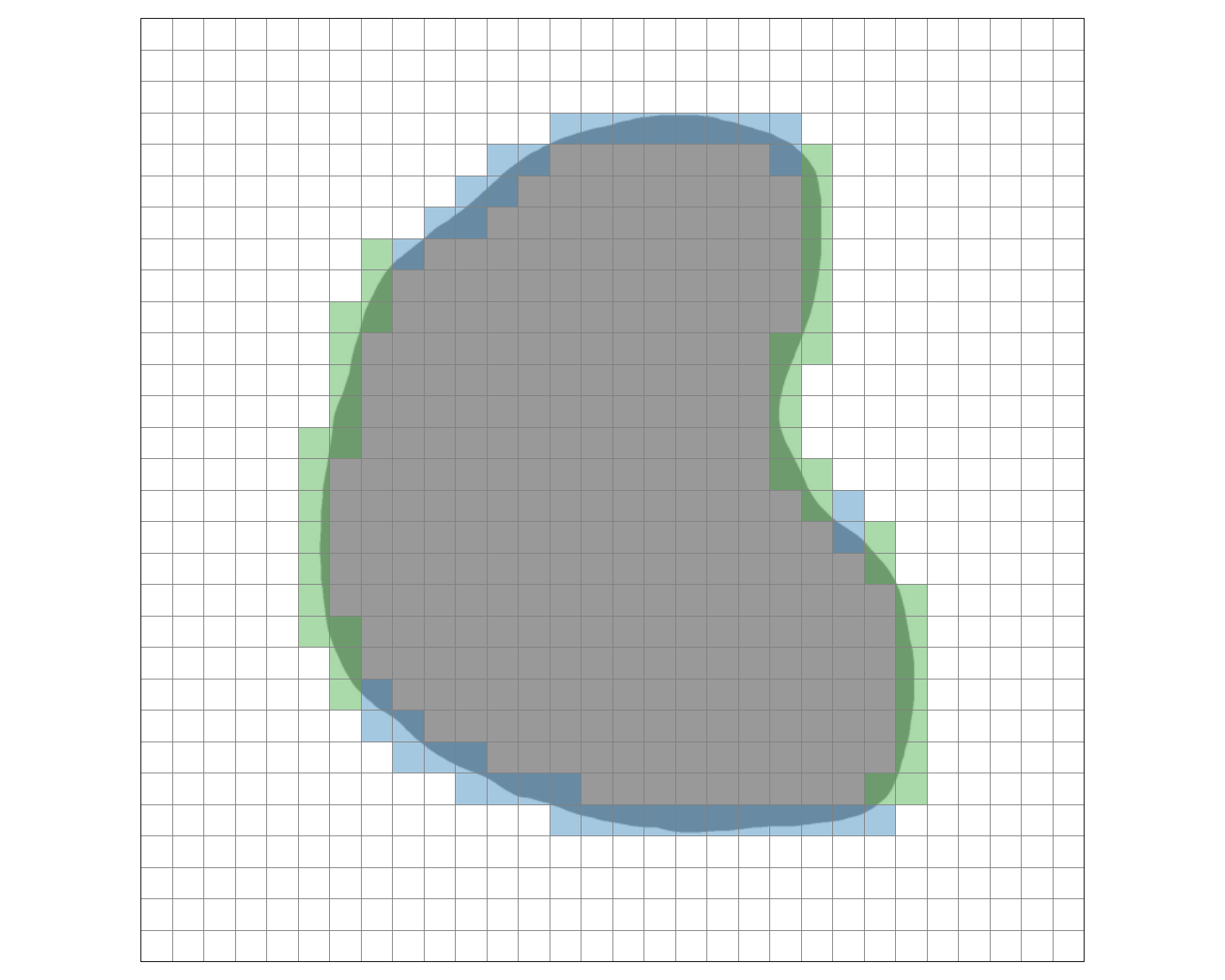}
    \caption{On the left $h=1/10$ and on the right $h=1/30$. The singular cells
idenfied as horizontally and vertically oriented are pictured
in blue and green. For $h=1/10$, there exists a singular cell $T$ (indicated in red)
for which the interface cannot be described by a graph $y=\psi(x)$ or $x=\psi(y)$ 
in any stencil $S_T$ of width $3$ that contains $T$. This is no more the case for $h=1/30$
when using an adaptive selection of the stencil.}
\label{fig:criticalscale}
\end{figure}

We thus need to identify for each $T$ an orientation $y=\psi(x)$ or $x=\psi(y)$ for the
exact interface over the stencil $S_T$. Let us immediately observe that this is only
possible if $h$ is below a certain critical resolution $h^*=h^*(\Omega)$ that depends
on the amount of variation of the tangent to the interface, as shown on Figure \ref{fig:criticalscale}
for the case $k=1$ (stencils of widths $3$).

More precisely, when $\Omega$ is a $\cC^s$ domain with $s> 1$, the variation
of the slope of the tangent to the interface between two points $z$ and $z'$ is 
controlled by a bound of the form $M|z-z'|^r$ where $r:=\min\{1,s-1\}$ 
and $M$ the bound on the $\cC^s$ norm of the functions $\psi$ that describe
the interface. Therefore a given orientation, say $y\leq \psi(x)$, can be maintained on a stencil $S_T$
of width $2k+1$ in the $x$ direction provided that $h\leq h^*\sim k^{-1}M^{-1/r}$.

\begin{figure}[!ht]
\begin{center}
\includegraphics[width=\textwidth]{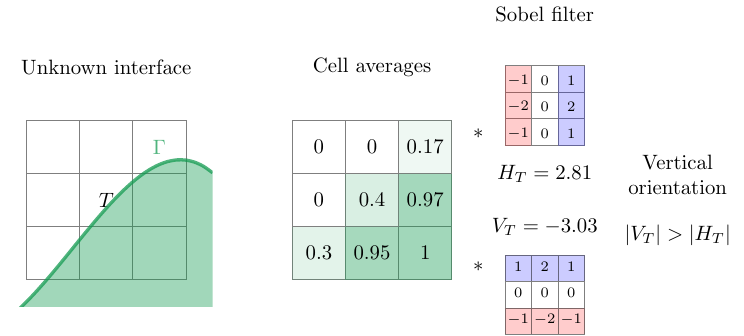}
\includegraphics[width=\textwidth]{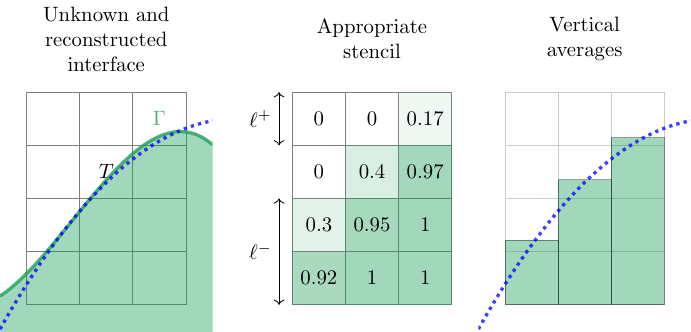}
\caption{\corr{Steps for computing AEROS method. First the preferable orientation is extracted from a $3\times 3$ stencil centered around $T$ by applying the Sobel filter (a convolution with the horizontal and vertical Sobel kernels). Then an appropriate stencil is found so that the interface crosses the sides of the stencil. Finally, the column averages are calculated and a polynomial is used to approximate the interface.}}
\label{fig:aerossteps}
\end{center}
\end{figure}

For identifying the orientation, we introduce a selection mechanism
based on a numerical gradient computed by the Sobel filter. We denote by
$T_\e$ with $\e=(\e_x,\e_y)\in \{-1,0,1\}^2$ the cells in the $3\times 3$ stencil centered
around $T=T_{0,0}$ where $\e_x$ and $\e_y$ indicate the amount of displacement by $h$
from $T$ in the $x$ and $y$ direction, respectively. We then define the numerical gradient
$$
G_T=(H_T,V_T),
$$
with horizontal component
$$
H_T:=2a_{T_{1,0}}+a_{T_{1,1}}+a_{T_{1,-1}}-(2a_{T_{-1,0}}+a_{T_{-1,1}}+a_{T_{-1,-1}})
$$
\corr{obtained by convolution between the cell averages $a_{T_\e}=a_{T_\e}(u)$ and the horizontal Sobel kernel (see Figure \ref{fig:aerossteps}). Similarly, the vertical component is defined as
$$
V_T:=2a_{T_{0,1}}+a_{T_{1,1}}+a_{T_{-1,1}}-(2a_{T_{0,-1}}+a_{T_{1,-1}}+a_{T_{-1,-1}}).
$$}

The selection mechanism is based on comparing the absolute size of $H_T$ and $V_T$
and examining their sign. More precisely:
\begin{enumerate}
\item
If $|V_T|\geq |H_T|$ and if $V_T\leq 0$, we search for a subgraph of the form $y\leq p(x)$.
\item
If $|V_T|\geq |H_T|$ and if $V_T> 0$, we search for an epigraph of the form $y\geq p(x)$.
\item
If $|V_T|<|H_T|$ and if $H_T\leq 0$, we search for a subgraph of the form $x\leq p(y)$.
\item
If $|V_T|< |H_T|$ and if $H_T> 0$, we search for an epigraph of the form $x\geq p(y)$.
\end{enumerate}

One important result is that this selection mechanism correctly detects the orientation
of the exact interface for $h$ sufficiently small, as also illustrated on Figure
\ref{fig:criticalscale}.

\begin{theorem}
\label{theosobel}
Let $\Omega$ be a $\cC^s$ domain for some $s>1$, then there exists $h^*=h^*(\Omega)$ such that
under the assumption $h<h^*$, the following holds for any $T\in \cS_h$: in each of the 
above cases $(1,2,3,4)$ of the selection mechanism, the exact domain $\Omega$ can be described by an equation of the same form with $p$ replaced by
a function $\psi\in \cC^s$ over the stencil $S_T$ centered at $T$
and of size $(2k+1)\times (2l+1)$
in case $(1,2)$ or $(2l+1)\times (2k+1)$ in case $(2,3)$ with $l=k+2$. The graph of $\psi$ remains
confined in $S_T$ in the following sense: denoting by
$I\times J:=\bigcup\{\t T \; : \; \t T\in S_T\}$ 
the total support of $S_T$, one has
$\psi(I)\subset J$ in case $(1,2)$ and $\psi(J)\subset I$ in case $(3,4)$.
\end{theorem}

We postpone the proof of this result to the Appendix, and proceed with the description
and error analysis of the recovery method. We place ourselves in case 1 without loss of generality since all other cases are dealt with similarly up to an obvious exchange of $x$ and $y$
or change of sign in one of these variable.

According to the above theorem,
we are ensured that $\Omega$ is characterized by the equation
$y\leq \psi(x)$ when $(x,y)\in S_T$ where $S_T$ is a stencil of size
$(2k+1)\times (2l+1)$ with $l=k+2$. The choice $l=k+2$ is conservative
and our numerical experiments show that it can sometimes
be reduced while maintaining
the property that the graph of $\psi$ remains
confined in $S_T$. In practice we use the following adaptive
strategy to use a stencil of minimal vertical side. 

By convention, we denote by $(i,j)$ the coordinates
of a generic cell $\t T$ when the lower left corner of $\t T$ is $(ih,jh)$.
Let $(i_T,j_T)$ be the coordinates of the singular cell $T$. We explore
the neighboring cells by defining 
$$
l^-:=\min\{l>0 \; : \; a_{\t T}(u)=1, \; i=i_T-k,\dots,i_T+k, \; j=j_T-l-1\},
$$
which is the smallest lower shift below which we find a row of non-singular cells.
Likewise, we define
$$
l^+:=\min\{l>0 \; : \; a_{\t T}(u)=0, \; i=i_T-k,\dots,i_T+k, \; j=j_T+l+1\}.
$$
Then, we take for $S_T$ the stencil that consists of cells $\t T$ of coordinates $(i,j)$ for 
$i=i_T-k,\dots,i_T+k$ and $j=j_T-l^{-},\dots,j_T+l^+$. This stencil has size
$(2k+1) \times (1+l^-+l^+)$ and is centered around the cell $T$
horizontally but not vertically. From the definition of $l^-$ and $l^+$
we have the guarantee that the graph of $\psi$ remains
confined in $S_T$ \corr{(see Figure \ref{fig:aerossteps})}. One option to further adapt the stencil $S_T$
is to allow that it is also not centered horizontally but still contains $T$.
This leads to $(2k+1) \times (1+l^-+l^+)$ stencils corresponding to values
$i=i_T-k^-,\dots,i_T+k^+$ and $j=j_T-l^{-},\dots,j_T+l^+$, where $k^-,k^+\geq 0$ are
such that $k^-+k^+=2k$ and are selected so to minimize the vertical size $1+l^-+l^+$.

Once the stencil $S_T$ has been selected, the polynomial $p_T\in \pP_{2k}$ 
is constructed as follows. For each $i=\corr{i_T-k^-},\dots,i_T+k^+$, we denote by
$R_{i}$ the column that consists of the cells $\t T\in S_T$ with first coordinate equal to $i$
and define the corresponding column average
$$
a_{i}(u)=h\left( \sum_{\t T\in R_{i}} a_{\t T}(u) \right).
$$
Since the graph of $\psi$ remains
confined in $S_T$, it follows that $a_{i}(u)$ can be identified
to the univariate cell average of $\psi$ on the interval $[ih,(i+1)h]$ after having subtracted the base
elevation of the stencil \corr{(see Figure \ref{fig:aerossteps})}, that is, 
$$
a_{i}(u)+b_T=\frac 1 h\int_{ih}^{(i+1)h} \psi(x)dx, \quad i=\corr{i_T-k^-},\dots,i_T+k^+, \quad b_T:=(j_T-l^{-})h.
$$
We then define $p_T\in \pP_{2k}$ as the unique polynomial of degree at most $2k$
that agrees with the observed averages of $\psi$, that is, such that
$$
\frac 1 h\int_{ih}^{(i+1)h} p_T(x)dx=a_{i}(u)+b_T \quad i=\corr{i_T-k^-},\dots,i_T+k^+.
$$
The polynomial $p_T$ is sometimes called the interpolant of averages, and its existence
and uniqueness is standard, similar to the more usual Lagrange interpolant of point values.
In particular, it is easily checked that the $\pP_{2k}$ interpolant of the averages of a function
$v$ on $2k+1$ adjacent intervals is the derivative of the $\pP_{2k+1}$ Lagrange
interpolant at the $2k+2$ interval endpoints for the primitive of $v$.

Quite remarkably, although we are locally approximating $u$ by a nonlinear family,
we observe that the recovery map
$$
L_T: \psi \mapsto p_T,
$$
is linear, and that the AEROS recovery approach amounts to solving a simple $(2k+1)\times (2k+1)$ linear system, resulting in substantial computational saving compared to the OBERA approach.

\subsection{Analysis of the recovery error}

In order to assess the recovery error, 
we first observe that the above described strategy has the property of exact recovery for polynomials
\be
\psi \in \pP_{2k} \implies p_T=\psi,
\ee
due to the uniqueness of the interpolant of averages. In other words,
AEROS recovers on $T$ the true interface if it is described by a polynomial
of degree $2k$ on $S_T$.

Our next observation is that the linear application
$L_T$ is stable in the max norm over 
the relevant interval $I_T=[(i_T-k^-)h,\dots,(i_T+k^+)h]$,
with stability constant that does not depend on $h$. This can be
proved by making the affine change of variable 
$$
x=\vp(\hat x)=h(i_T+\hat x), 
$$
that maps the reference interval $\hat I:=[-k^-,\dots,k^+]$ onto $I_T$.
Then, it is readily checked that
$$
L \psi \circ \vp=\hat L (\psi\circ \vp),
$$
where $\hat L$ is the average interpolant for the intervals of size $1$
contained in $\hat I$. Therefore, we may
write 
$$
\|L \psi \|_{L^\infty(I_T)}
=\|L \psi\circ \vp \|_{L^\infty(\hat I)} 
=\|\hat L (\psi\circ \vp) \|_{L^\infty(\hat I)} 
\leq C \|\psi\circ \vp \|_{L^\infty(\hat I)} 
=\hat C \|\psi \|_{L^\infty(I_T)} 
$$
where the constant $\hat C$ is the norm of $\hat L$ acting on $L^\infty(\hat I)$. This
constant only depends on $k$.

We are now in position to obtain an error estimate by writing
for all $p\in \pP_{2k}$,
$$
\|\psi-p_T\|_{L^\infty(I_T)} \leq \|\psi-p\|_{L^\infty(I_T)} +\|p_T-p\|_{L^\infty(I_T)} 
=\|\psi-p\|_{L^\infty(I_T)} +\|L(\psi-p)\|_{L^\infty(I_T)} \leq (1+\hat C)\|\psi-p\|_{L^\infty(I_T)},
$$
where we have combined exact recovery of polynomials and uniform stability. Since $p$ is arbitrary
we have obtained the following result.

\begin{theorem}
The AEROS recovery of the interface based on 
polynomials of degree $2k$ satisfies for each singular cell the near optimality property
\be
\|\psi-p_T\|_{L^\infty(I_T)} \leq (1+\hat C)\min_{p\in \pP_{2k}} \|\psi-p\|_{L^\infty(I_T)}.
\ee
\end{theorem}

This error bound of $\psi$ in $L^\infty$
readily induces an $L^1(T)$ error bound between the recovery
$$
R_T(a_S(u))=u_T:=\Chi_{\{y\leq p_T(x)\}}
$$
and $u$ by multiplying by the width of the cell. This gives
\be
\|u-R_T(a_S(u))\|_{L^1(T)}\leq (1+\hat C)h\min_{p\in \pP_{2k}} \|\psi-p\|_{L^\infty(I_T)}.
\label{localgraphbound}
\ee
Note that the right-side is not $e_n(u)_S$ and we thus have not 
obtained the near-optimal recovery property in the form \iref{localrecbound}.
Nevertheless we can derive from \iref{localgraphbound} the same convergence rates
since these are based on the Taylor polynomial approximation error \iref{polapp},
which thus yields the following result.

\begin{theorem}
Let $\Omega$ be a $\cC^s$ domain for some $s\geq 1$. The AEROS recovery of the interface based on 
polynomial of degree $2k$ satisfies for each singular cell $T$ a local error bound of the form
\be
\|u-R_T(a_S(u))\|_{L^1(T)}\leq Ch^{r+1},\quad  r:=\min\{s,2k+1\},
\label{localorderbound}
\ee
and the global error bound \iref{globerror} of order $\cO(h^r)$ for the same value of $r$.
\label{th:aeros-convergence}
\end{theorem}

\section{Numerical experiments}

\subsection{Recovery of smooth domains}

\begin{figure}
     \begin{subfigure}[b]{0.32\textwidth}
         \centering
         \includegraphics[width=\linewidth]{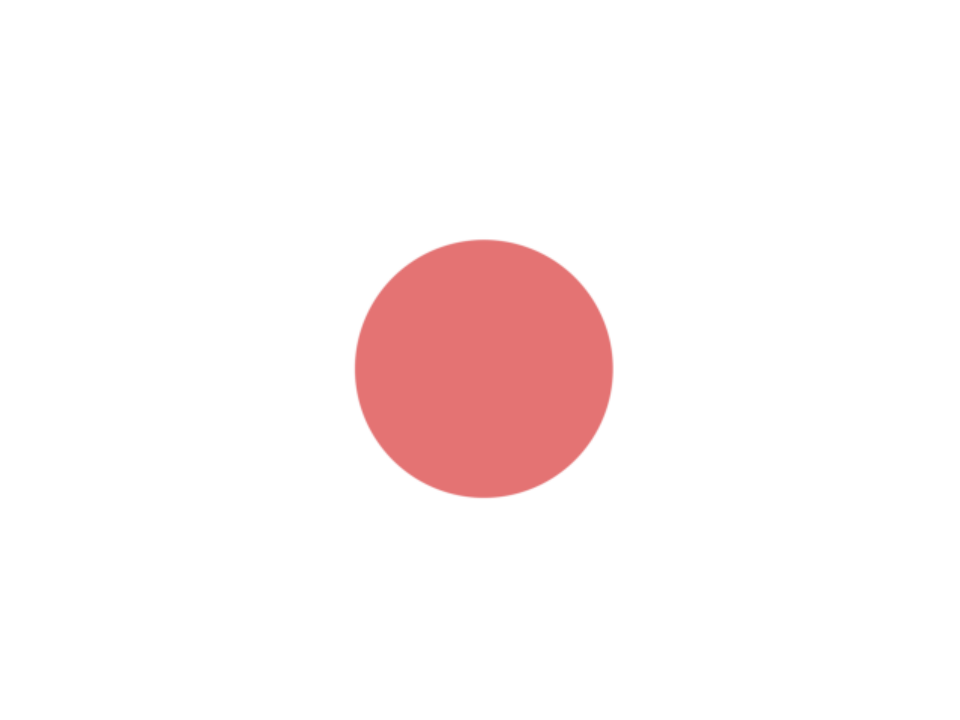}
         \caption{Circle.}
         \label{fig:circle}
     \end{subfigure}
     \hfill
     \begin{subfigure}[b]{0.32\textwidth}
         \centering
         \includegraphics[width=\linewidth]{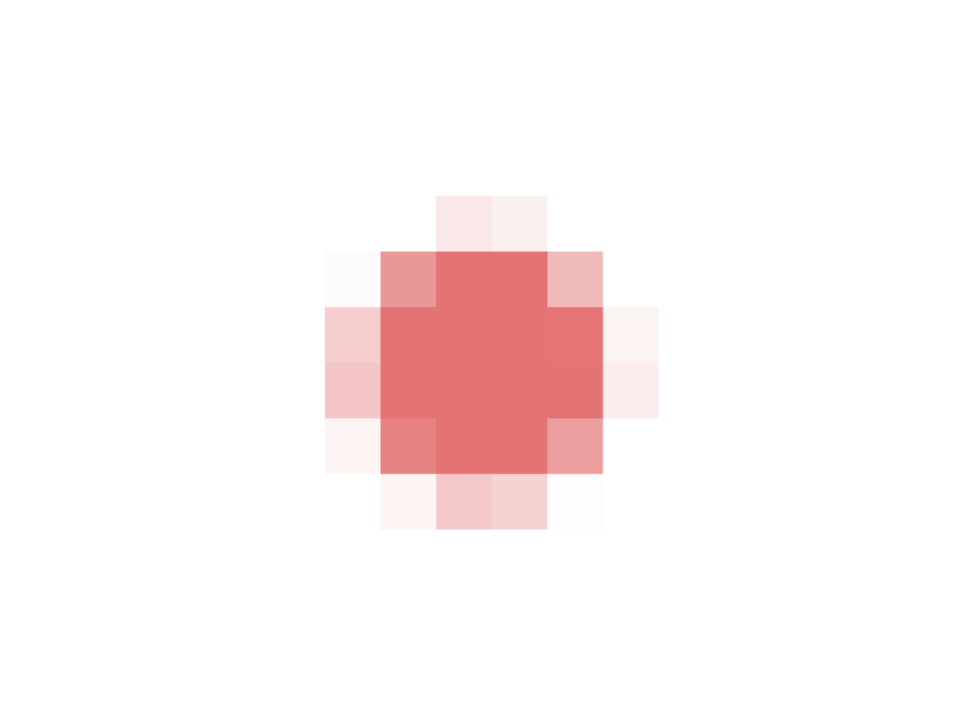}
         \caption{$h=1/10$.}
         \label{fig:circleavg1}
     \end{subfigure}
     \hfill
     \begin{subfigure}[b]{0.32\textwidth}
         \centering
         \includegraphics[width=\linewidth]{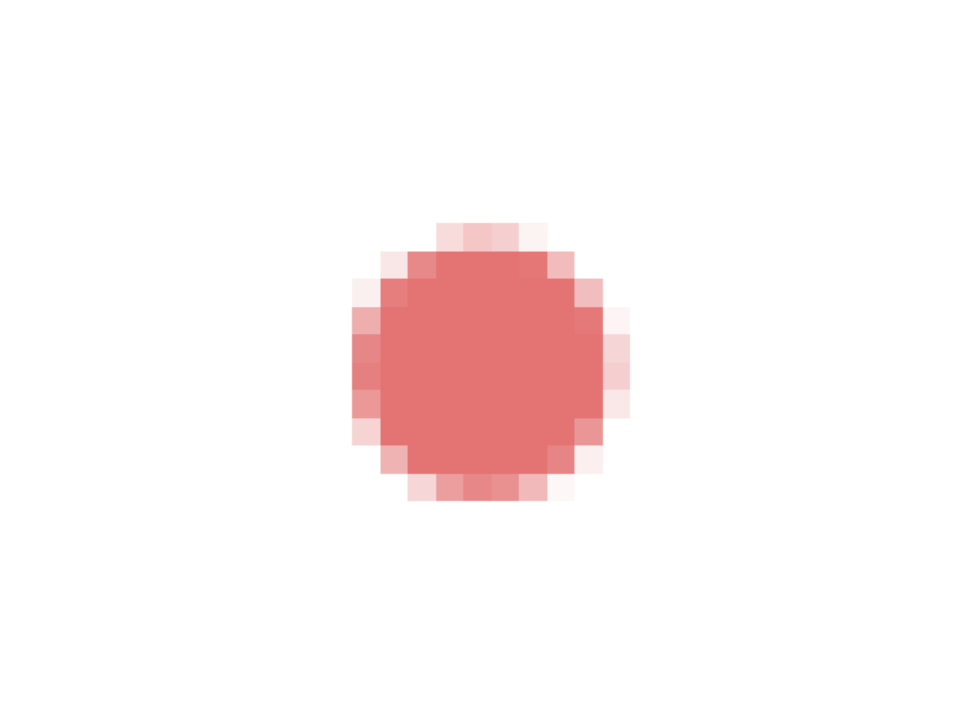}
         \caption{$h=1/20$.}
         \label{fig:circleavg2}
     \end{subfigure}
    \caption{Circular domain and its cell average data at different scales of refinement.}
    \label{fig:circdomain}
\end{figure}

In this section, we compare various recovery strategies in terms of:
\begin{enumerate}
\item
visual aspects,
\item
quantitative rate of convergence,
\item
computational time.
\end{enumerate}
In order to draw the second comparison, we first consider
the simple case of a domain $\Omega$ with circular boundary (see Figure \ref{fig:circdomain})
for which the recovery error can be computed within machine precision.

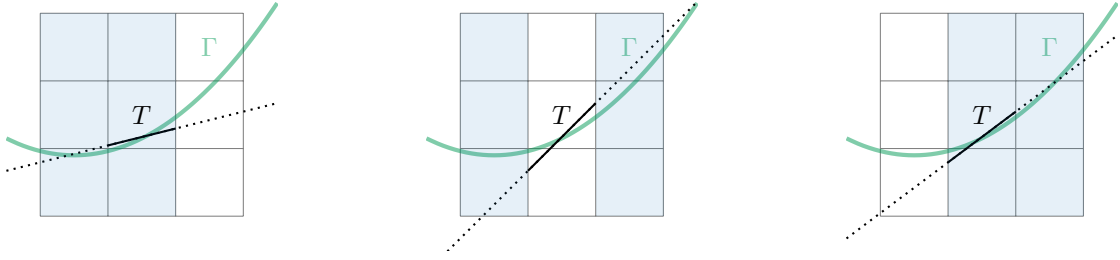
\begin{figure}[ht]
\begin{center}
\begin{minipage}{.3\textwidth}
\centering
\begin{tikzpicture}[scale=0.9]
	  \draw[step=1cm, opacity=0.5] (0,0) grid (3,3);

	  \node[ForestGreen, opacity=0.5] at (2.5, 2.5) {$\Gamma$};

	  \clip (-0.5, -0.5) rectangle (3.5, 3.5);

	  \draw[solid, ultra thick, ForestGreen, opacity=0.5] plot[smooth,domain=-0.5:3.5] (\x, {0.5*(\x-1.5)+1.35+0.25*((\x-1.5)^2)-0.2});
	  
	  \draw[dotted, thick] plot[smooth,domain=-0.5:3.5] (\x, {1.3708333+0.25*(\x-1.5)-0.2});
	  \draw[solid, thick] plot[smooth,domain=1:2] (\x, {1.3708333+0.25*(\x-1.5)-0.2});
	  
      \fill[fill=NavyBlue, fill opacity=0.1] (0, 0) rectangle (2, 3);

	  \node at (1.5, 1.5) {$T$};
	\end{tikzpicture}
\end{minipage}
\hspace{0.25cm}
\begin{minipage}{.3\textwidth}
\centering
\begin{tikzpicture}[scale=0.9]
	  \draw[step=1cm, opacity=0.5] (0,0) grid (3,3);

	  \node[ForestGreen, opacity=0.5] at (2.5, 2.5) {$\Gamma$};

	  \clip (-0.5, -0.5) rectangle (3.5, 3.5);

	  \draw[solid, ultra thick, ForestGreen, opacity=0.5]  plot[smooth,domain=-0.5:3.5] (\x, {0.5*(\x-1.5)+1.35+0.25*((\x-1.5)^2)-0.2});
	  
	  \draw[dotted, thick] plot[smooth,domain=-0.5:3.5] (\x, {1.3708333+(\x-1.5)-0.2});	  
	  \draw[solid, thick] plot[smooth,domain=1:2] (\x, {1.3708333+(\x-1.5)-0.2});
	  
      \fill[fill=NavyBlue, fill opacity=0.1] (0, 0) rectangle (1, 3);
      \fill[fill=NavyBlue, fill opacity=0.1] (2, 0) rectangle (3, 3);

	  \node at (1.5, 1.5) {$T$};
	\end{tikzpicture}
\end{minipage}
\hspace{0.25cm}
\begin{minipage}{.3\textwidth}
\centering
\begin{tikzpicture}[scale=0.9]
	  \draw[step=1cm, opacity=0.5] (0,0) grid (3,3);

	  \node[ForestGreen, opacity=0.5] at (2.5, 2.5) {$\Gamma$};

	  \clip (-0.5, -0.5) rectangle (3.5, 3.5);

	  \draw[solid, ultra thick, ForestGreen, opacity=0.5] plot[smooth,domain=-0.5:3.5] (\x, {0.5*(\x-1.5)+1.35+0.25*((\x-1.5)^2)-0.2});
	  
	  \draw[dotted, thick] plot[smooth,domain=-0.5:3.5] (\x, {1.3708333+0.75*(\x-1.5)-0.2});	  
	  \draw[solid, thick] plot[smooth,domain=1:2] (\x, {1.3708333+0.75*(\x-1.5)-0.2});
	  
      \fill[fill=NavyBlue, fill opacity=0.1] (1, 0) rectangle (3, 3);

	  \node at (1.5, 1.5) {$T$};
	\end{tikzpicture}
\end{minipage}
\caption{\corr{ELVIRA $3$ cases for the vertical orientation. The slope is estimated using the differences between the two first column averages (left), the first and third (center) and the last two columns (right). }}
\label{fig:elviracalculations}
\end{center}
\end{figure}

The following eight reconstruction strategies are compared:
\begin{enumerate}
    \item \textbf{\perplexityinsert{piecewise-constant}}: as mentioned in the introduction, the simplest linear method that one can come up with is $\t u:=\sum_{T\in \cT_h} a_T(u) \Chi_T$, that is, on each cell $T$ the value given by the observed cell average $a_T(u)$.

    \item \textbf{\perplexityinsert{linear-obera}}: 
 we apply the minimization strategy described in \S\ref{sec:obera} for \textit{linear interfaces} as in Example 1, with loss function $\cL(u, v)=\|a_S(u)-a_S(v)\|$ using $\ell^1$ norm and a $3\times 3$ stencil $S$.
    
  \item \textbf{\perplexityinsert{linear-obera-w}}: 
  we apply the same approach but enforcing area consistency on $T$ through the weighted loss function \iref{eq:kloss} with $K=100.$
    
    \item \textbf{\perplexityinsert{elvira}}: following \cite{Pilliod2004}, the ELVIRA method consists on replacing, still for \textit{linear interfaces}, the OBERA continuous optimization strategy by providing only $6$ combinations of parameters $\mu$ from which to choose $\mu^*$, the minimizer of $\cL$ with $\ell^2$ norm instead of $\ell^1$. The $6$ alternatives are obtained by proposing, for the interface's slope, the $6$ possible finite differences estimations of a $3\times 3$ stencil \corr{(see Figure \ref{fig:elviracalculations})}.
    
    \item \textbf{\perplexityinsert{elvira-w-oriented}}: 
  we apply ELVIRA but choosing first an {\it orientation}, as in AEROS, which allows to reduce the choice to $3$ alternatives. In addition, we work with the modified
  {\it weighted} loss function, with $K=100$, to favor area consistency on $T$. 

    \item \textbf{\perplexityinsert{quadratic-obera-non-adaptive}}: we apply the minimization strategy for \textit{quadratic interfaces} after choosing an orientation to have a Cartesian parametrization of the interface, that is, $v_\mu=\Chi_P$, as in Example 3, and $p\in \pP_2$ the space of univariate polynomials of degree $2$. Here, we also use a $3 \times 3$ fixed stencil and enforce area consistency on $T$ through the same modified loss function \iref{eq:kloss} with $K=100$.

    \item \textbf{\perplexityinsert{quadratic-aero}}: 
    we apply the AEROS reconstruction strategy for \textit{quadratic interfaces} with the adaptive method, described in \S\ref{subsec:aeros}, to build stencils of width $3$ minimal height.
      
    \item \textbf{\perplexityinsert{quartic-aero}}: 
    we apply the  AEROS strategy now with polynomials of degree $4$, therefore with stencils of width $5$.
\end{enumerate}

\begin{figure}
     \begin{subfigure}[b]{0.22\textwidth}
         \centering
         \includegraphics[width=\linewidth]{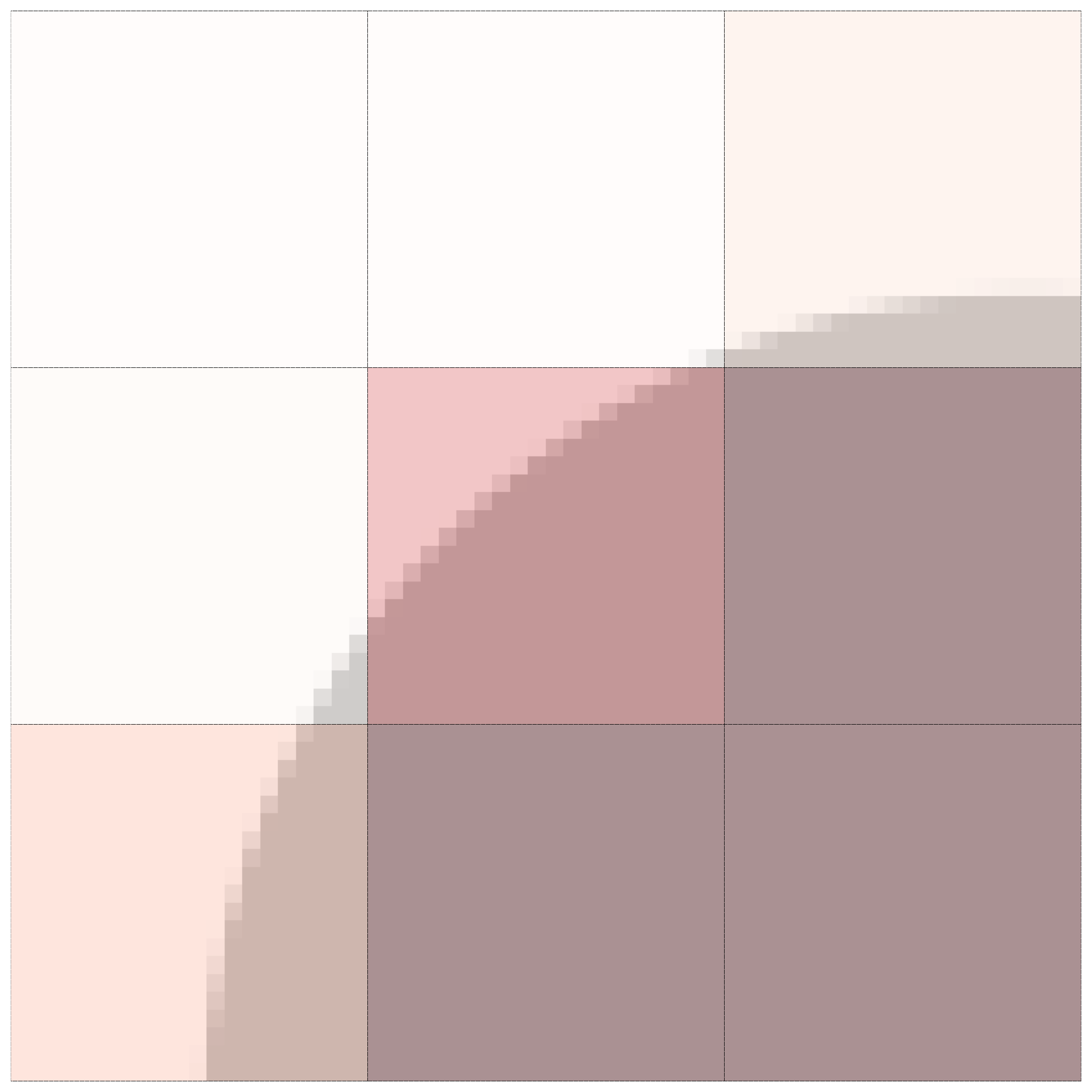}
         \caption{\perplexityinsert{piecewise-constant}.}
         \label{fig:pwc}
     \end{subfigure}
     \hfill
     \begin{subfigure}[b]{0.22\textwidth}
         \centering
         \includegraphics[width=\linewidth]{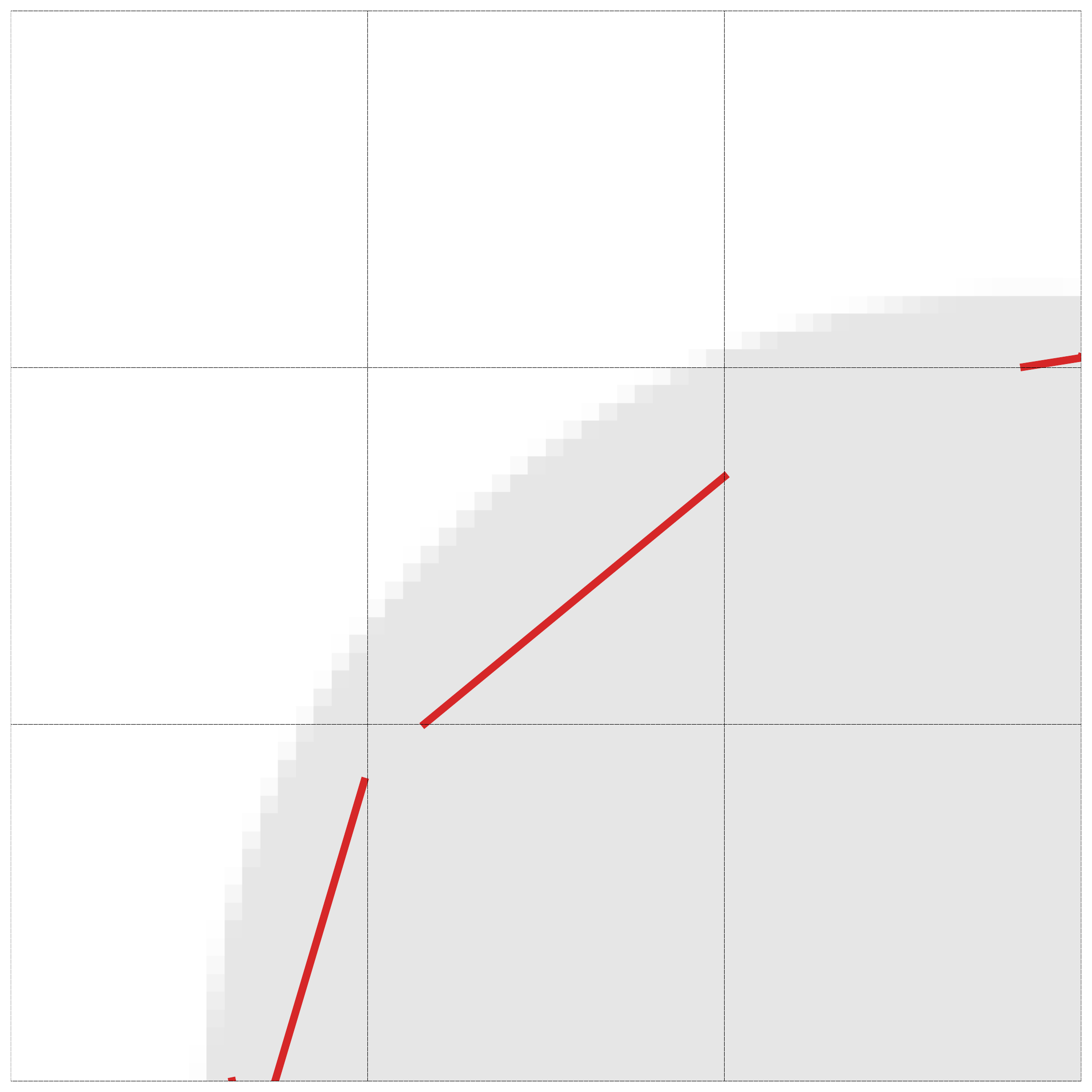}
         \caption{\perplexityinsert{linear-obera}.}
         \label{fig:obera-linear}
     \end{subfigure}
     \hfill
     \begin{subfigure}[b]{0.22\textwidth}
         \centering
         \includegraphics[width=\linewidth]{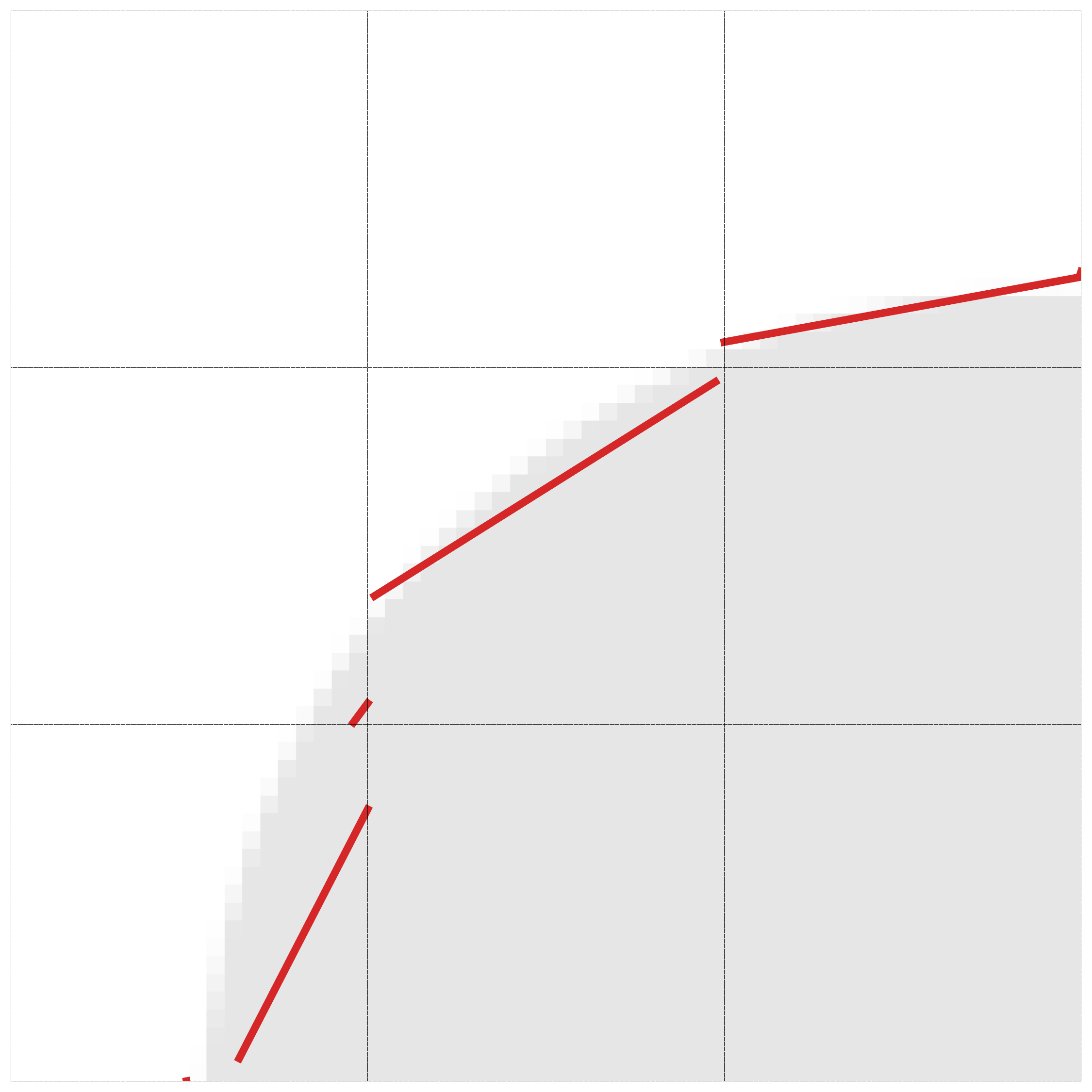}
         \caption{\perplexityinsert{elvira}.}
         \label{fig:elvira}
     \end{subfigure}
     \hfill
     \begin{subfigure}[b]{0.22\textwidth}
         \centering
         \includegraphics[width=\linewidth]{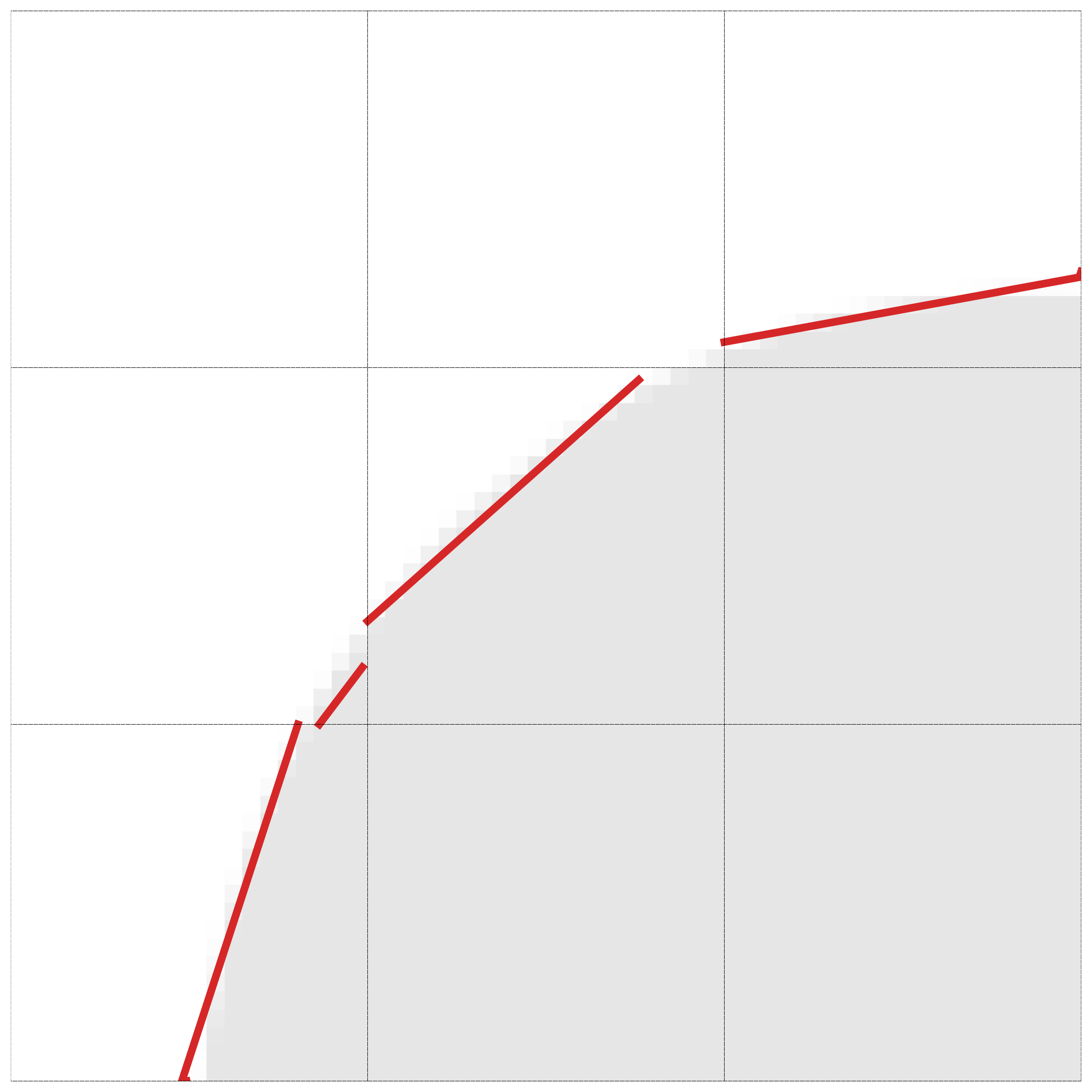}
         \caption{\perplexityinsert{elvira-w-oriented}.}
         \label{fig:elviraw}
     \end{subfigure}

     \begin{subfigure}[b]{0.22\textwidth}
         \centering
         \includegraphics[width=\linewidth]{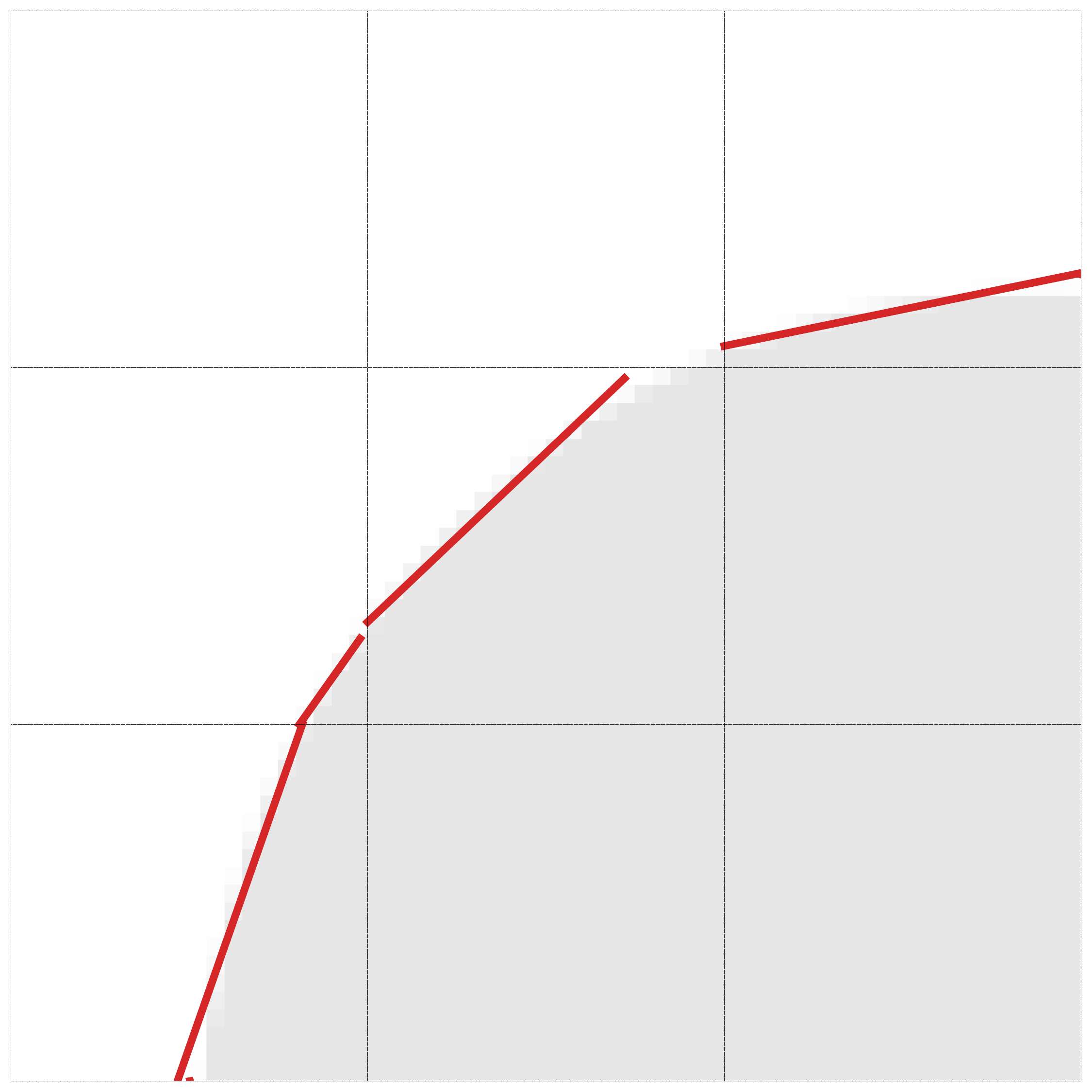}
         \caption{\perplexityinsert{linear-obera-w}.}
         \label{fig:obera-linear-w}
     \end{subfigure}
     \hfill
     \begin{subfigure}[b]{0.22\textwidth}
         \centering
         \includegraphics[width=\linewidth]{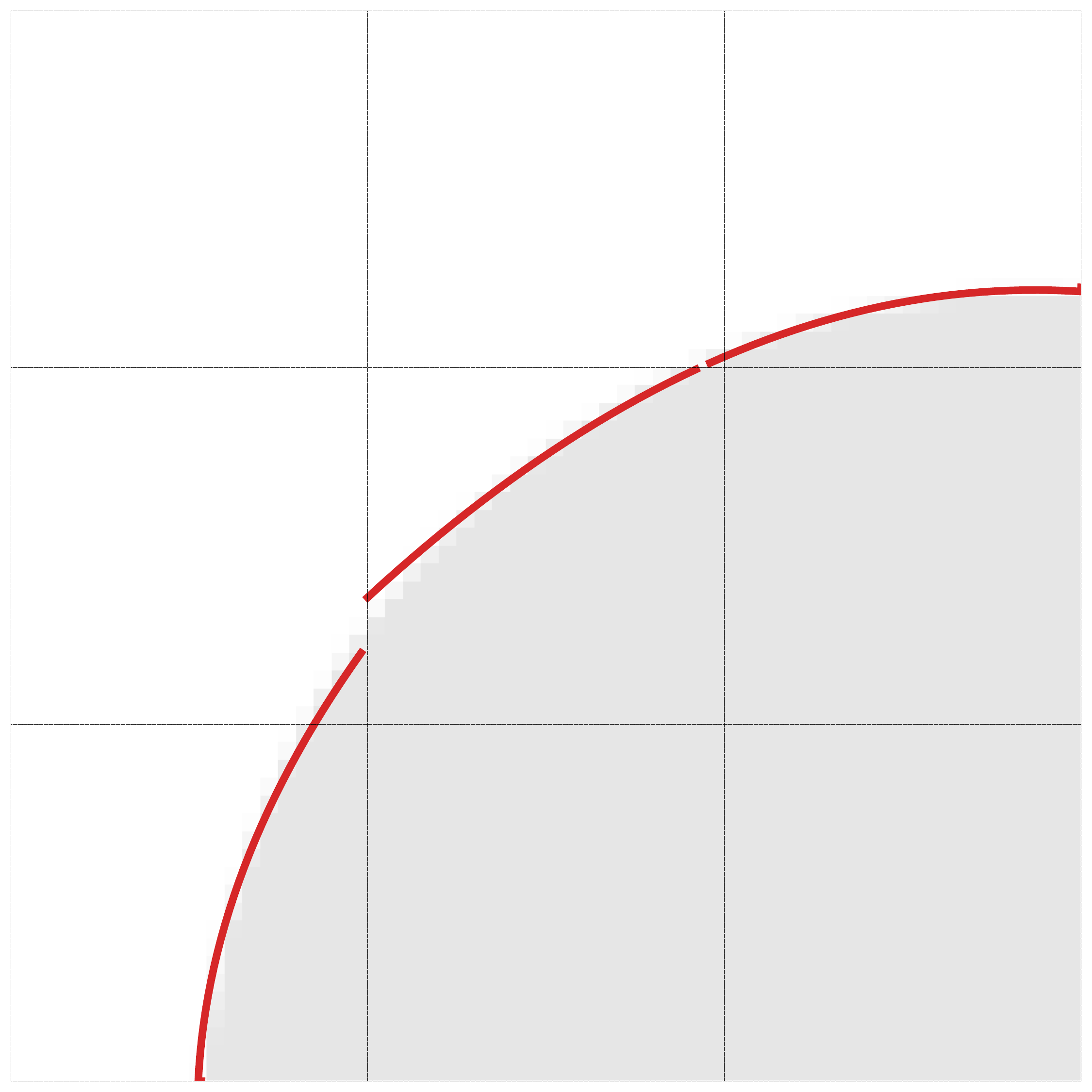}
         \caption{\perplexityinsert{quadratic-aero}.}
         \label{fig:aero-quadratic}
     \end{subfigure}
     \hfill
     \begin{subfigure}[b]{0.22\textwidth}
         \centering
         \includegraphics[width=\linewidth]{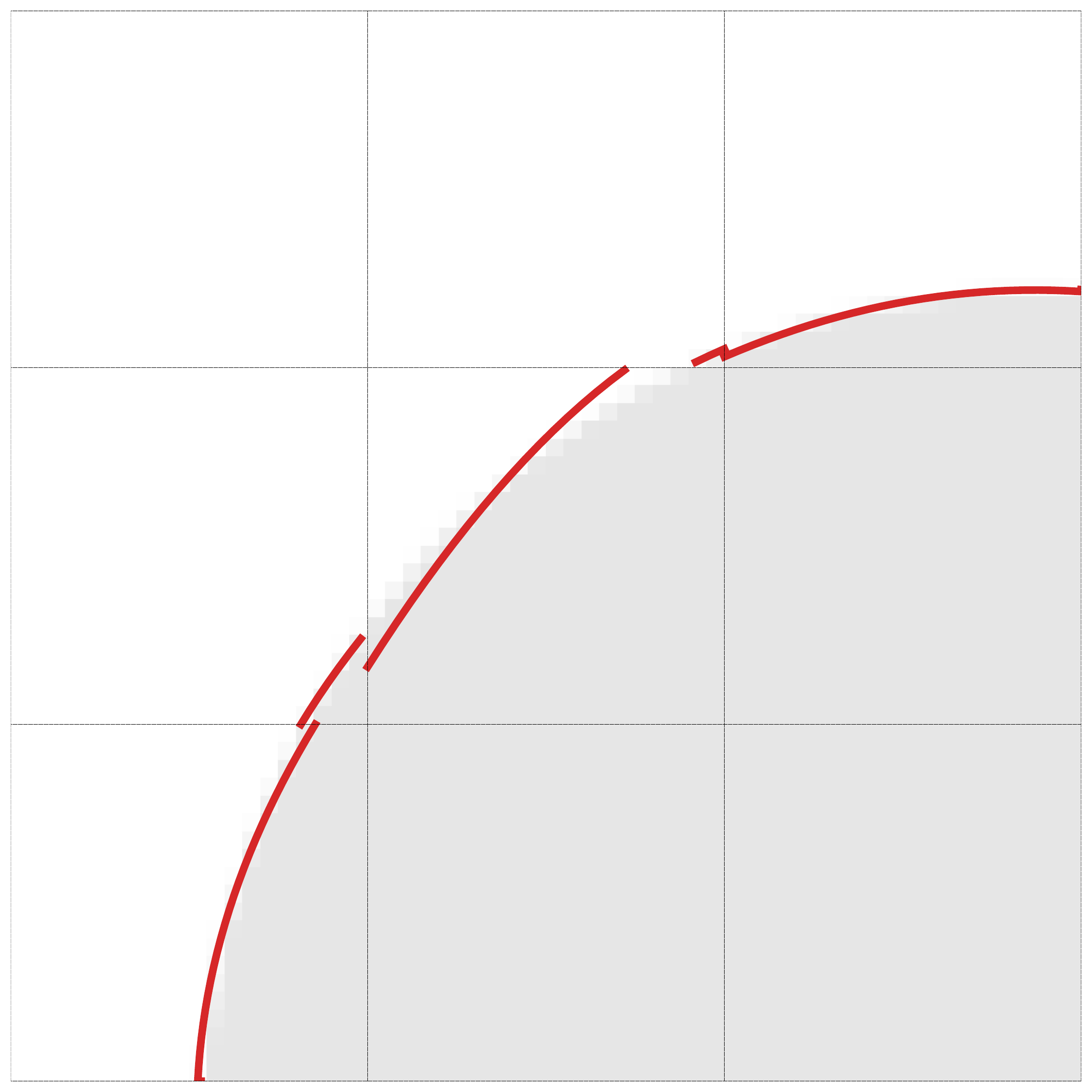}
         \caption{\perplexityinsert{quadratic-obera-non-adaptive}.}
         \label{fig:obera-quadratic}
     \end{subfigure}
     \hfill
     \begin{subfigure}[b]{0.22\textwidth}
         \centering
         \includegraphics[width=\linewidth]{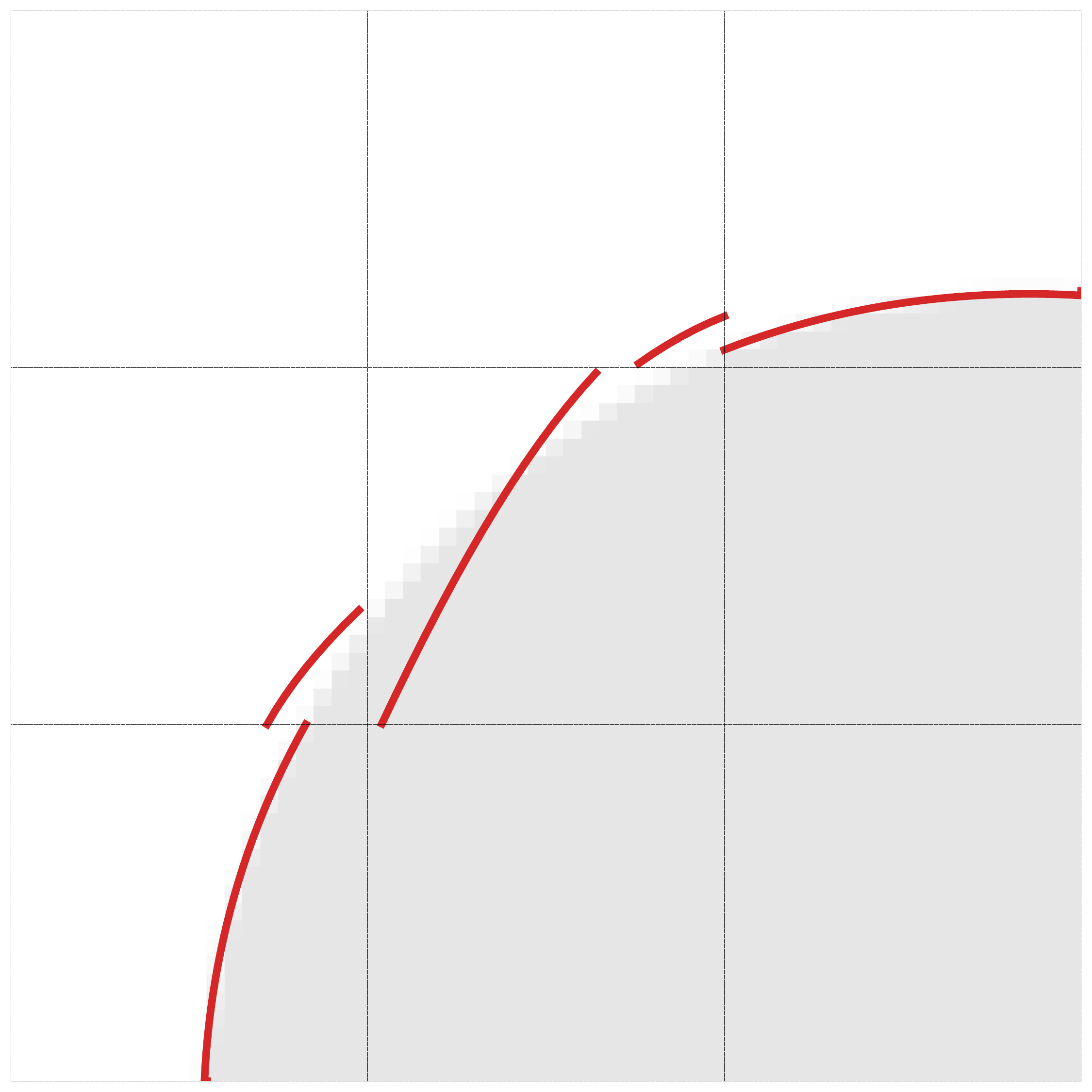}
         \caption{\perplexityinsert{quartic-aero}.}
         \label{fig:aero-quartic}
     \end{subfigure}
    \caption{Reconstruction of a portion of the circle by different methods for a scale of $h=1/10$.}
     \label{fig:visual10}
\end{figure}

\begin{figure}
     \begin{subfigure}[b]{0.22\textwidth}
         \centering
         \includegraphics[width=\linewidth]{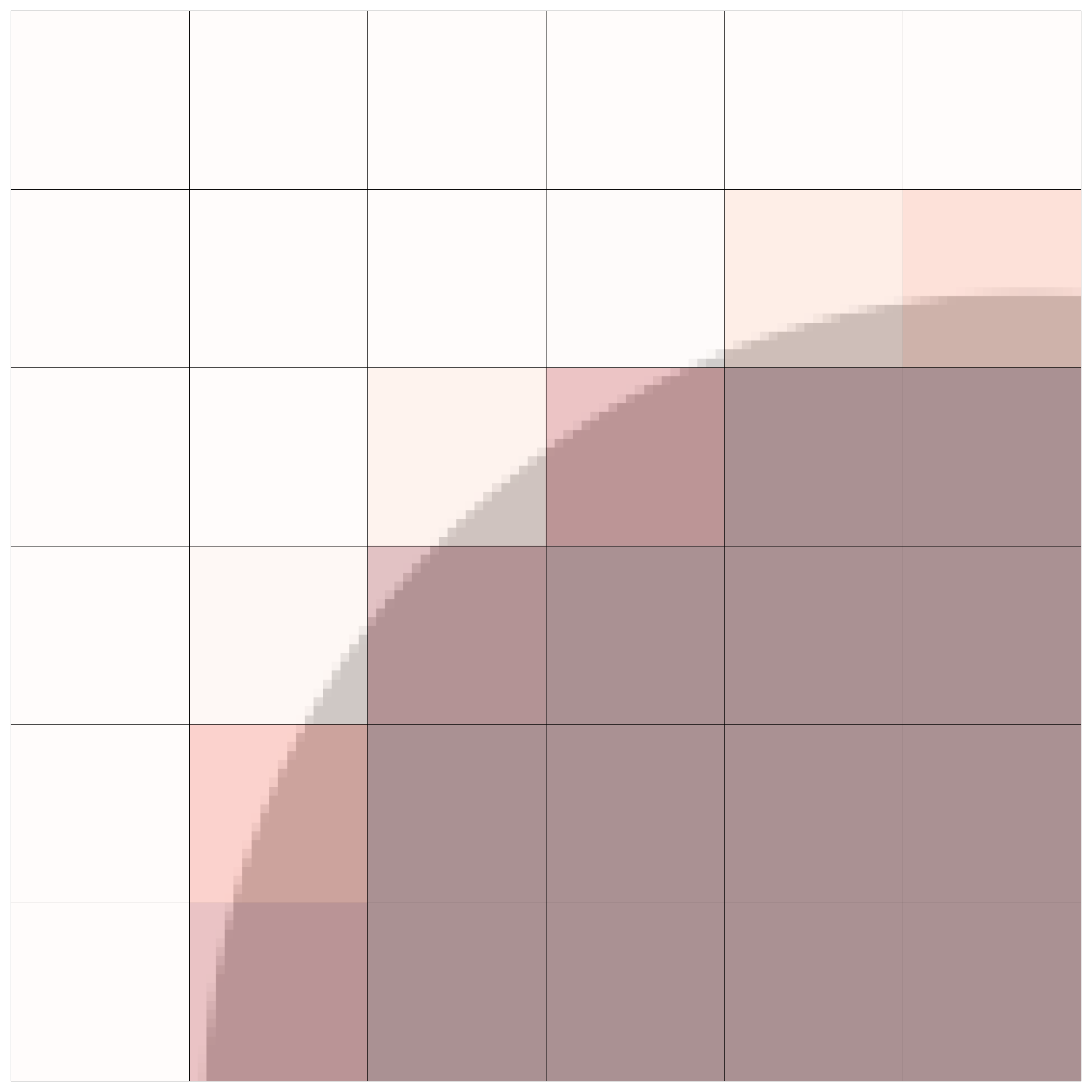}
         \caption{\perplexityinsert{piecewise-constant}.}
         \label{fig:pwc}
     \end{subfigure}
     \hfill
     \begin{subfigure}[b]{0.22\textwidth}
         \centering
         \includegraphics[width=\linewidth]{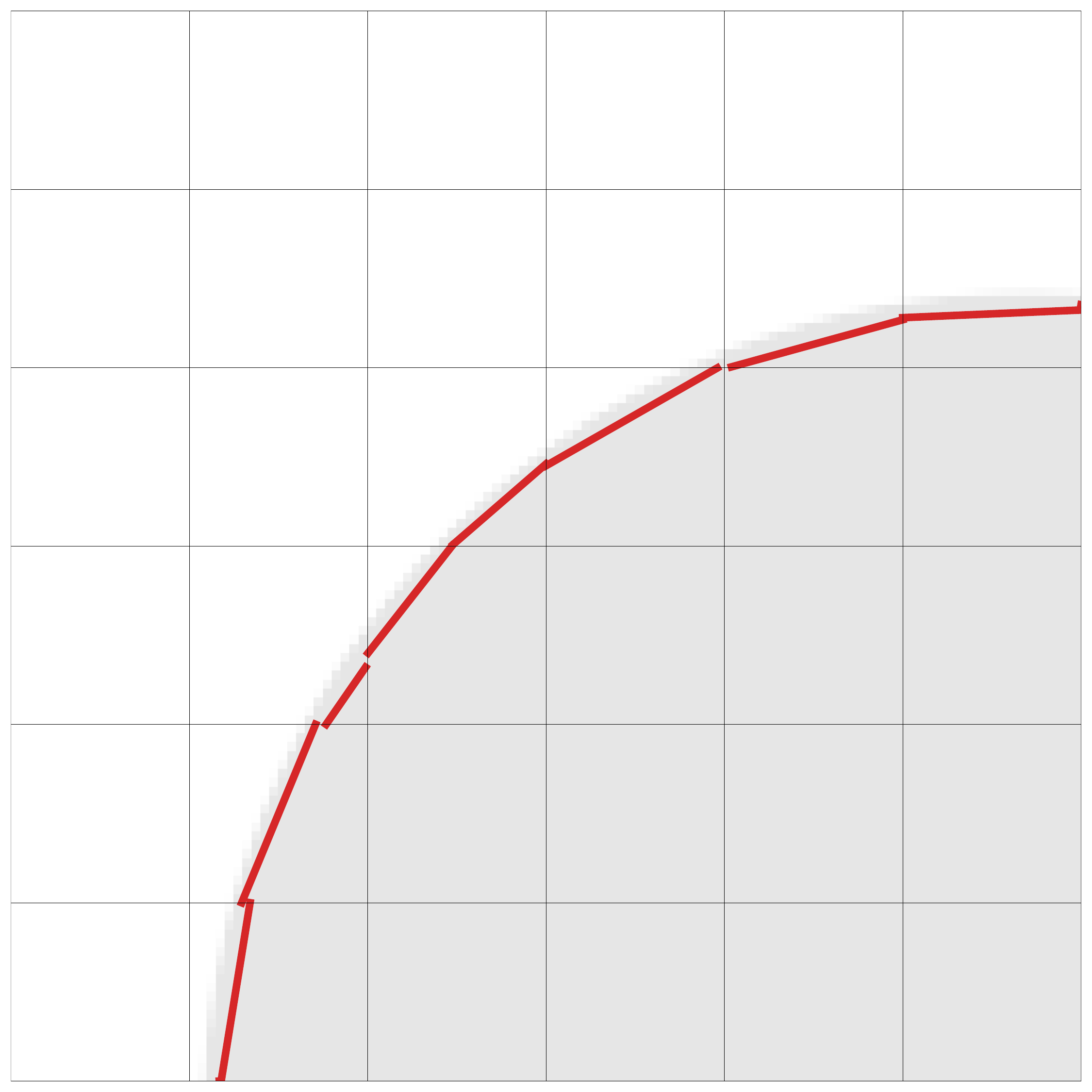}
         \caption{\perplexityinsert{linear-obera}.}
         \label{fig:obera-linear}
     \end{subfigure}
     \hfill
     \begin{subfigure}[b]{0.22\textwidth}
         \centering
         \includegraphics[width=\linewidth]{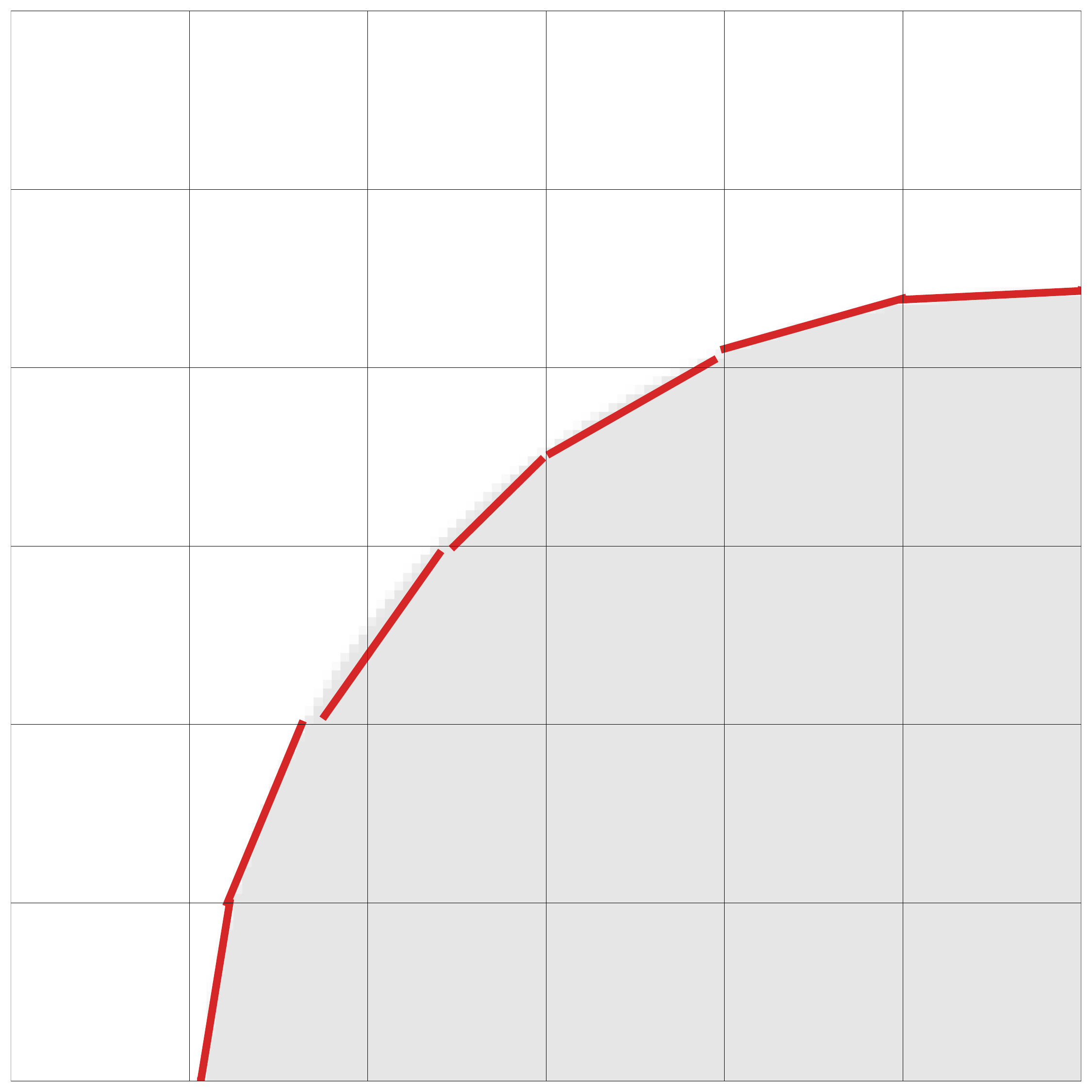}
         \caption{\perplexityinsert{elvira}.}
         \label{fig:elvira}
     \end{subfigure}
     \hfill
     \begin{subfigure}[b]{0.22\textwidth}
         \centering
         \includegraphics[width=\linewidth]{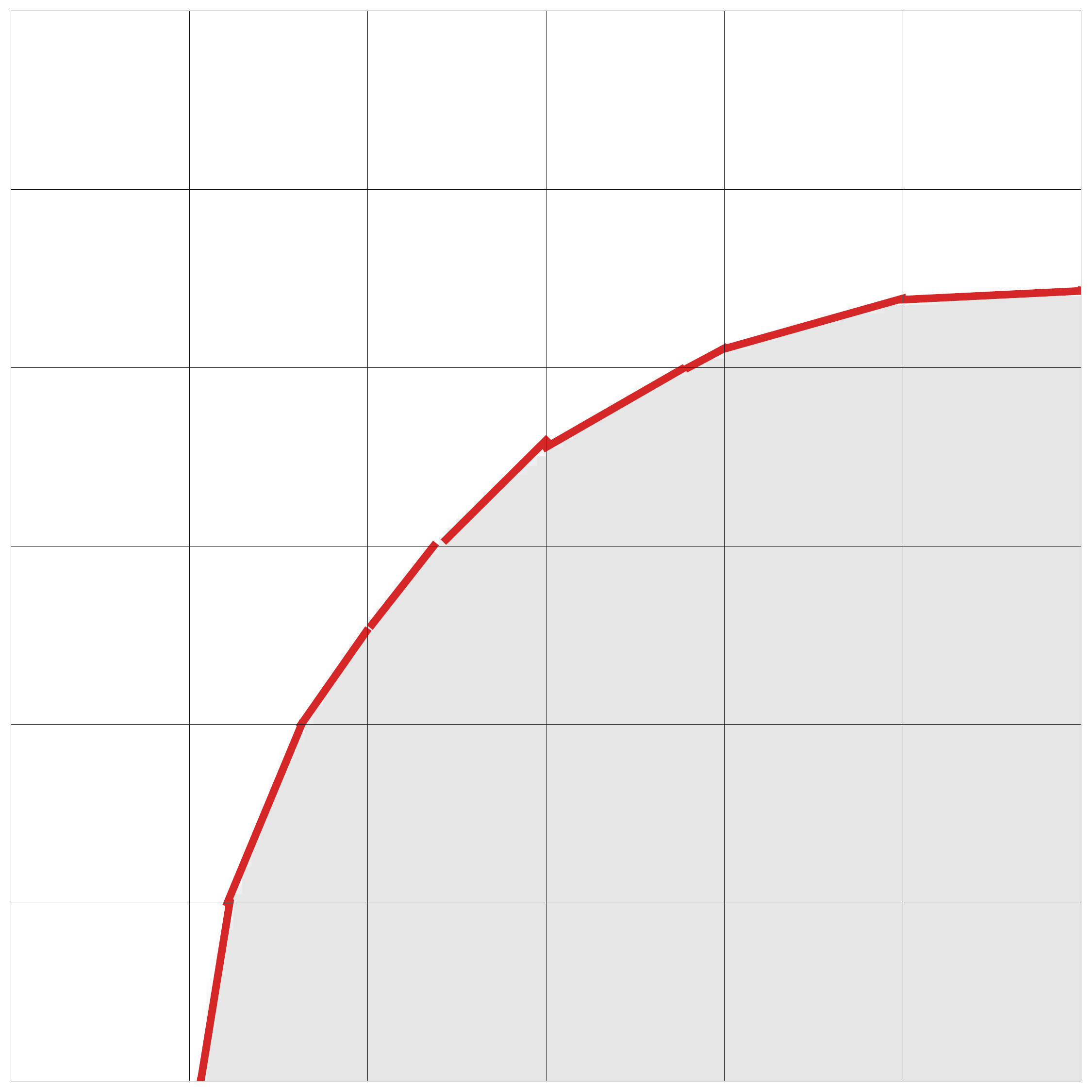}
         \caption{\perplexityinsert{elvira-w-oriented}.}
         \label{fig:elviraw}
     \end{subfigure}

     \begin{subfigure}[b]{0.22\textwidth}
         \centering
         \includegraphics[width=\linewidth]{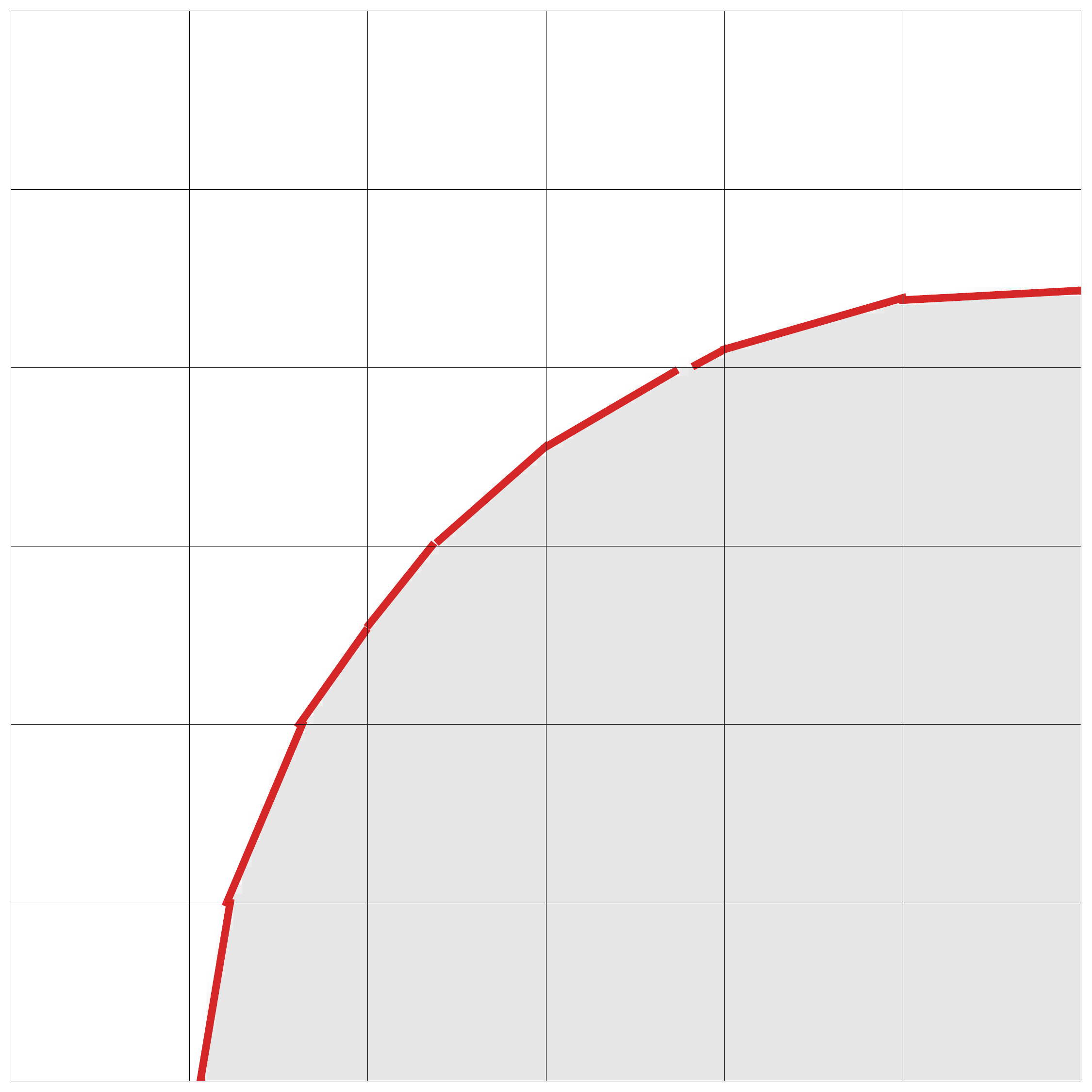}
         \caption{\perplexityinsert{linear-obera-w}.}
         \label{fig:obera-linear-w}
     \end{subfigure}
     \hfill
     \begin{subfigure}[b]{0.22\textwidth}
         \centering
         \includegraphics[width=\linewidth]{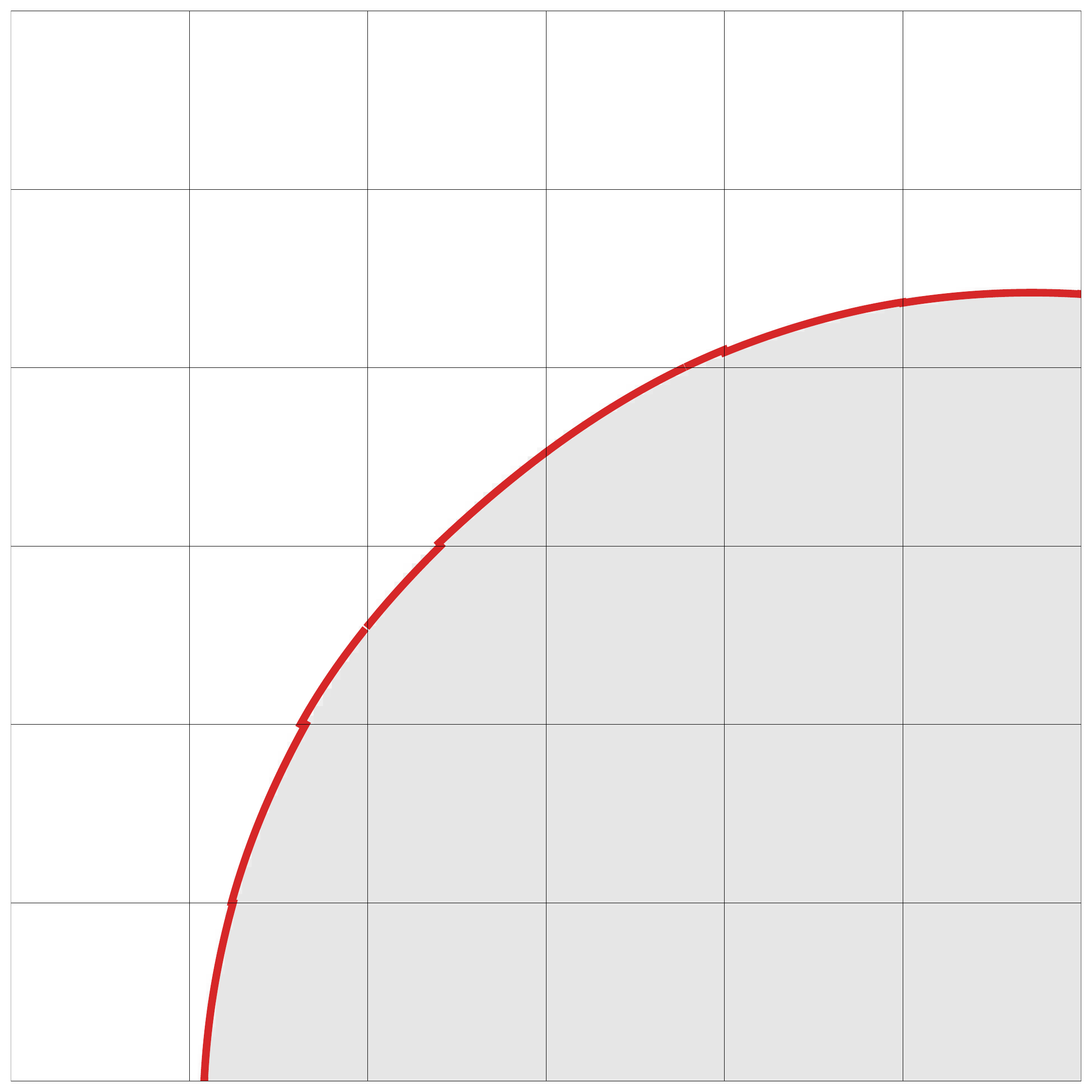}
         \caption{\perplexityinsert{quadratic-aero}.}
         \label{fig:aero-quadratic}
     \end{subfigure}
     \hfill
     \begin{subfigure}[b]{0.22\textwidth}
         \centering
         \includegraphics[width=\linewidth]{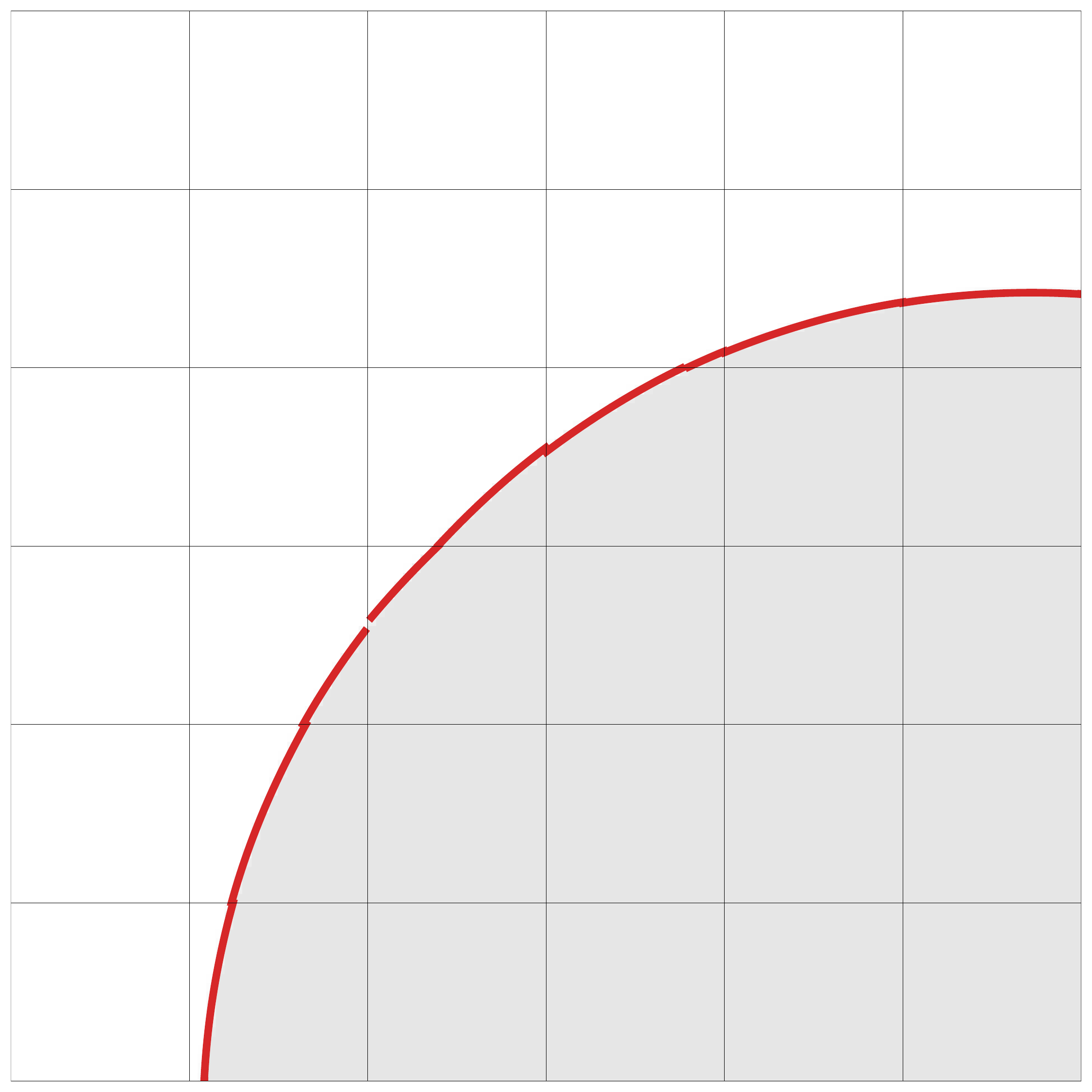}
         \caption{\perplexityinsert{quadratic-obera-non-adaptive}.}
         \label{fig:obera-quadratic}
     \end{subfigure}
     \hfill
     \begin{subfigure}[b]{0.22\textwidth}
         \centering
         \includegraphics[width=\linewidth]{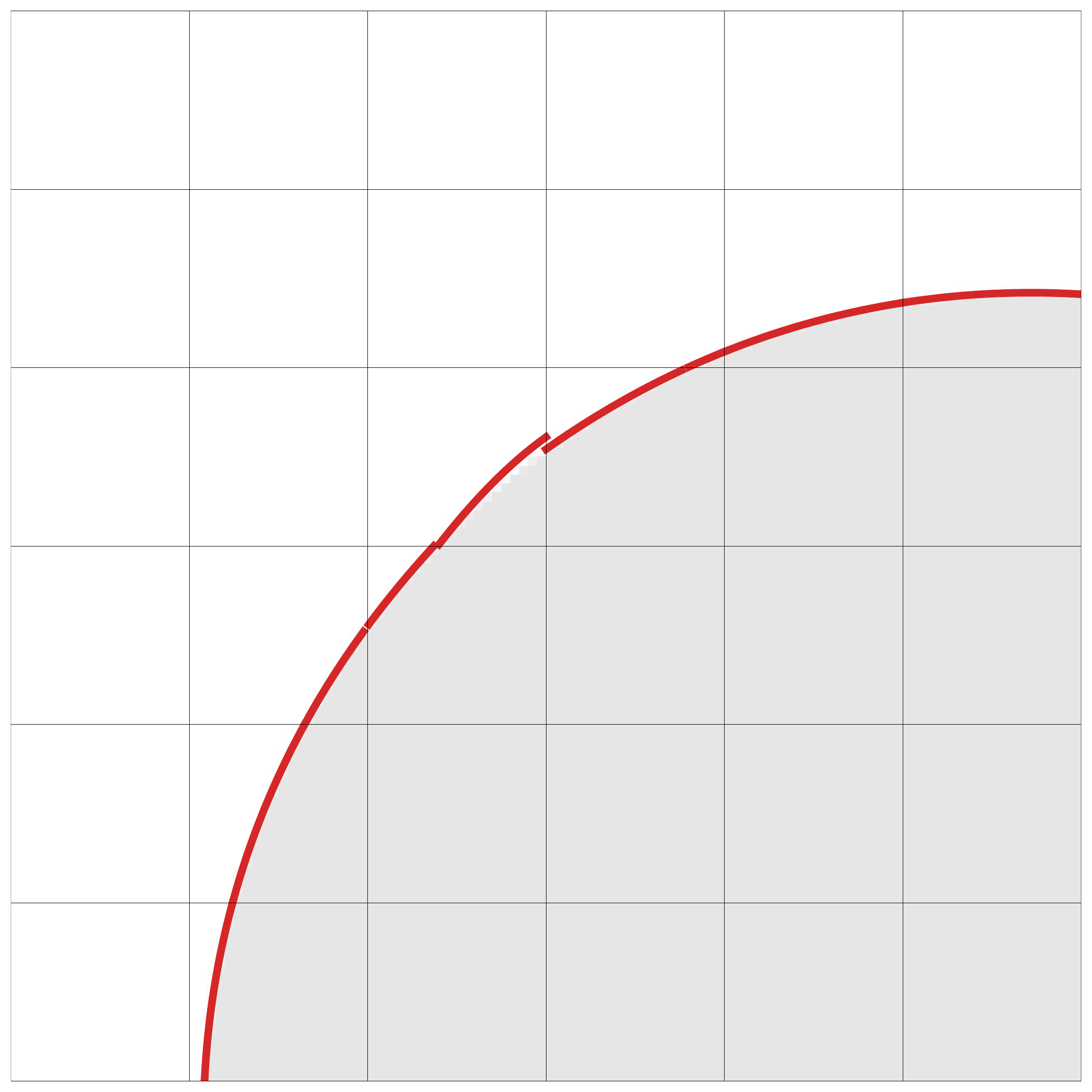}
         \caption{\perplexityinsert{quartic-aero}.}
         \label{fig:aero-quartic}
     \end{subfigure}
    \caption{Reconstruction of a portion of the circle by different methods for a scale of $h=1/20$.}
     \label{fig:visual20}
\end{figure}

Figure \ref{fig:visual10} and Figure \ref{fig:visual20}.
display a detail of the recovery with $h=1/10$ and $h=1/20$
respectively in order to compare the visual quality.

\begin{figure}[!h]
     \centering
     \includegraphics[width=0.8\linewidth]{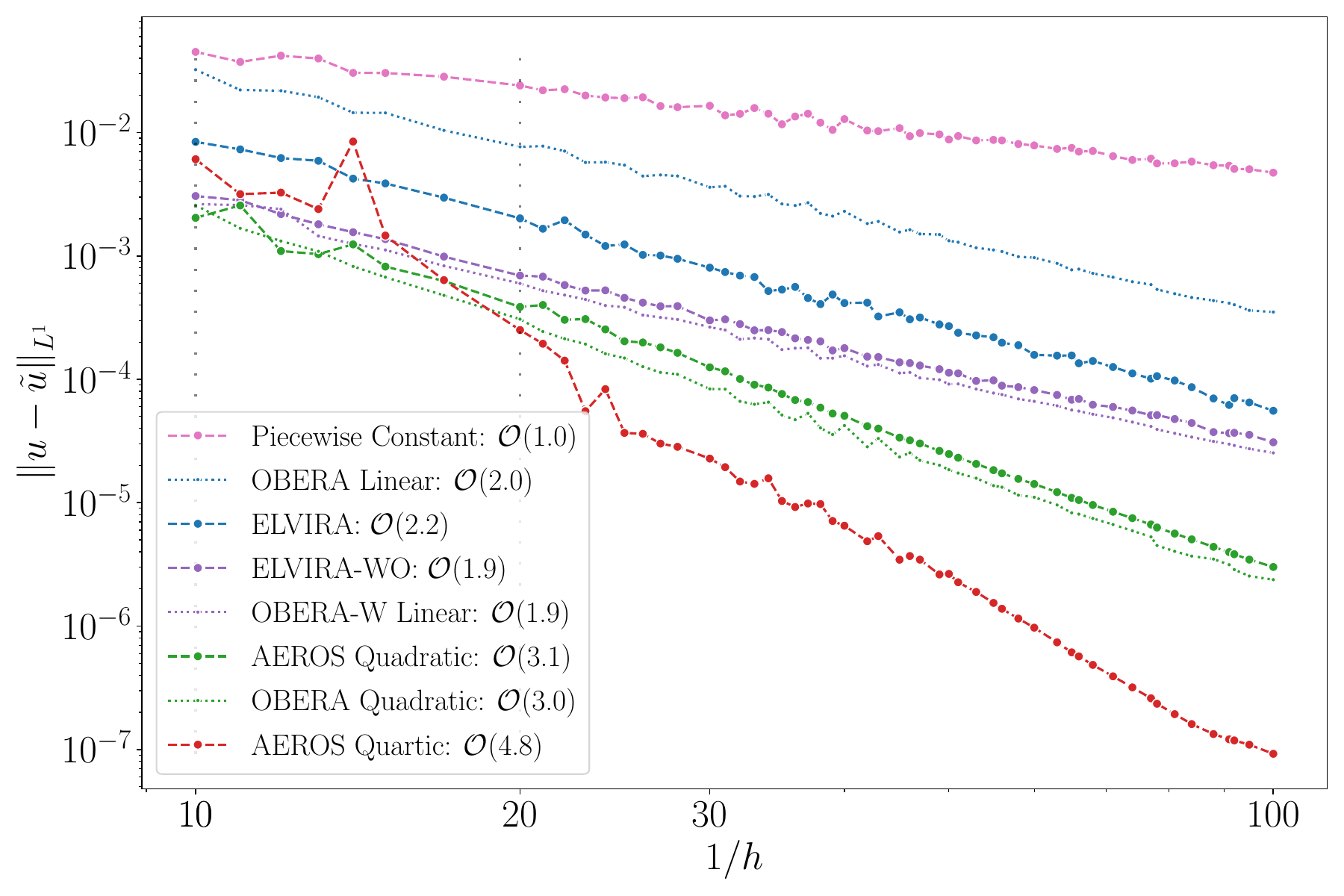}
     \caption{Convergence for different reconstruction models. The convergence rates (in parenthesis) are estimated using values $1/h>30$.}
     \label{fig:convergence}
\end{figure}

As to linear interfaces, we clearly notice the 
relevance of enforcing area consistency 
on the cell $T$ of interest
by appropriately modifying the loss function $\cL$. 
Although the four methods 
benefit from similar convergence rates of $\cO(h^2)$ as expected from
\iref{col:lvira-convergence} and \iref{th:aeros-convergence}, see Figure \ref{fig:convergence}, 
the reconstruction error improvement is of an order of magnitude 
by this simple change. When area consistency is not imposed, the linear interfaces are
pushed towards the interior of the circle as the curvature of the original domain points inwardly.

As to the AEROS strategies, we notice that
for $h=1/10$, they have difficulties to reconstruct the interface 
on cells where $\Gamma$ 
is rapidly passing from a situation where an horizontal orientation
is preferred to another in which a vertical one is better.
This effect if particulary evident for the case of \perplexityinsert{quartic-aero}
as the method needs a stencil with $5$ columns. At this scale
we are above the critical scale (see Figure \ref{fig:criticalscale}) in which for
some cells there is no stencil of the needed width allowing the curve to be described as a graph.
This problem disappears for $h=1/20$, for which these
methods have the best visual quality.

On Figure \ref{fig:convergence}, we see that in terms of convergence we obtain for both AEROS the
expected rates from \iref{th:aeros-convergence}: for quadratics we get $\cO(h^3)$ and we almost get $\cO(h^5)$
for polynomials of degree $4$. As mentioned before and graphically shown in Figure \ref{fig:visual10},
\perplexityinsert{quartic-aero} breaks 
down when the scale of the discretization is above the
critical scale which, for this particular, example is around $h=1/20$. 

Regarding the computational time per cell taken by each algorithm we observe on Table \ref{tab:timesmodels} that
OBERA strategies are two orders of magnitude slower than any AEROS approach. At the same time,
although ELVIRA methods are faster than OBERA, they still are an order slower than AEROS
due to the bottleneck of having to compute the stencil cell averages to compare the
$6$ or $3$ alternatives under evaluation. This limit is justified by the fact that choosing
an orientation, as in \perplexityinsert{elvira-w-oriented}, cuts by half the computing time of the
overall algorithm, while the improvement in accuracy is achieved by modifying the loss function. 

In summary, all three comparisons in terms of visual aspect, order of convergence and computational time are in favor of the AEROS
strategy provided that $h$ is below the critical scale.

\begin{figure}
     \begin{subfigure}[b]{0.45\textwidth}
         \centering
         \includegraphics[width=\linewidth]{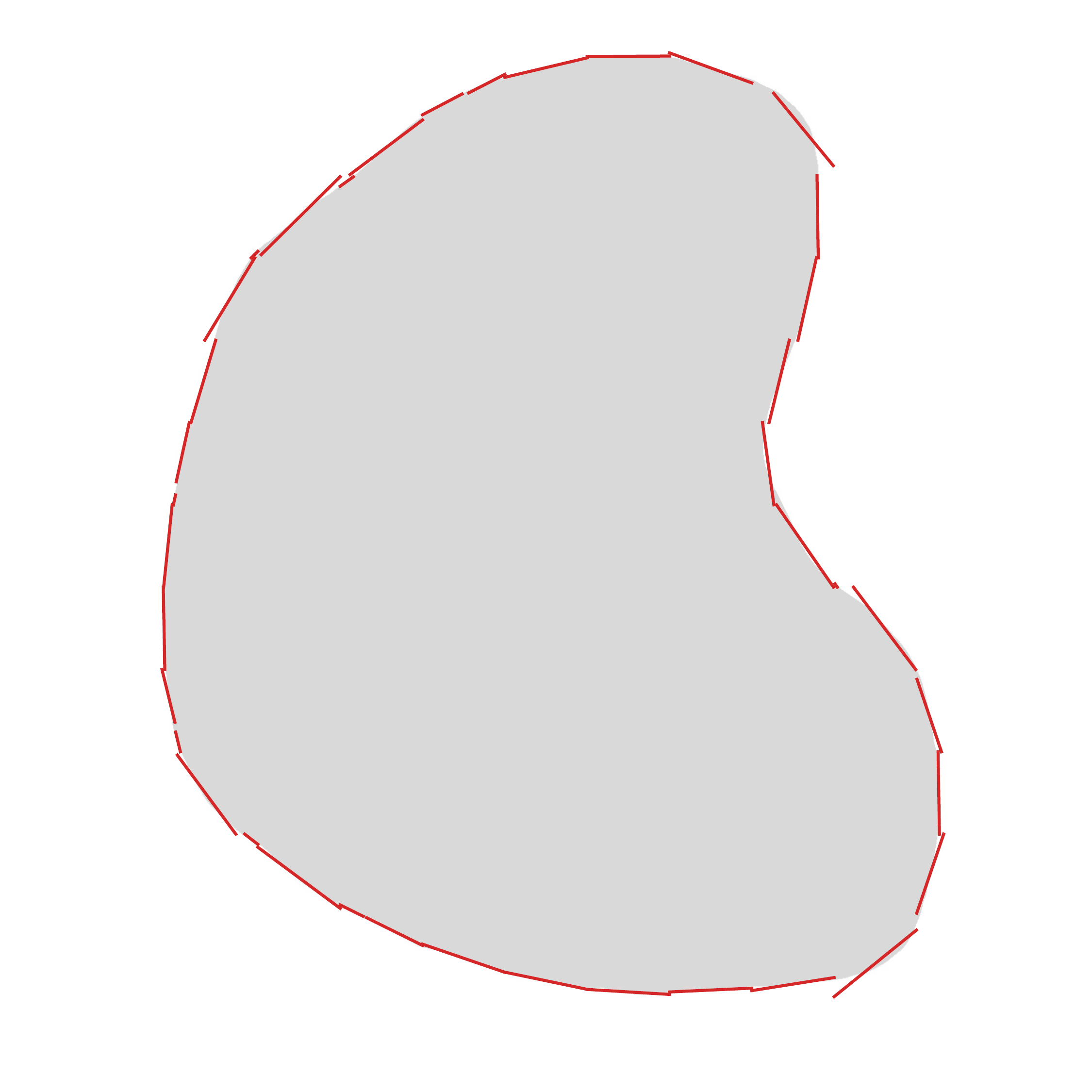}
         \caption{\perplexityinsert{linear-obera-w}.}
         \label{fig:batata-obera-linear-w}
     \end{subfigure}
     \hfill
     \begin{subfigure}[b]{0.45\textwidth}
         \centering
         \includegraphics[width=\linewidth]{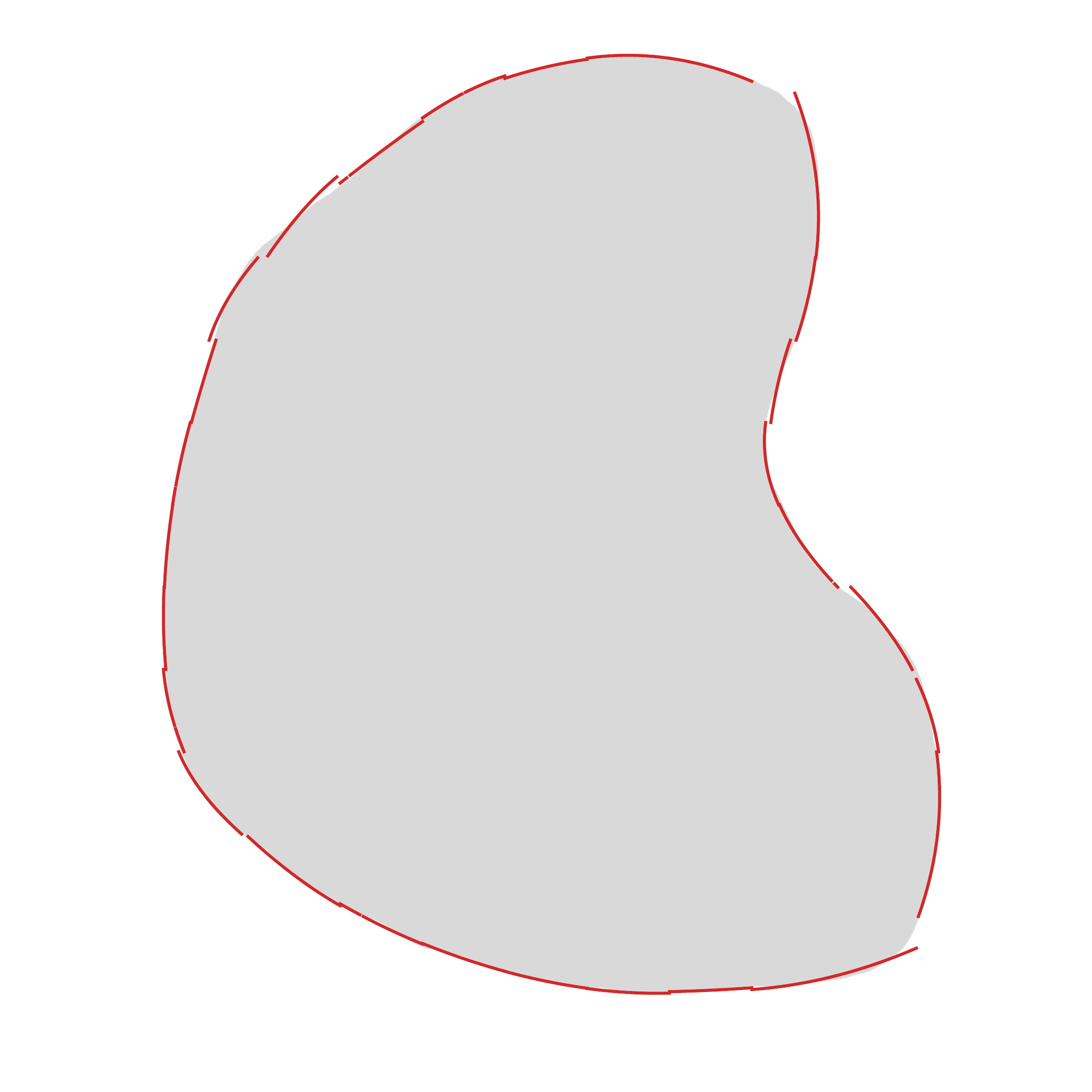}
         \caption{\perplexityinsert{quadratic-aero}.}
         \label{fig:batata-aero-quadratic}
     \end{subfigure}
     \hfill
     \begin{subfigure}[b]{0.45\textwidth}
         \centering
         \includegraphics[width=\linewidth]{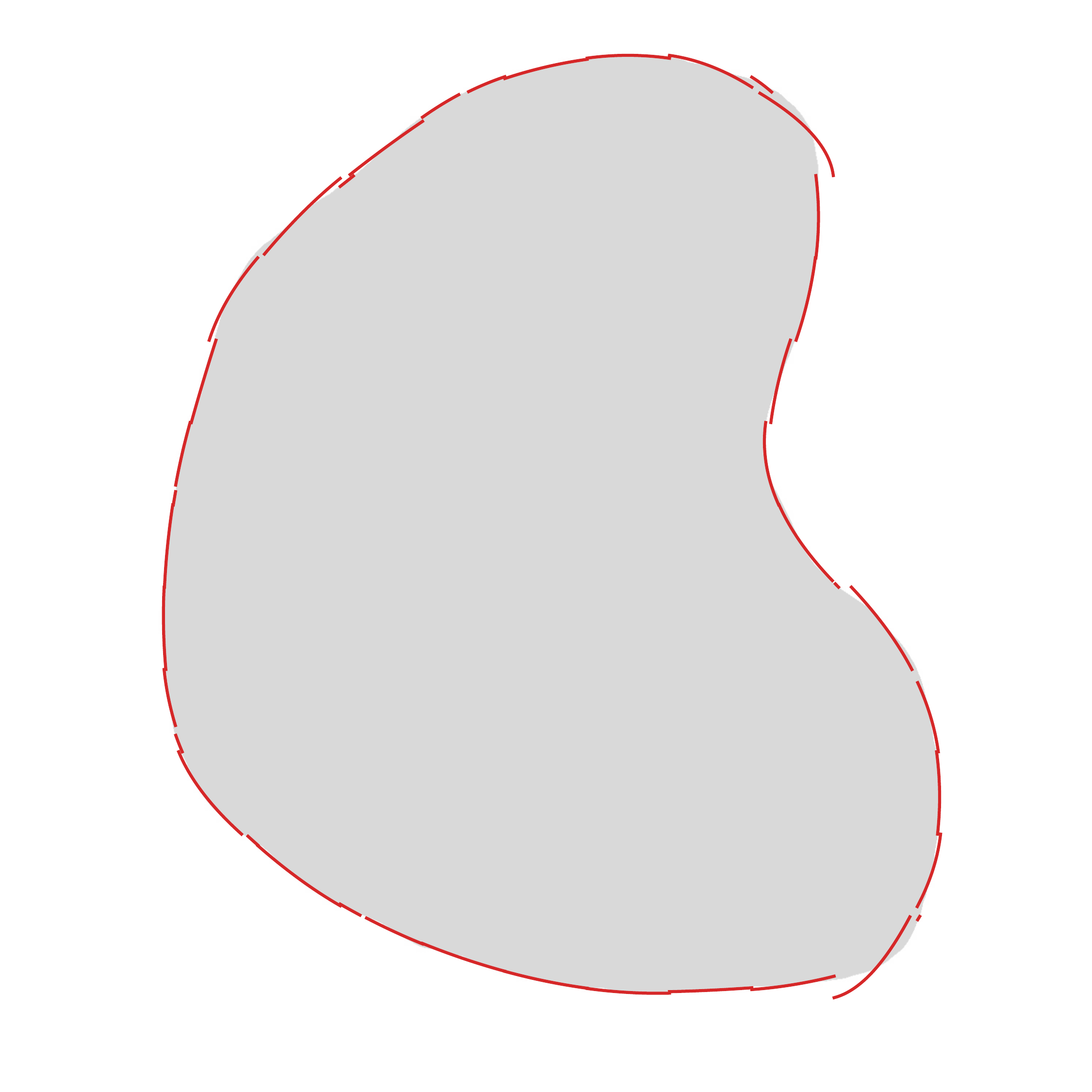}
         \caption{\perplexityinsert{quadratic-obera-non-adaptive}.}
         \label{fig:batata-obera-quadratic}
     \end{subfigure}
     \hfill
     \begin{subfigure}[b]{0.45\textwidth}
         \centering
         \includegraphics[width=\linewidth]{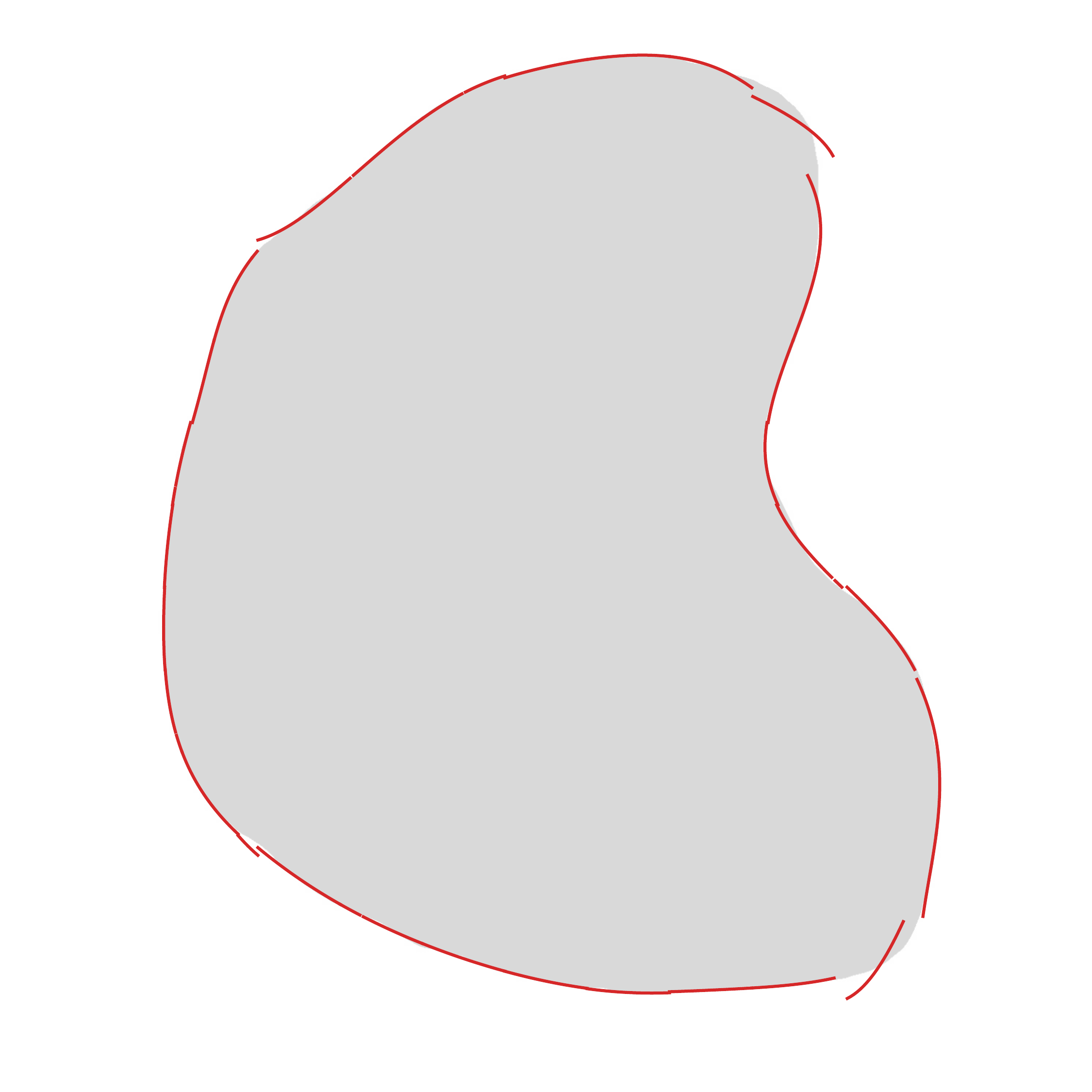}
         \caption{\perplexityinsert{quartic-aero}.}
         \label{fig:batata-aero-quartic}
     \end{subfigure}
    \caption{Reconstruction of a smooth domain by different methods for a scale of $h=1/15$. \corr{The higher the order the smoother the reconstruction except in the transitions of orientation where higher order methods like AEROS Quartic struggle because there is not enough information or the stencil is too much de-centered.}}
     \label{fig:smoothrec}
\end{figure}

Finally, we show in Figure \ref{fig:smoothrec} how the best linear interface method, \perplexityinsert{linear-obera-w}, and the three higher order methods compare when used to reconstruct an arbitrary, still smooth, domain. We see that by passing from the linear interface method to \perplexityinsert{quadratic-aero} the reconstruction becomes smoother while still suffering from some imperfections, notably in regions where there is a stronger change in the orientation of the curve. This is slightly improved by the optimization done in \perplexityinsert{quadratic-obera-non-adaptive} method. An even smoother result is obtained with \perplexityinsert{quartic-aero} at the expense of some small deviations again in regions where the orientation is rapidly changing.

\FloatBarrier

\FloatBarrier

\begin{table}[!h]
\centering
\begin{tabular}{ c | l }
 \perplexityinsert{linear-obera}  &
 $\mathbf{\perplexityinsert{linear-obera-time}}$
\\

 \perplexityinsert{linear-obera-w} &
 $\mathbf{\perplexityinsert{linear-obera-w-time}}$
\\

 \perplexityinsert{quadratic-obera-non-adaptive}  &
 $\mathbf{\perplexityinsert{quadratic-obera-non-adaptive-time}}$
\\
 \hline

 \perplexityinsert{elvira} &
 $\mathbf{\perplexityinsert{elvira-time}}$
\\

 \perplexityinsert{elvira-w-oriented} &
 $\mathbf{\perplexityinsert{elvira-w-oriented-time}}$
\\
 \hline

 \perplexityinsert{quadratic-aero} &
 $\mathbf{\perplexityinsert{quadratic-aero-time}}$
\\

 \perplexityinsert{quartic-aero} &
 $\mathbf{\perplexityinsert{quartic-aero-time}}$
\\
\end{tabular}
\caption{Average time (in seconds) taken to find the parameters of the interface by the different tested models. \corr{The average is taken over all instances in which each algorithm was called (to
perform a local approximation) to produce Figure \ref{fig:convergence} (which is in the order of the $4000$ per method).}}
\label{tab:timesmodels}
\end{table}

\FloatBarrier

\FloatBarrier

\subsection{The treatment of corner domains}
\label{sec:vertices}

\begin{figure}
         \centering
         \includegraphics[width=0.8\linewidth]{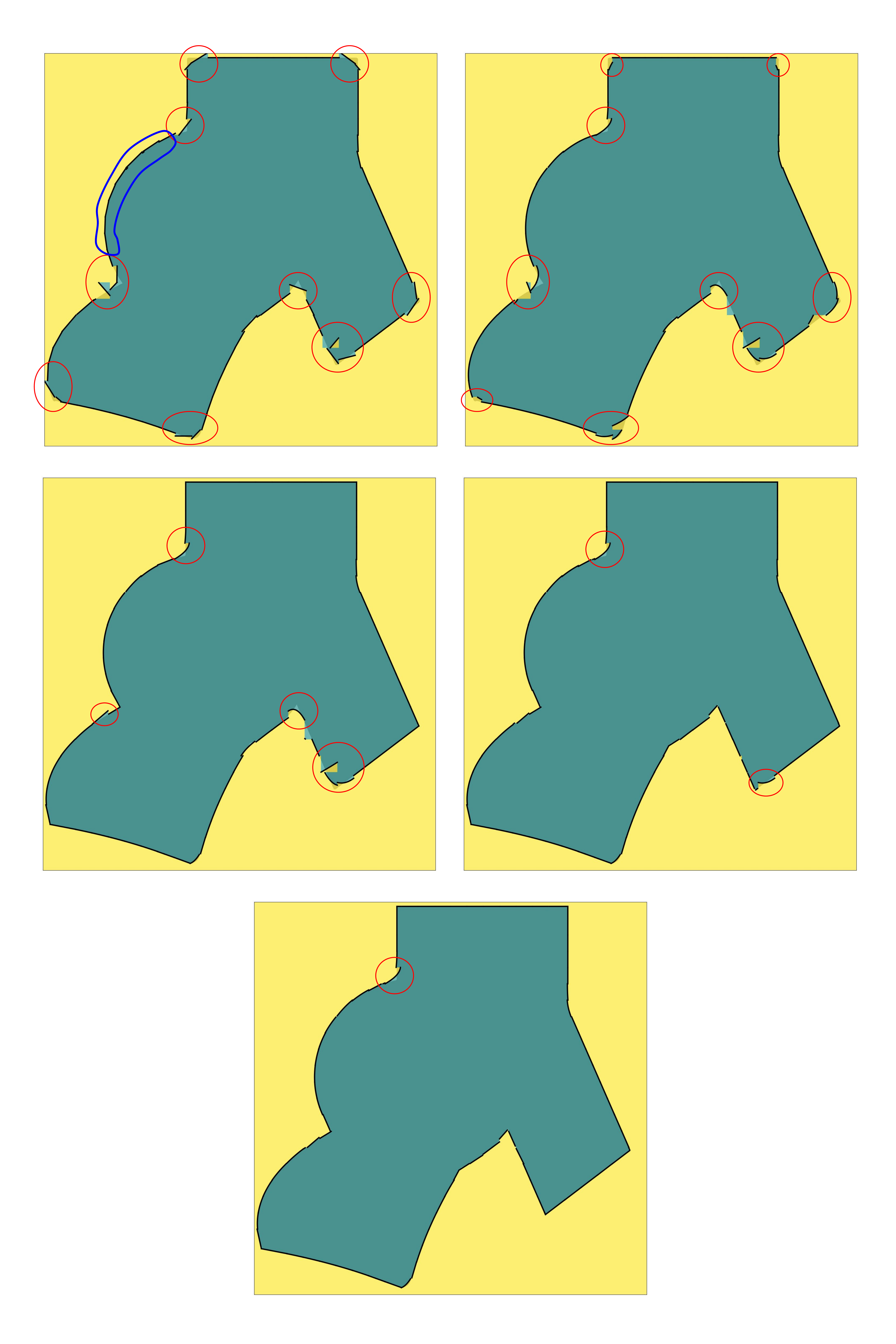}
    \caption{Reconstruction of a domain with corners 
    from cell-averages of scale $h=1/30$. 
    Reconstructions made with  
	ELVIRA-WO (up-left), 
    AEROS Quadratic (up-right), 
    AEROS Quadratic + TEM (center-left),
    AEROS Quadratic + TEM + AEROS Vertex (center-right) and finally 
    ELVIRA-WO + AEROS Quadratic + TEM + AEROS Vertex (down). 
    The red markings show the problems around vertices and how they are progressively resolved when one puts all the methods to work together.
    }
     \label{fig:cornerrec}
\end{figure}

The previously tested strategies achieve to reconstruct, at their corresponding orders, different smooth domains with H\"older smoothness $s>1$. 
This excludes the case of domains with piecewise smooth
boundaries for which the presence of corners call for a specific treatment.
The simplest option that we study here is to use local recovery by the approximation family $V_4$ of piecewise linear
interface from Example 4. We propose, in what follows, two methods to deal with vertices of any angle,  bearing in mind the limitations expressed in Figure \ref{fig:injparab}.
\newline
\newline
\textbf{AEROS Vertex}: The first method applies the AEROS strategy for $v\in \hat V_4$ a restriction of $V_4$ where $\pi/2<\theta_1<3\pi/2$ and $-\pi/2<\theta_2<\pi/2$, \textit{i.e.} elements $v$ whose interface can be written as an oriented graph, which is a particular instance of Example 3 where instead of searching $p$ in a space of univariate polynomials, we pick it into $W_4$ the space of piecewise functions with one breakpoint contained in $T$
or its immediate neighbours. This excludes the possibility of reconstructing a rectangle whose sides are parallel to the mesh but it applies to the case of the same rectangle slightly rotated. Under this restrictions it is possible, although lengthy, to extract explicit equations that allow us to derive a finite set of admissible parameters $\mu=(x_1, x_2, \theta_1, \theta_2)$ 
from the observed cell averaged vector $a_S$. We compare each proposed local approximation $v\in \hat V_4$ as in ELVIRA or OBERA, that is, by means of their associated loss $\cL(u, v)$ while retaining at the end the one achieving the minimal value between the many possibilities. This same model selection strategy can be used to aggregate other competing models, like quadratic interfaces. This has the effect of keeping higher order models when the interface is locally smooth, while taking corners into account
as illustrated on Figure \ref{fig:cornerrec}. By this approach we avoid defining a vertex detection mechanism at the expense of computational overhead as we now need to compute for each cell many losses, which was already the time bottleneck for ELVIRA method. 	
\newline
\newline	
\textbf{Tangent Extension Method (TEM)}: The above restriction that the interface with a vertex needs to be an oriented graph could be limiting in some applications but it can be removed at the expense of complexifying the reconstruction procedure. Our second proposed method deals with this aspect and consists in the following steps:
	\begin{enumerate}
	\item Associate to each singular cell $T\in \cS_h$ some reconstruction $v_T$ stemming from any of the local interface reconstruction methods discussed discussed so far.
	\item For each cell where the presence of a vertex is suspected (eventually for all $T\in \cS_h$) we search for two singular cells $T_1,T_2\in \cS_h$ satisfying the following
		\begin{itemize}
		\item $T\neq T_1\neq T_2\neq T$
		\item $S_{T_1} \cap \{T\}=\emptyset$
		\item $S_{T_2} \cap \{T\}=\emptyset$
		\end{itemize}
	where $S_{v_{T_i}}$ denotes the stencil $S$ used by a given
	smooth interface reconstruction method (for example linear or quadratic) to produce local approximations $v_{T_i}$. 	
	\item Take the parameterized interfaces $\Gamma_1$ and $\Gamma_2$, associated to $v_{T_1}$ and $v_{T_2}$ respectively, and do an order $1$ Taylor expansion at an intermediate point between cells $(T, T_1)$ and $(T, T_2)$ respectively. This yields the parameters of the two half planes $H_1$ and $H_2$ of Example 4 needed to define $v\in V_4$.
	\item Finally, we compare the new local approximation $v\in V_4$ with the existing one $v_T$, as explained above, retaining only the one whose associated loss, $\cL(u, v)$ or $\cL(u, v_T)$, is minimal.
	\end{enumerate}

	This procedure will reconstruct exactly corners when the interface is a line along both directions, but it will not produce area-consistent reconstructions on cell $T$ otherwise. In this regard, AEROS Vertex, being based on AEROS strategy, will yield interfaces that are, though not cell-consistent as one could get with OBERA, at least column-consistent as long as we remain in the interpolation case. This is ensured if the stencil width equals the number of parameters of the approximating class which in the case of $\hat V_4$ is guaranteed by using $4$-width stencils.

Figure \ref{fig:cornerrec} displays the successive improvements in the reconstruction when combining the different strategies described so far:
	\begin{itemize}
	\item \corr{On Figure \ref{fig:cornerrec} (up-left), we use the \perplexityinsert{elvira-w-oriented} method. We observe that it 
	recovers the interface in a satisfactory manner only far enough from corners.}
	\item On Figure \ref{fig:cornerrec} (up-right), we use the \perplexityinsert{quadratic-aero} method. We observe that it 
	recovers the interface in a satisfactory manner only far enough from corners.
	\item On Figure \ref{fig:cornerrec} (center-left), we first find a curve for each singular cell using \perplexityinsert{quadratic-aero} and then TEM. In this case, some of the problems are solved, in particular corners with a $90^\circ$ angle and parallel to the grid.
	\item On Figure \ref{fig:cornerrec} (center-right), we add to the previous method the first proposed approach, based on AEROS for vertices. We obtain almost perfect results except for some cells where the quadratic approximation of the interface given by \perplexityinsert{quadratic-aero} was not replaced by a better one.
	\item On Figure \ref{fig:cornerrec} (down), this last issue is addressed by aggregating, before applying any of the vertex mechanisms described before, an ELVIRA-WO strategy to offer an alternative when \perplexityinsert{quadratic-aero} is too much affected by the presence of a nearby corner.
	\end{itemize}

\FloatBarrier

\subsection{Finite volume evolution in time}
\label{sec:fvs}

Finally, we use the interface recovery methods presented so far as a constituent part of a finite volume solver with the objective of reducing numerical dissipation. We study the particular case of a simple 
linear transport PDE:
$$\frac{\partial u}{\partial t}+b\cdot \nabla u=0,$$
in the unit square domain $D=[0, 1]^2$ 
with periodic boundary conditions and initial condition 
being a piecewise constant function $u^0:D\rightarrow \{0, 1\}$ with $\Omega$ limited by a smooth interface as in Figure \ref{fig:smoothrec} or an interface having corners as in Figure \ref{fig:cornerrec}. 

As in a large class finite volume schemes, the reconstruction is used
at each step to compute the flux that updates the averages at the next 
step. For simplicity of the presentation we have set a constant (both in space and time) velocity field $b=(h/4,0)$ and worked with unit time steps $\Delta t=1$ and coarse grids of size $h=1/30$ \corr{or $h=1/60$}, so that the CFL condition is maintained. In this case, the numerical flux induced by a local
reconstruction $u_{i,j}^k$ on a cell $T$ of coordinate $(i,j)$ at time step $k$ takes the form
$$
\cF(u_{i,j}^k):=\frac 1 {|R_T|} \int_{R_T} u_{i,j}^k(x)dx
$$
where $R_T=[(i+1)h-b, (i+1)h]\times [jh, (j+1)h]$. The finite volume
approximation at the next time step $k+1$ is then given by
the updated cell-average
$$
a_{i,j}^{k+1} = a_{i,j}^{k}+\cF(\t u_{i-1,j}^k)-\cF(\t u_{i,j}^k).
$$

%

\corr{
Figure \ref{fig:scheme-corners-30} and \ref{fig:scheme-corners-60} show the effect of the presence of corners in the interface of $\Omega$ in terms of a slow but accumulative deterioration in both methods \perplexityinsert{elvira-w-oriented}, \perplexityinsert{quadratic-aero} as they are not designed to treat vertices. In contrast, \perplexityinsert{aero-qelvira-vertex} barely increase its error when the mesh is coarse: $h=1/30$. When the mesh is fine enough $(h=1/60)$ it manages to keep the error approximately on the same level for all the times tested.
}

\begin{figure}
     \centering
     \begin{subfigure}[b]{0.495\textwidth}
         \centering
         \includegraphics[width=\textwidth]{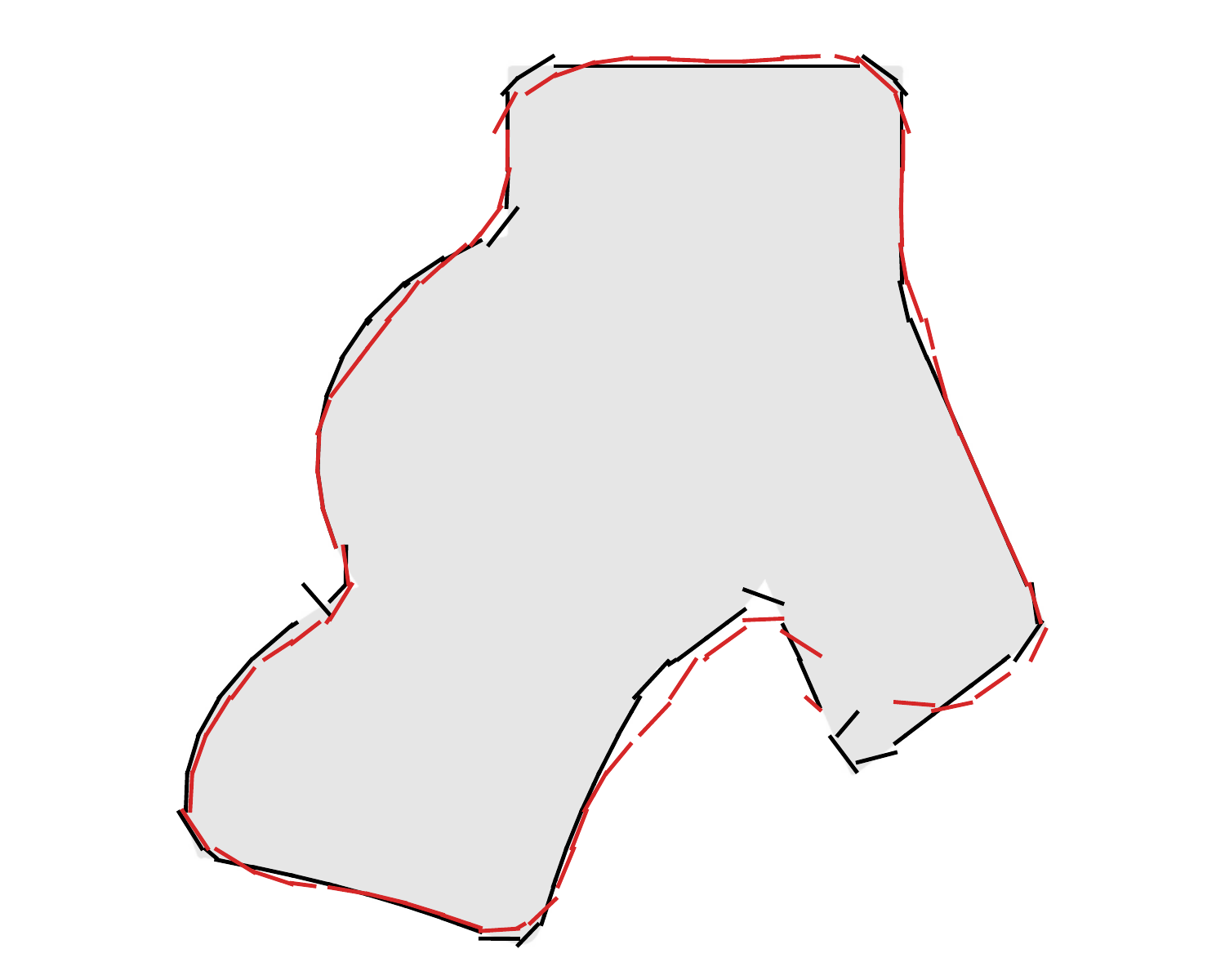}
         \caption{\corr{\perplexityinsert{elvira-w-oriented} and $h=1/30$.}}
         \label{fig:scheme-elvira}
     \end{subfigure}
     \hfill
     \begin{subfigure}[b]{0.495\textwidth}
         \centering
         \includegraphics[width=\textwidth]{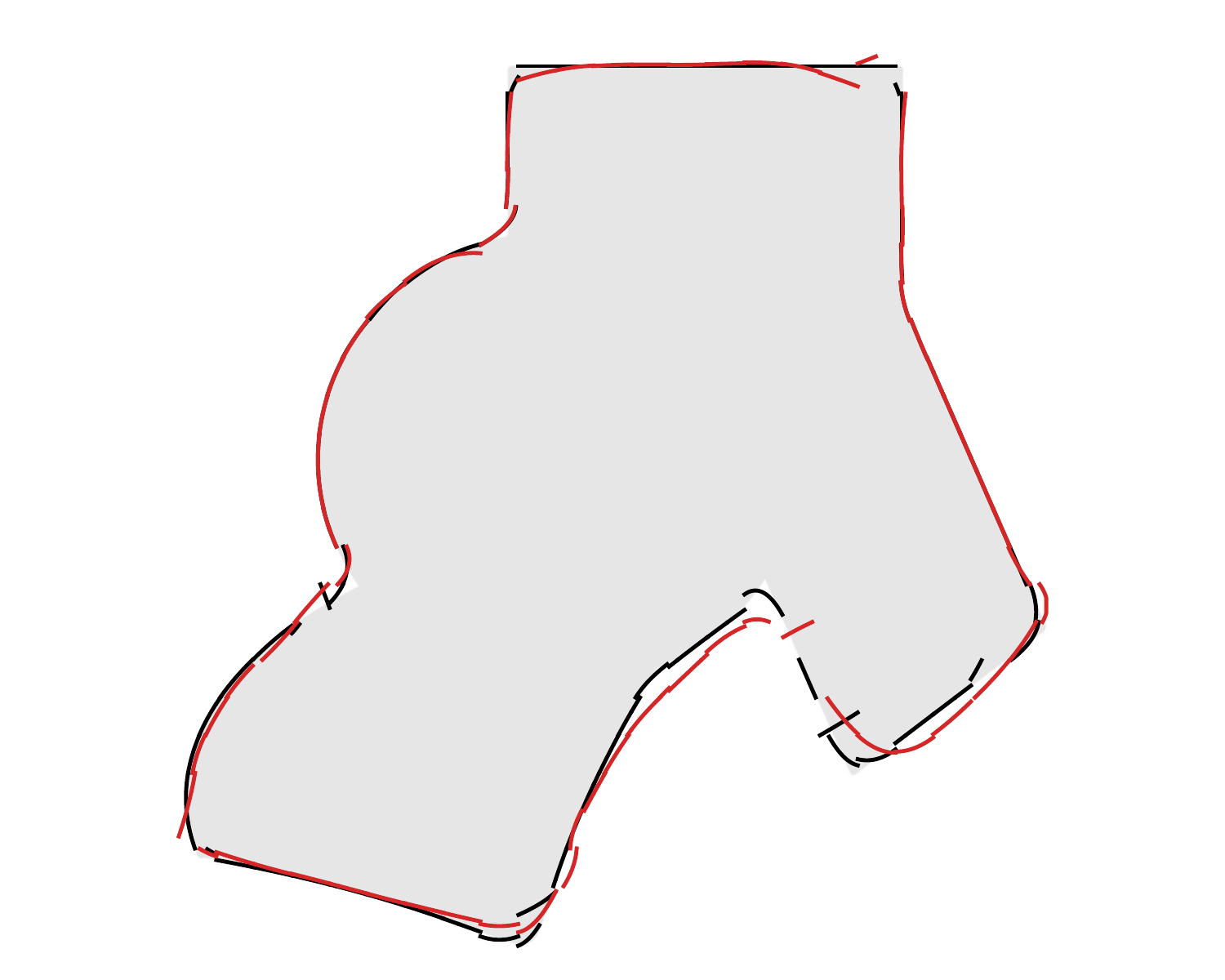}
         \caption{\corr{\perplexityinsert{quadratic-aero} and $h=1/30$.}}
         \label{fig:scheme-aeros}
     \end{subfigure}
     \vfill
     \begin{subfigure}[b]{0.495\textwidth}
         \centering
         \includegraphics[width=\textwidth]{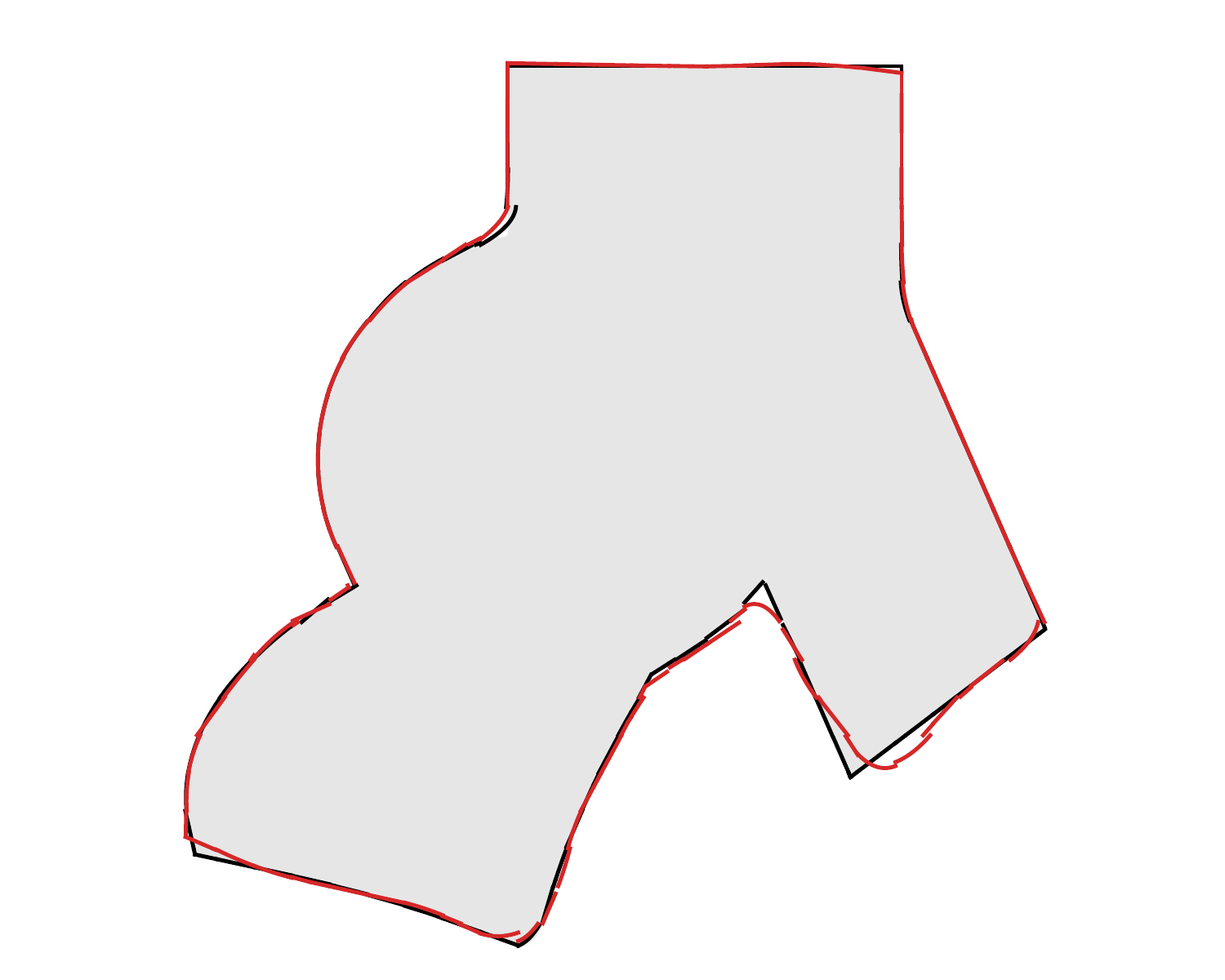}
         \caption{\corr{\perplexityinsert{aero-qelvira-vertex} and $h=1/30$.}}
         \label{fig:scheme-vertex-30}
     \end{subfigure}
     \hfill
     \begin{subfigure}[b]{0.495\textwidth}
         \centering
         \includegraphics[width=\textwidth]{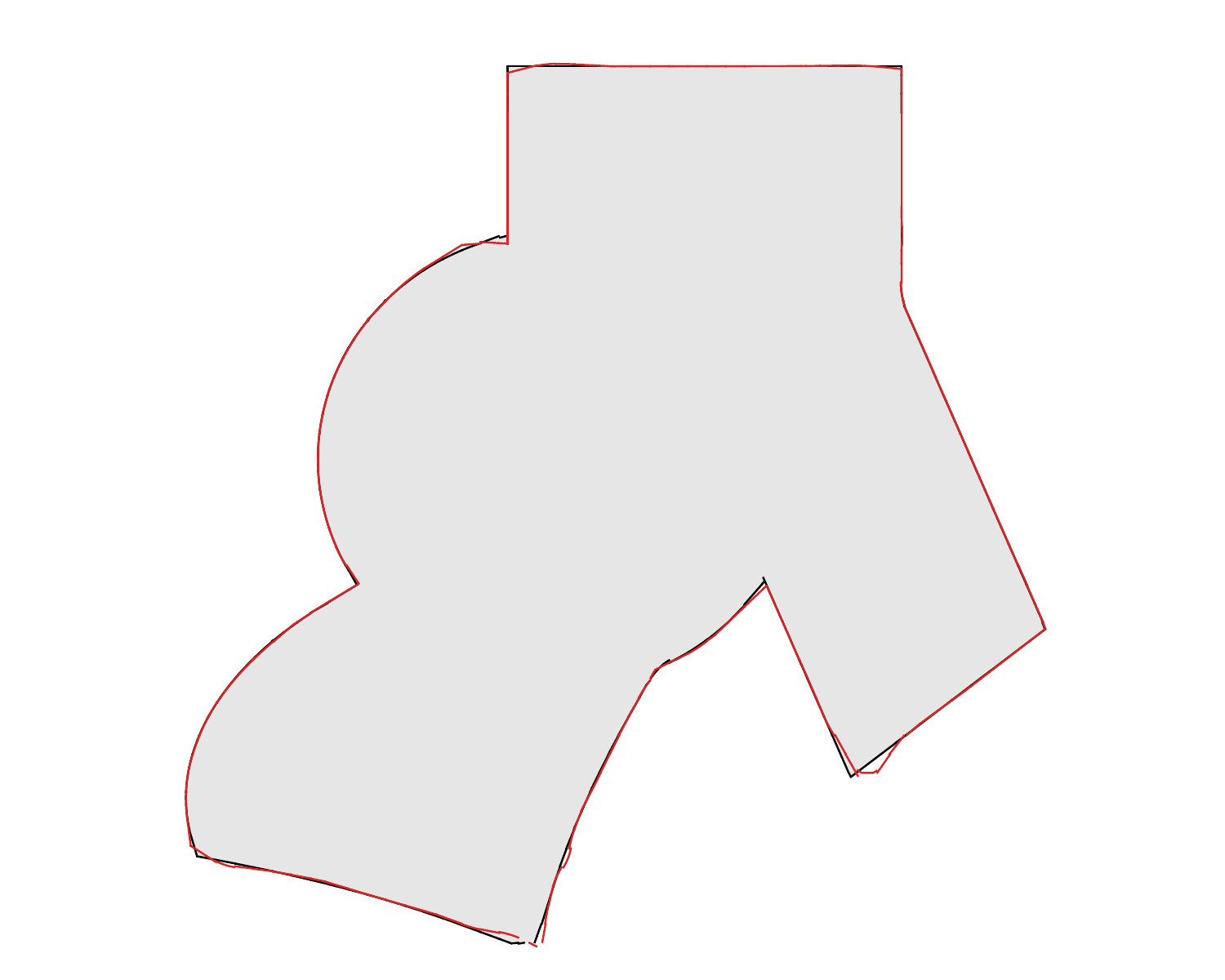}
         \caption{\corr{\perplexityinsert{aero-qelvira-vertex} and $h=1/60$.}}
         \label{fig:scheme-vertex-60}
     \end{subfigure}
	\caption{\corr{Reconstructions of a domain with corners at initial time (black) and final time (red).}}
    \label{fig:scheme-reconstruction}
\end{figure}

\corr{
In Figure \ref{fig:scheme-reconstruction} we can see the reconstruction of the interface at initial time and after transporting the shape presenting corners so that it coincides with the starting position: $120$ iterations for $h=1/30$ and $240$ iterations for $h=1/60$. We see that both \perplexityinsert{elvira-w-oriented} and \perplexityinsert{quadratic-aero} approximately keep the shape of the interface when $\partial \Omega$ is smooth with \perplexityinsert{quadratic-aero} being slightly better. However, as expected, both fail near corners. This issue is addressed by \perplexityinsert{aero-qelvira-vertex}: in Figure \ref{fig:scheme-vertex-30} the reconstruction suffers in some corners but the overall shape is maintained. In Figure \ref{fig:scheme-vertex-60} there are almost no artefacts remaining. For other comparisons see Appendix \ref{sec:schemes-rec}.
}

\begin{figure}
     \centering
     \begin{subfigure}[b]{\textwidth}
         \centering
         \includegraphics[width=\textwidth]{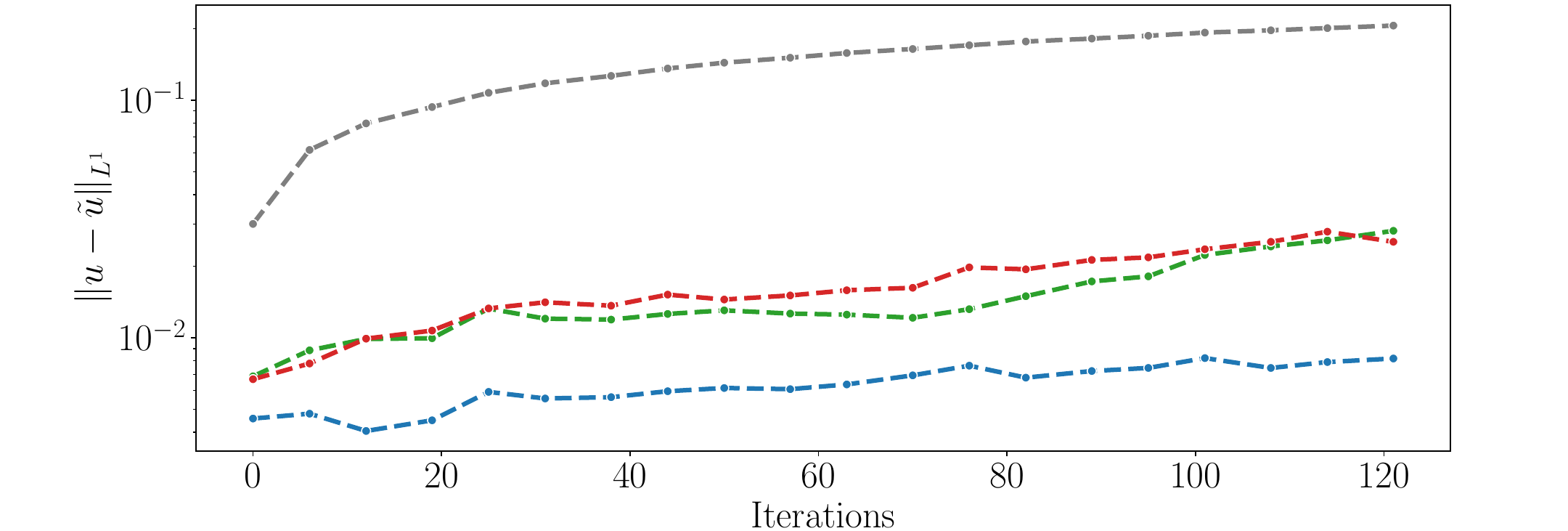}
         \caption{\corr{$h=1/30$}.}
         \label{fig:scheme-corners-30}
     \end{subfigure}
     \hfill
     \begin{subfigure}[b]{\textwidth}
         \centering
         \includegraphics[width=\textwidth]{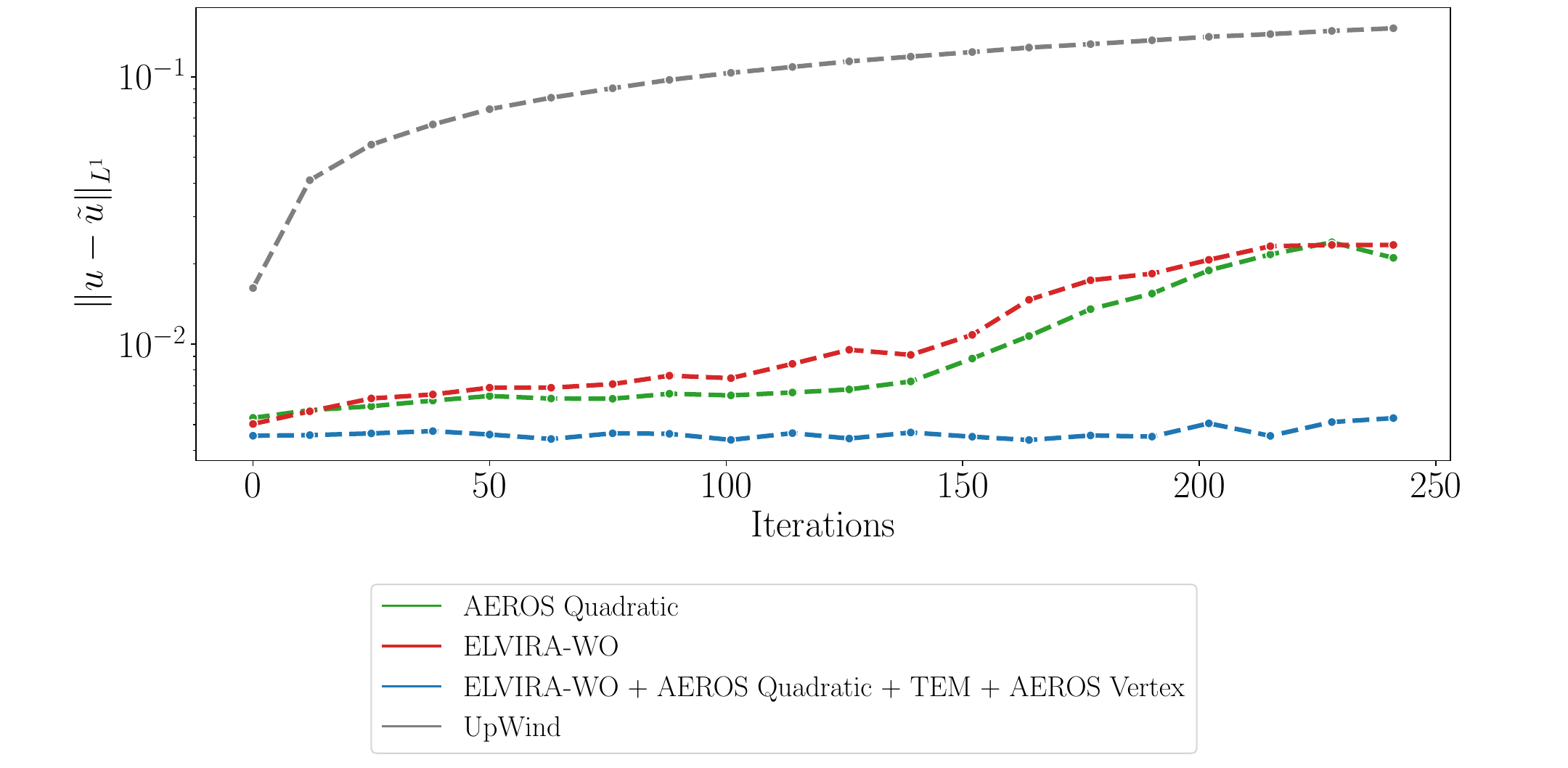}
         \caption{\corr{$h=1/60$}.}
         \label{fig:scheme-corners-60}
     \end{subfigure}
	\caption{Time evolution of the finite volume scheme $L^1$ error for a domain with corners.}
    \label{fig:scheme-error}
\end{figure}

\section{Conclusion and perspectives}
\label{sec:perspectives}

In this work, we have presented several interface recovery methods.
For the two main classes OBERA and AEROS, we have provided general
analysis strategies for establishing convergence rates that depend on the geometric smoothness of the interface. From a practical perspective, the methods can be combined with the aggregation strategy outlined in the previous section, and we have made them available in the open-source python package \footnote{\url{https://github.com/agussomacal/SubCellResolution}}.

Several natural perspectives are foreseen (and we have explored some of them already
 in the open-source package):
\begin{enumerate}
\item
Address the reconstruction of more general piecewise smooth functions with jump discontinuities across geometrically smooth or piecewise smooth interfaces. This requires a proper adaptation of the 
interface recovery strategies, combined with a high-order treatment of the smooth part of the function corresponding to the cell $T\notin \cS_h$.
The latter can be done by using
polynomial reconstructions on stencils not containing cells of $\cS_h$ following the standard ENO strategy.
\item
Study the use of machine learning techniques trained on 
sufficiently rich sets of interfaces for performing certain tasks
in an automated and hopefully more efficient manner. Such tasks
include, for example, the fast reconstruction of the
parameter $\mu$ from the cell averages, the identification of 
cells that may contain vertices, or the direct access to the numerical flux
in the case of finite volume scheme, as proposed for example in
\cite{Despres2020} for vertices forming angles of $\phi=90^{\circ}$. 
It should however be noted that, as opposed to the approaches
that we developed in this paper, the machine learning-based
approach does not offer rigorous convergence guarantees.
Also, our attempts to beat the AEROS method with machine learning strategies were so far unsuccessful, 
both from the accuracy point of view,
and the runtime point of view. We provide our implementation
of these strategies in the Python package.
\end{enumerate}

\appendix
\section{The orientation test}

In this apprendix we give the proof of Theorem 
\ref{theosobel}, which is based on:
\begin{enumerate}
\item First studying the case where the interface $\Gamma$ is a line over the $3\times 3$ stencil where
the numerical gradient $G_T=(H_T,V_T)$ is computed.
\item Second applying a perturbation argument in the case of a general $\cC^s$
interface which locally deviates from a line
in a quantitatively controlled manner.
\end{enumerate}

\subsection{The case of a linear interface}

We assume here that, over the $3\times 3$ stencil $S$
centered at $T$, the interface $\Gamma$ is a line
crossing $T$. Therefore, the restriction of $u$
to $S$ is of the form
$$
u_{|S}=v_{r,\theta}:=\Chi_{\{\<z-z_T,e_\theta^\perp\>\leq r\} }.
$$
where $z_T$ is the center of $T$
and $e_\theta=(-\sin(\theta),\cos(\theta))$ with $\theta$ the
angle between $\Gamma$ and the horizontal line, that is, $e_\theta$ is the unit normal
vector to $\Gamma$ pointing to the outward direction where $u_{|S}=0$.
The following result shows that the orientation test discriminates exactly
if the direction of $\Gamma$ is closer to horizontal or vertical.
Its proof uses elementary geometrical arguments, 
which are only sketched using pictures in order to avoid cumbersome
analytic developments.

\begin{theorem}
If $G_T=(H_T,V_T)$ is the numerical gradient based on the Sobel
filter for the above function $v_{r,\theta}$, then the following holds:
\begin{itemize}
\item $|V_T|>|H_T|$ if and only if $\theta\in [0,\pi/4[\cup ]3\pi/4,5\pi/4[\cup ]7\pi/4,2\pi[$ 
\item $|H_T|>|V_T|$ if and only if $\theta\in ]\pi/4,3\pi/4[\cup ]5\pi/4,7\pi/4[$ 
\item $|H_T|=|V_T|$ if and only if $\theta\in \{\pi/4,3\pi/4,5\pi/4,7\pi/4\}$ 
\end{itemize}
In addition
\begin{itemize}
\item
$V_T>0$ if and only if $\theta \in ]\pi/2,3\pi/2[$
and $V_T< 0$ if and only if $\theta\in [0,\pi/2[\cup ]3\pi/2,2\pi[$
\item
$H_T>0$ if and only if $\theta\in ]0,\pi[$ and 
$H_T<0$ if and only if $\theta \in ]\pi,2\pi[$.
\end{itemize}
\label{theorientline}
\end{theorem}

\noindent
{\bf Proof:} Without loss of generality, we only consider the case 
where $\theta\in [0,\pi/4]$ since all other cases $[k\pi/4,(k+1)\pi/4]$ 
for $k=1,\dots,7$ are treated in a similar way. In order to understand
the effect of $\theta$ on the values of $H_T$ and $V_T$, we parametrize
the function $v_{r,\theta}$ differently: we fix $\o z_T$ to be
the point crossed  by $\Gamma$ on the descending diagonal of $T$
(which exists and is unique when $\theta\in [0,\pi/4]$) and study $H_T$ and $V_T$ for the function
$$
v_{\theta}:=\Chi_{\{\<z-\o z_T,e_\theta^\perp\>\leq 0\} },
$$
as we let $\theta$ vary.  By scale invariance, we may assume that
we work with cells of side-length equal to $1$ without affecting $H_T$
and $V_T$.

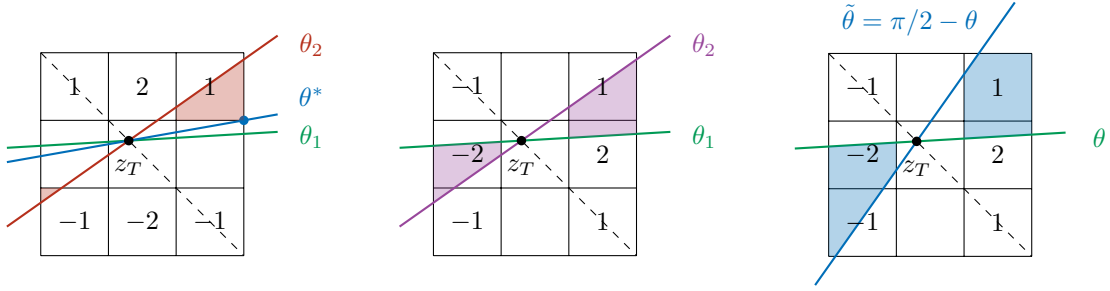
\begin{figure}[ht]
\begin{center}
\begin{minipage}{.3\textwidth}
\centering
\begin{tikzpicture}[scale=0.9]
  \draw[step=1cm] (0,0) grid (3,3);
  \draw[dashed](0,3)--(3,0); 
  
  \draw[solid, ForestGreen, thick] plot[smooth,domain=-0.5:3.5] (\x, {1/17*(\x-1.3)+1.7}); 
  \draw[solid, NavyBlue, thick] plot[smooth,domain=-0.5:3.5] (\x, {(3/17)*(\x-1.3)+1.7}); 
  \fill[NavyBlue] (3,2) circle (0.07cm);
  \draw[solid, BrickRed, thick] plot[smooth,domain=-0.5:3.5] (\x, {(12/17)*(\x-1.3)+1.7}); 
  
  \fill [fill=BrickRed, fill opacity=0.3] (2, 2)--(2, 2.1941176470588237)--(3, 2.9)--(3, 2)--cycle;
  \fill [fill=BrickRed, fill opacity=0.3] (0, 1)--(0.31, 1)--(0, 0.7823529411764705)--cycle;
    
  \node[ForestGreen] at (4, 1.7+1/17) {$\theta_1$};
  \node[NavyBlue] at (4, 1.7+11/17) {$\theta^*$};
  \node[BrickRed] at (4, 1.7+24/17) {$\theta_2$};
  
  \node at (1.3, 1.3) {$z_T$};
  \fill (1.3,1.7) circle (0.07cm);
  
  \node at (0.5, 0.5) {$-1$};
  \node at (1.5, 0.5) {$-2$};
  \node at (2.5, 0.5) {$-1$};
  \node at (0.5, 2.5) {$1$};
  \node at (1.5, 2.5) {$2$};
  \node at (2.5, 2.5) {$1$};
\end{tikzpicture}
\end{minipage}
\begin{minipage}{.3\textwidth}
\centering
\begin{tikzpicture}[scale=0.9]
  \draw[step=1cm] (0,0) grid (3,3);
  \draw[dashed](0,3)--(3,0); 
  
  \draw[solid, ForestGreen, thick] plot[smooth,domain=-0.5:3.5] (\x, {1/17*(\x-1.3)+1.7}); 
  \draw[solid, Purple, thick] plot[smooth,domain=-0.5:3.5] (\x, {(12/17)*(\x-1.3)+1.7}); 
  
  \fill [fill=Purple, fill opacity=0.3] (2, 1.7411764705882353)--(3, 1.8)--(3, 2.9)--(2, 2.1941176470588237)--cycle;
  \fill [fill=Purple, fill opacity=0.3] (0, 0.7823529411764705)--(1, 1.4882352941176469)--(1, 1.6823529411764706)--(0, 1.6235294117647059)--cycle;
    
  \node[ForestGreen] at (4, 1.7+1/17) {$\theta_1$};
  \node[Purple] at (4, 1.7+24/17) {$\theta_2$};
  
  \node at (1.3, 1.3) {$z_T$};
  \fill (1.3,1.7) circle (0.07cm);
  
  \node at (0.5, 0.5) {$-1$};
  \node at (0.5, 1.5) {$-2$};
  \node at (0.5, 2.5) {$-1$};
  \node at (2.5, 0.5) {$1$};
  \node at (2.5, 1.5) {$2$};
  \node at (2.5, 2.5) {$1$};
\end{tikzpicture}
\end{minipage}
\begin{minipage}{.3\textwidth}
\centering
\begin{tikzpicture}[scale=0.9]
  \draw[step=1cm] (0,0) grid (3,3);
  \draw[dashed](0,3)--(3,0); 
  
  \draw[solid, ForestGreen, thick] plot[smooth,domain=-0.5:3.5] (\x, {1/17*(\x-1.3)+1.7}); 
  \draw[solid, NavyBlue, thick] plot[smooth,domain=-0.2:2.75] (\x, {(17/12)*(\x-1.3)+1.7}); 
  
  \fill [fill=NavyBlue, fill opacity=0.3] (2, 2.6916666666666664)--(2.2176470588235295, 3)--(3, 3)--(3, 1.8)--(2, 1.7411764705882353)--cycle;
    \fill [fill=NavyBlue, fill opacity=0.3] (0, 0)--(0.1, 0)--(1, 1.275)--(1, 1.6823529411764706)--(0, 1.6235294117647059)--cycle;
    
  \node[ForestGreen] at (4, 1.7+1/17) {$\theta$};
  \node[NavyBlue] at (1.2, 3.5) {$\t \theta=\pi/2-\theta$};
  
  \node at (1.3, 1.3) {$z_T$};
  \fill (1.3,1.7) circle (0.07cm);
  
  \node at (0.5, 0.5) {$-1$};
  \node at (0.5, 1.5) {$-2$};
  \node at (0.5, 2.5) {$-1$};
  \node at (2.5, 0.5) {$1$};
  \node at (2.5, 1.5) {$2$};
  \node at (2.5, 2.5) {$1$};
\end{tikzpicture}
\end{minipage}
\end{center}
\caption{Dependence of $V_T$ (left), $H_T$ (center), $|V_T|-|H_T|$ (right)
as $\theta$ varie in $[0,\pi/4]$.}
\label{figVH}
\end{figure}

Figure \ref{figVH} (left) pictures the value of $V_T$ as 
the difference between areas of the portions of 
cells from the upper and lower rows crossed by the half-plane
below $\Gamma$ with weight $2$ for central cells
and $1$ for left and right cells. This difference is strictly
negative for all $\theta\in [0,\pi/4]$. Its value at $\theta=0$
is equal $-4$. As $\theta$ grows towards $\pi/4$ it first stays equal to 
$-4$ until it starts strictly increasing for some value $\theta^*\in [0,\pi/4[$
that depends on the position of $\o z_T$ on the diagonal. This monotonic growth can be checked
by observing that for $0\leq \theta_1<\theta_2\leq\pi/4$, one
has $V_T(v_{\theta_2})-V_T(v_{\theta_1})=V_T(v_{\theta_2}-v_{\theta_1})$
and the function $v_{\theta_2}-v_{\theta_1}$ is supported in a symmetric cone 
$K_{\theta_1,\theta_2}$ centered
at $\o z_T$ and has value $1$ on the right and $-1$ on the left. 
Thus $V_T(v_{\theta_2})-V_T(v_{\theta_1})$
is the sum of the areas of the portions of 
cells from the upper and lower rows intersected by $K_{\theta_1,\theta_2}$
with weight $2$ for central cells and $1$ for for left and right cells,
which is strictly positive if $\theta_2>\theta^*$.

Figure \ref{figVH} (center) pictures the value of $H_T$ as 
the difference between areas of the portions of 
cells from the right and left columns crossed by the half-plane
below $\Gamma$ with weight $2$ for central cells
and $1$ for lower and upper cells. This difference is null when $\theta=0$
and increases strictly as $\theta$ grows from $0$ to $\pi/4$.
Once again, the strictly monotonic growth is due to the fact
that $V_T(v_{\theta_2})-V_T(v_{\theta_1})$
is the sum of the areas of the portions of 
cells from the left and right columns crossed by $K_{\theta_1,\theta_2}$
with weight $2$ for central cells and $1$ for for left and right cells,
which is strictly positive whenever $0\leq \theta_1<\theta_2\leq\pi/4$.

This demonstrates the second statement regarding signs of $V_T$ and $H_T$.

On the other hand, by symmetry, we note that $V_T(v_{\theta})=-H_T(v_{\t \theta})$
where $\t \theta=\pi/2-\theta$ with same base point $z_T$. Therefore 
$|V_T|-|H_T|=-V_T-H_T=H_T(v_{\t \theta}-v_{\theta})$
and $v_{\t \theta}-v_{\theta}$ is supported in a symmetric cone $K_{\theta,\t \theta}$ centered
at $\o z_T$ with value $1$ on the right and $-1$ on the left.
As pictured on Figure \ref{figVH} (right), this quantity
is the sum of the areas of the portions of 
cells from the right and left column crossed by $K_{\theta,\t \theta}$
with weight $2$ for central cells and $1$ for lower and upper cells.
These quantities are null when $\theta=\pi/4$ since the cone is then
restricted to a line, and increases strictly 
as $\theta$ decreases from $\theta/4$ to $0$ since
the cone is opening.

This demonstrates the first statement regarding comparison
between $|V_T|$ and $|H_T|$.
\hfill $\Box$.
\newline

We will use the following direct consequence of this result, which is obtained
by compactness since $V_T$ and $H_T$ are continuous with respect to $\theta$
and $\o z_T$: for any $0<\delta<\pi/4$, there exists a 
$\gamma=\gamma(\delta)>0$ such that 
\be
\theta \in [\pi/4+\delta,3\pi/4-\delta]\cup  [5\pi/4+\delta,7\pi/4-\delta]
\implies |H_T|\geq |V_T|+\gamma
\label{strict1}
\ee
and
\be
\theta \in [3\pi/4-\delta,5\pi/4+\delta]
\implies V_T\geq \gamma,
\label{strict2}
\ee

\subsection{A perturbation analysis}

We next turn to the proof of Theorem 
\ref{theosobel}, focusing without loss of generality on case $1$.
Let us assume that the interface $\Gamma=\partial \Omega$ has
$\cC^s$ smoothness for some $s>1$, and let $T$
be a singular cell crossed by $\Gamma$ and 
$S$ the $3\times 3$ stencil centered at $T$.

\begin{figure}[ht]
\begin{center}
\begin{minipage}{\textwidth}
\centering
\begin{tikzpicture}[scale=0.9]
  \draw[step=1cm] (0,0) grid (3,3);
  
  \draw[solid, ForestGreen, thick] plot[smooth,domain=-0.5:3.5] (\x, {4/17*(\x-1.3)+1.7}); 
  \draw[dashed, ForestGreen] plot[smooth,domain=-3.1:6.1] (\x, {4/17*(\x-1.3)+1.7}); 
  \draw[->, ForestGreen, thick] (1.3, 1.7)--(1.1855, 2.18); 
  \node[ForestGreen] at (0.7855, 2.18) {$e_\theta$}; 
  \draw[solid, NavyBlue, thick] plot[smooth,domain=-0.5:3.5] (\x, {(4/17)*(\x-1.3)+1.7-0.15*(\x-1.3)^2+0.05*(\x-1.3)^3}); 
  \draw[dashed, NavyBlue] plot[smooth,domain=-1:5] (\x, {(4/17)*(\x-1.3)+1.7-0.15*(\x-1.3)^2+0.05*(\x-1.3)^3}); 
  
  
  \node[ForestGreen] at (-2.9, 1.15) {$L$}; 
  \node[ForestGreen] at (-1.6, 1.65) {$v_\theta=0$}; 
  \node[ForestGreen] at (-1.6, 0.25) {$v_\theta=1$}; 
  
  \node[NavyBlue] at (4.75, 2.25) {$\Gamma$}; 
  \node[NavyBlue] at (4, 3) {$u=0$}; 
  \node[NavyBlue] at (4, 1.5) {$u=1$}; 
  
  \node at (1.3, 1.3) {$\t z_T$};
  \fill (1.3,1.7) circle (0.07cm);
\end{tikzpicture}
\end{minipage}
\end{center}
\caption{Approximation of a smooth interface by its tangent}
\label{pertu}
\end{figure}
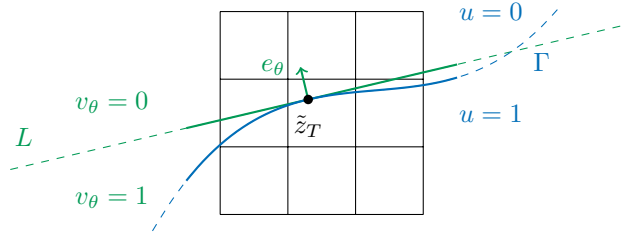

We now fix $\t z_T$ to be one point of $\Gamma\cap T$ and $e_\theta^\perp=(-\sin(\theta),\cos(\theta))$ be
the outer unit normal to $\Omega$ at $\t z_T$. The function
$$
v_{\theta}:=\Chi_{\{\<z-\t z_T,e_\theta^\perp\>\leq 0\} },
$$
is a perturbation of $u$ as pictured on Figure \ref{pertu}
where the line interface $L$ is the tangent to $\Omega$ at $\t z_T$.
Since $\Omega$ has $\cC^s$ smoothness, the deviation between
the curved interface $\Gamma$ and its tangent $L$
has area of order $\cO(h^{s+1})$ over $S$, and therefore
$$
|a_{\t T}(u)-a_{\t T}(v_{\theta})|\leq Ch^{s-1}, \quad \t T\in S,
$$
for some fixed constant $C$ that only depends on the $\cC^s$ norm
of the graph that locally characterizes $\Gamma$. In turn, up to multiplying
the constant $C$ by $8$, one has
\be
|H_T(u)-H_T(v_{\theta})|\leq Ch^{s-1} \quad {\rm and}
\quad |V_T(u)-V_T(v_{\theta})|\leq Ch^{s-1}
\label{pert}
\ee

\medskip

We now take any $0<\delta<\pi/4$ arbitrarily small and consider
the quantity $\gamma=\gamma(\delta)$ such that \iref{strict1} and \iref{strict2} 
are valid. For $h\leq h_0$ small enough, we are ensured that
$$
Ch^{s-1} \leq \frac \gamma 3,
$$
where $C$ is the constant in \iref{pert}. Therefore, if 
$\theta\notin [0,\pi/4+\delta] \cup [7\pi/4-\delta,2\pi[$,
or equivalently $\theta\notin [-\pi/4-\delta,\pi/4+\delta]$,
we obtain either by \iref{strict1} that
$$
|H_T(u)|\geq |H_T(v_{\theta})|-\frac \gamma 3 \geq
|V_T(v_{\theta}|+\frac {2\gamma} 3 \geq |V_T(u)|+
\frac {\gamma} 3>|V_T(u)|,
$$
or by \iref{strict2} that
$$
V_T(u)\geq V_T(v_{\theta})-\frac \gamma 3 
\geq \frac {2\gamma} 3>0.
$$
Therefore, we conclude for case $1$ that if
$|H_T(u)|\leq |V_T(u)|$ and $V_T(u)\leq 0$, the angle $\theta$
of the tangent line $L$ necessarily lies in $[-\pi/4-\delta,\pi/4+\delta]$.
In other words, the points $z=(x,y)$ in the half plane $\{\<z-\t z_T,e_\theta^\perp\>\leq 0\}$ can 
be characterized by an equation of the form
$$
y\leq a(x),
$$
where $a$ is affine and $|a'(x)|\leq 1+\e$ such that $1+\e=\tan(\pi/4+\delta)$,
with $0\leq \e\leq \delta$ and $\e\sim \delta$ when $\delta$ is small.
On the other hand the interface $\Gamma$ is described on $T$ by an
equation of the form 
$$
y\leq \psi(x),
$$
where due to the $\cC^s$ smoothness of $\Gamma$, one has an estimate
of the form
$$
|\psi'(x)-a'(x)|\leq Ch^{s-1},
$$
and the same will hold on a stencil $S_T$ of width $2k+1$ up to enlarging the value of $C$, so
that
$$
|\psi'(x)|\leq 1+\e+Ch^{s-1}, \quad x\in I,
$$
where $I$ is the horizontal support of $S_T$. 
Therefore, we can find $h^*=h^*(\Omega)$ such that if $h\leq h^*$ we are thus ensured that
$$
|\psi'(x)| \leq \frac {k+2}{k+1}, \quad x\in I.
$$
This implies that the graph of $\psi$ remains
confined in $S_T$ if we choose it to be of size
$(2k+1)\times (2l+1)$ with $l=k+2$, which concludes 
the proof of the theorem.

\appendix
\section{Finite volumes reconstructions}
\label{sec:schemes-rec}

\FloatBarrier

\begin{figure}
     \centering
     \begin{subfigure}[b]{\textwidth}
         \centering
         \includegraphics[width=\textwidth]{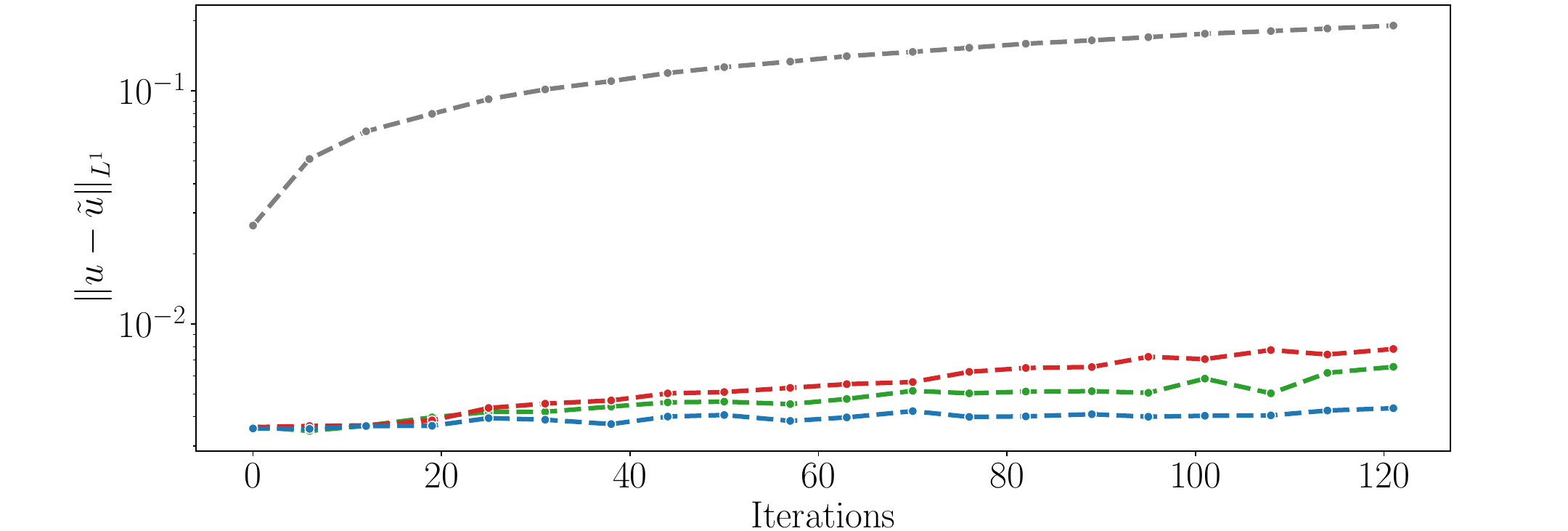}
         \caption{\corr{$h=1/30$}.}
         \label{fig:scheme-corners-30}
     \end{subfigure}
     \hfill
     \begin{subfigure}[b]{\textwidth}
         \centering
         \includegraphics[width=\textwidth]{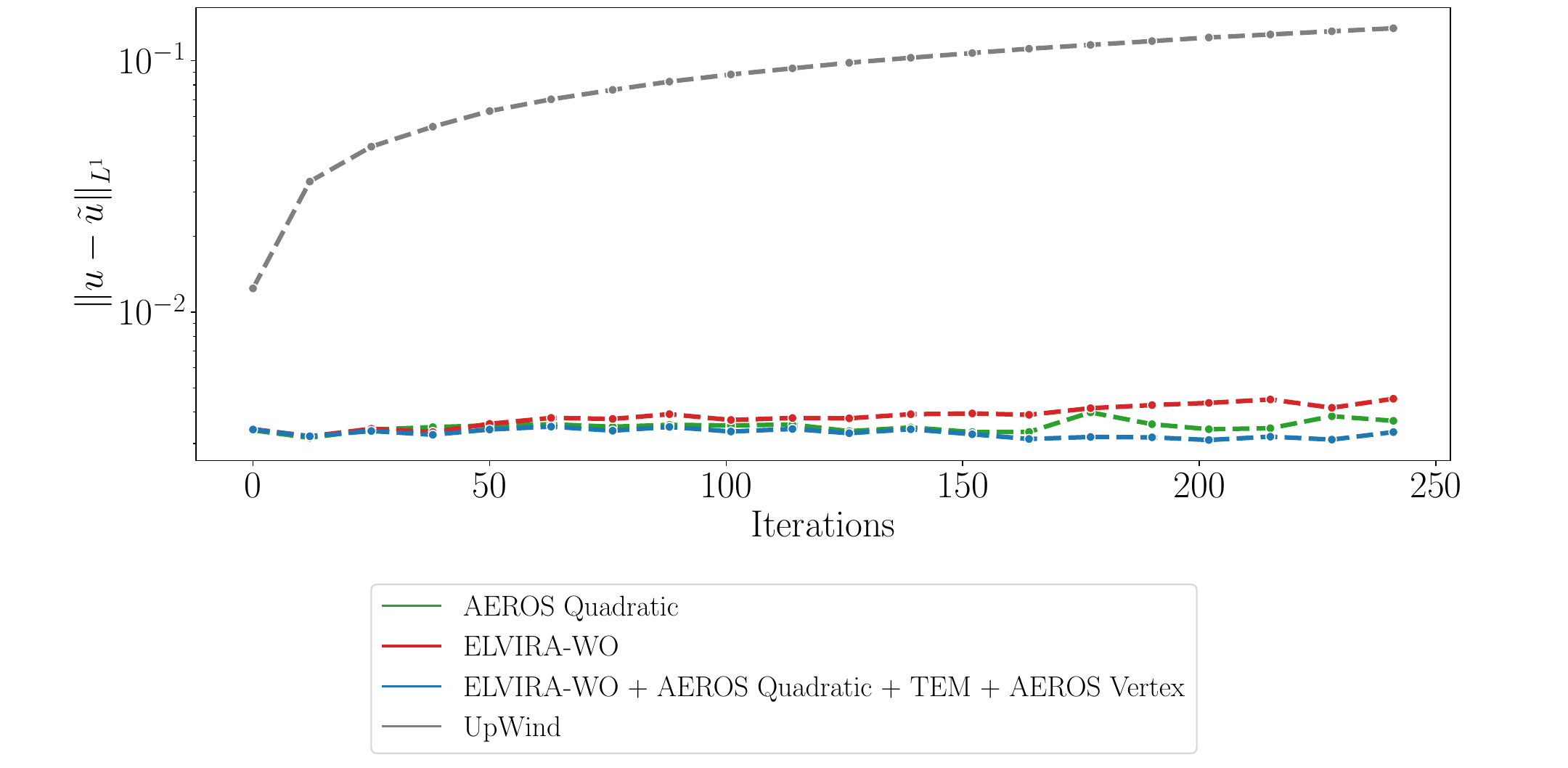}
         \caption{\corr{$h=1/60$}.}
         \label{fig:scheme-corners-60}
     \end{subfigure}
	\caption{Time evolution of the finite volume scheme $L^1$ error for a smooth domain.}
    \label{fig:scheme-error-batata}
\end{figure}

\begin{figure}
     \centering
     \begin{subfigure}[t]{0.32\textwidth}
         \centering
         \includegraphics[width=\textwidth]{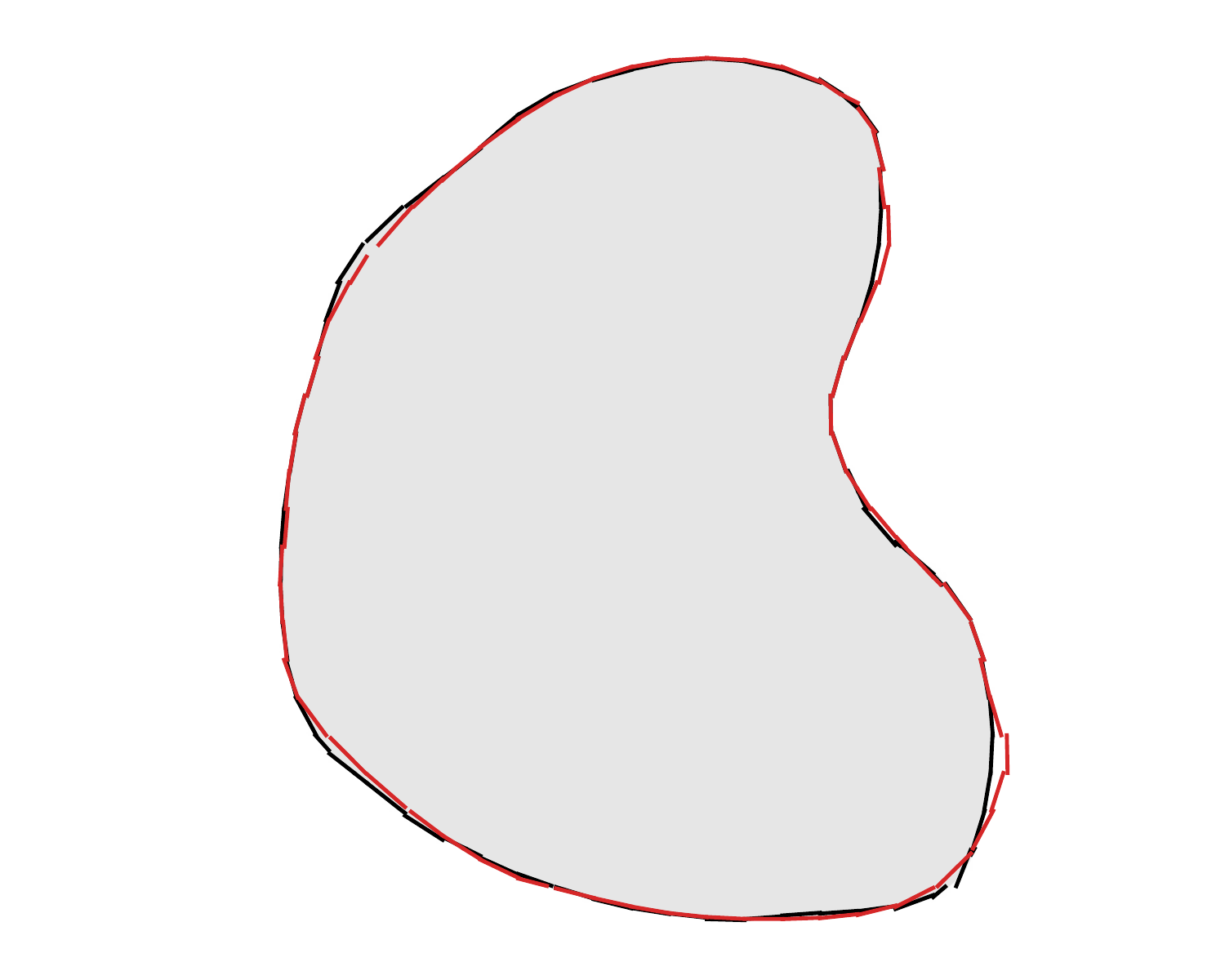}
         \caption{\corr{\perplexityinsert{elvira-w-oriented} and $h=1/30$.}}
     \end{subfigure}
     \hfill
     \begin{subfigure}[t]{0.32\textwidth}
         \centering
         \includegraphics[width=\textwidth]{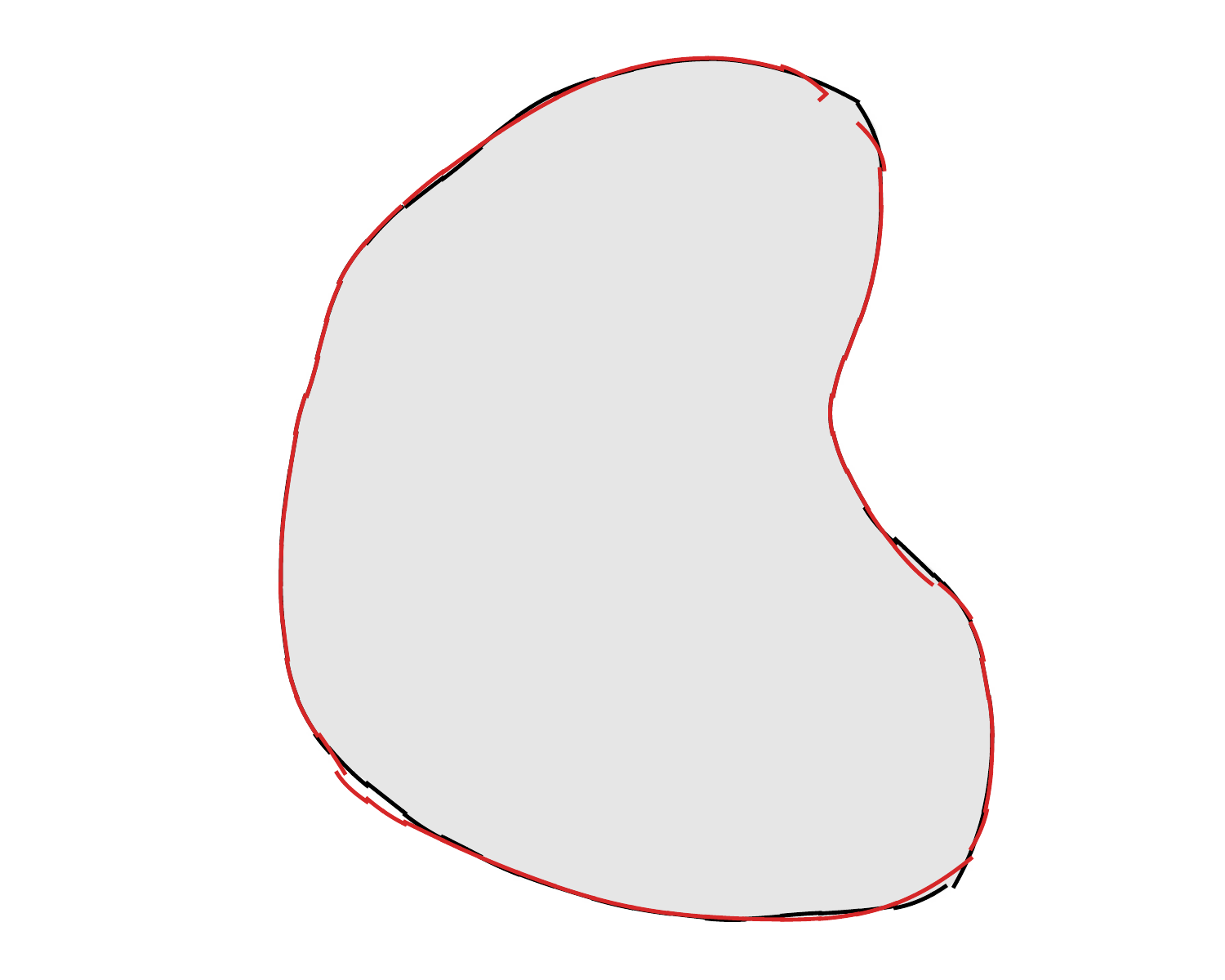}
         \caption{\corr{\perplexityinsert{quadratic-aero} and $h=1/30$.}}
     \end{subfigure}
     \hfill
     \begin{subfigure}[t]{0.32\textwidth}
         \centering
         \includegraphics[width=\textwidth]{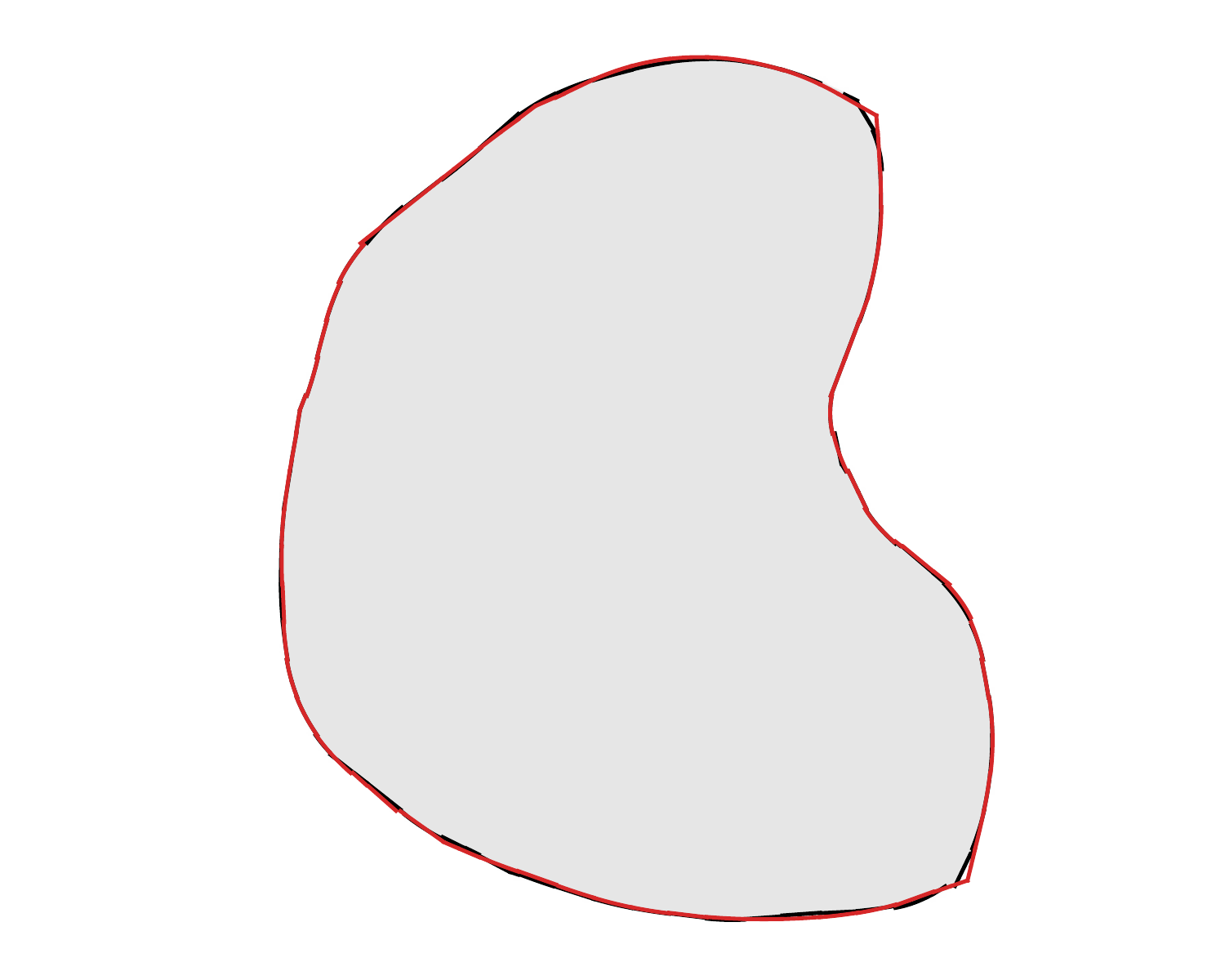}
         \caption{\corr{\perplexityinsert{aero-qelvira-vertex} and $h=1/30$.}}
     \end{subfigure}
     \vfill
          \begin{subfigure}[t]{0.32\textwidth}
         \centering
         \includegraphics[width=\textwidth]{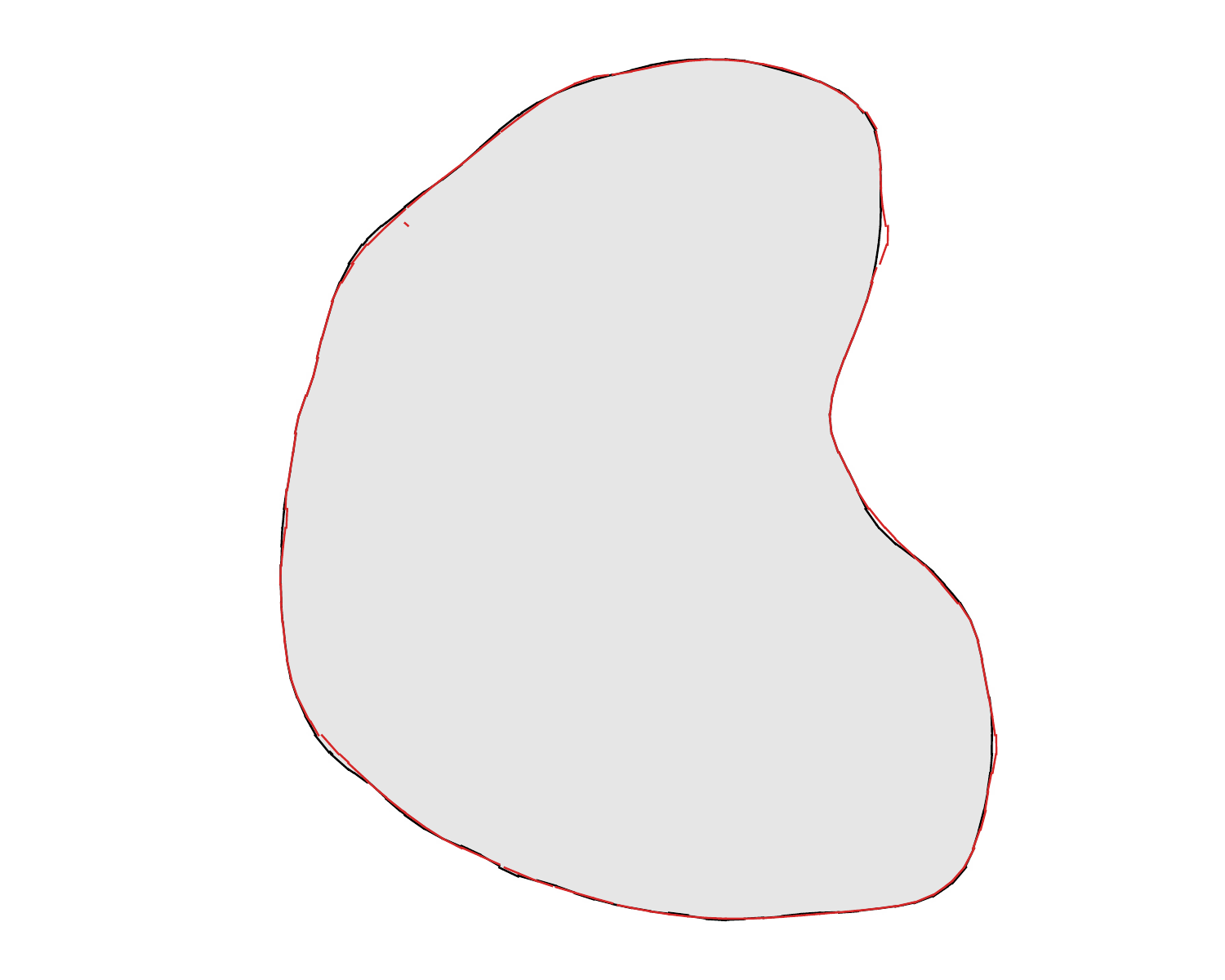}
         \caption{\corr{\perplexityinsert{elvira-w-oriented} and $h=1/60$.}}
     \end{subfigure}
     \hfill
     \begin{subfigure}[t]{0.32\textwidth}
         \centering
         \includegraphics[width=\textwidth]{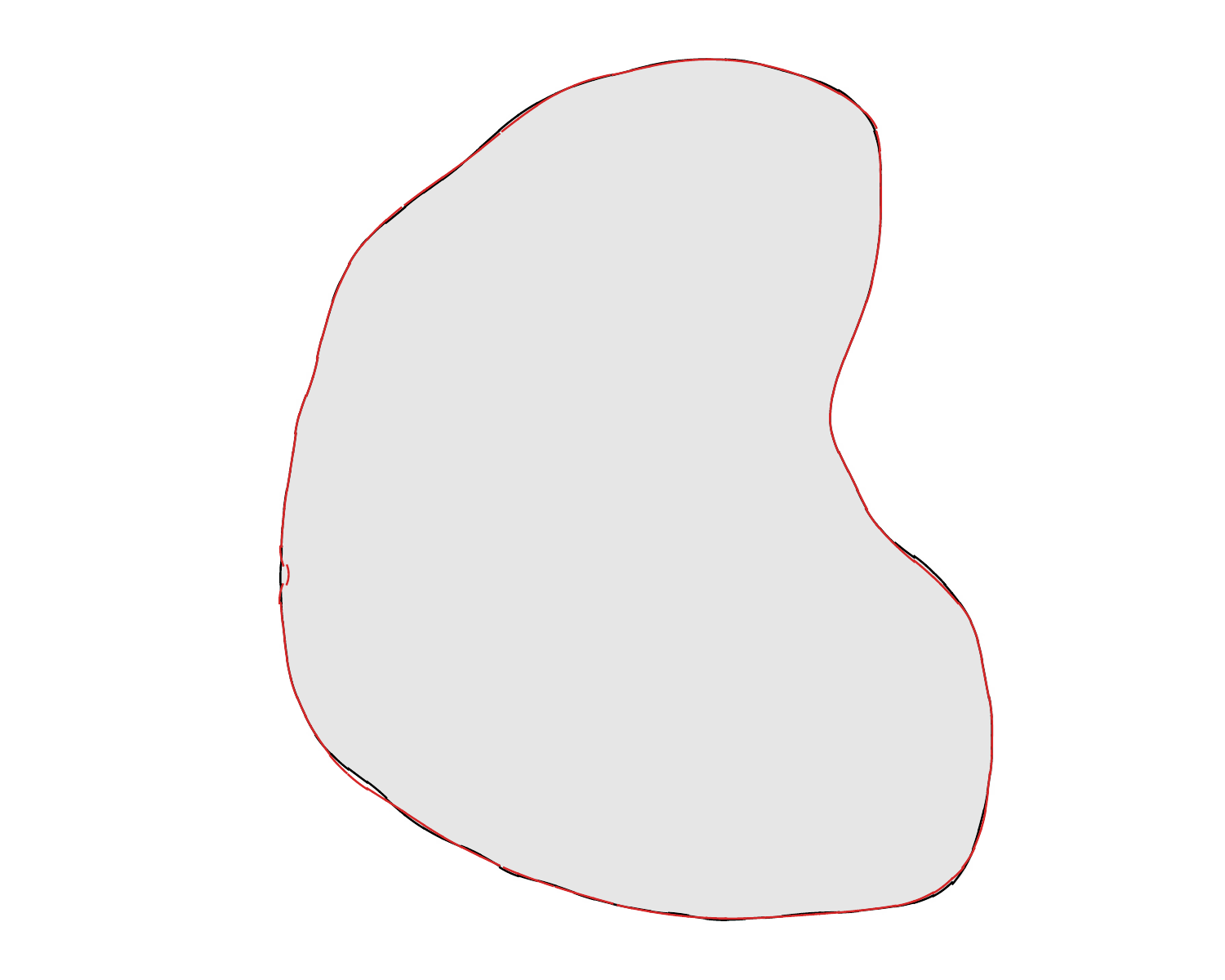}
         \caption{\corr{\perplexityinsert{quadratic-aero} and $h=1/60$.}}
     \end{subfigure}
     \hfill
     \begin{subfigure}[t]{0.32\textwidth}
         \centering
         \includegraphics[width=\textwidth]{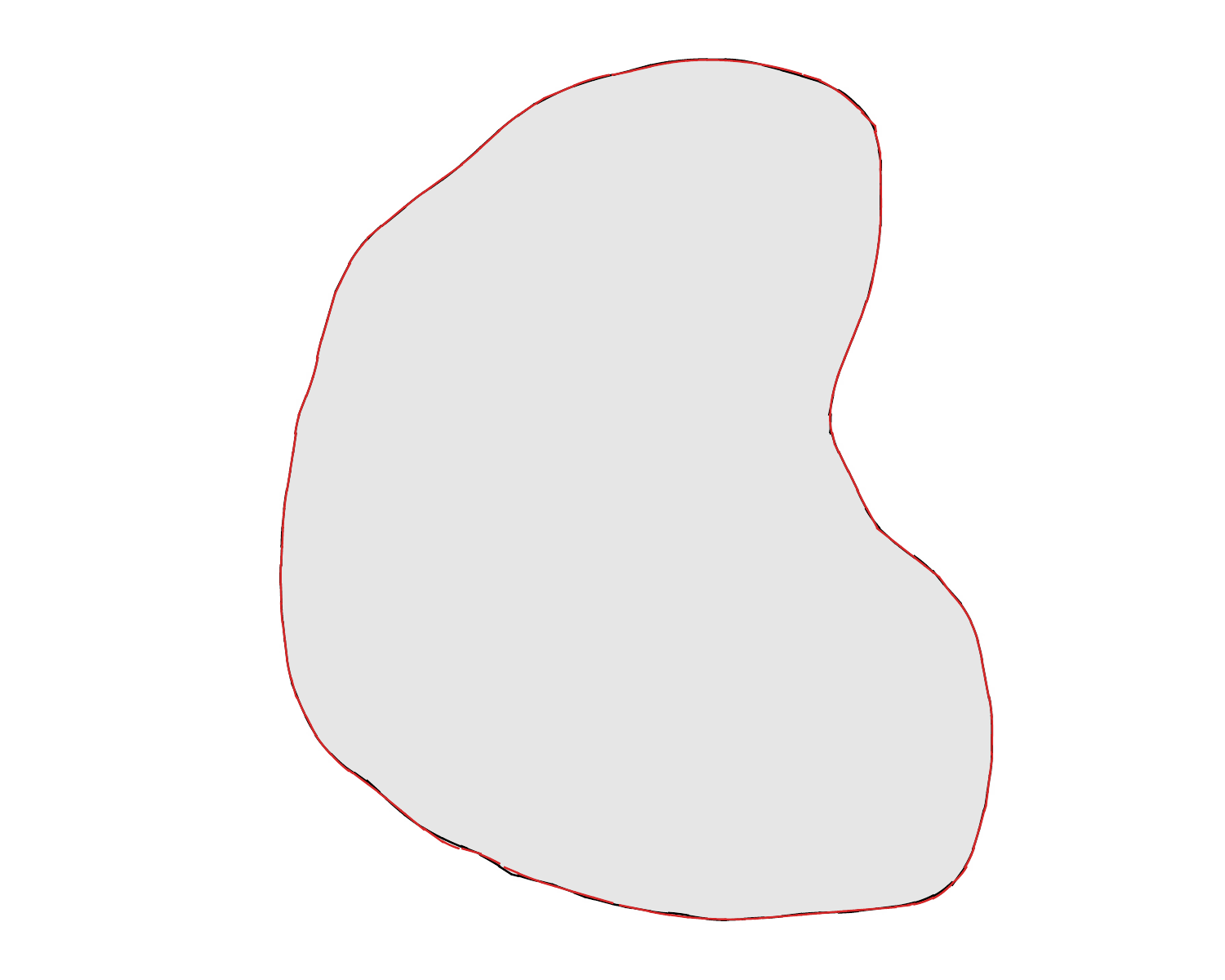}
         \caption{\corr{\perplexityinsert{aero-qelvira-vertex} and $h=1/60$.}}
     \end{subfigure}
	\caption{\corr{Reconstructions of a smooth domain at initial time (black) and final time (red).}}
	\label{fig:scheme-batata}
\end{figure}

\begin{figure}
     \centering
     \begin{subfigure}[t]{0.32\textwidth}
         \centering
         \includegraphics[width=\textwidth]{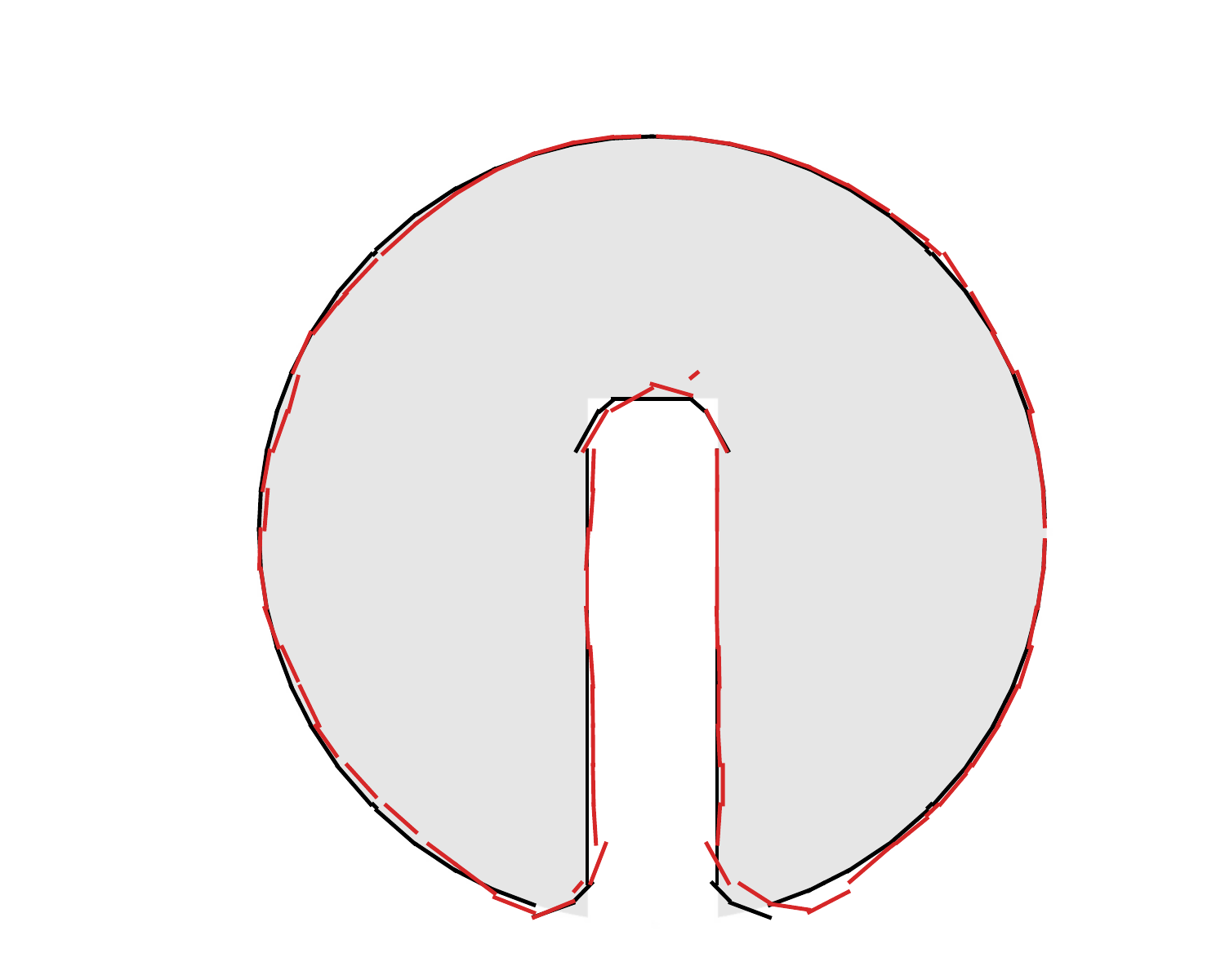}
         \caption{\corr{\perplexityinsert{elvira-w-oriented} and $h=1/30$.}}
     \end{subfigure}
     \hfill
     \begin{subfigure}[t]{0.32\textwidth}
         \centering
         \includegraphics[width=\textwidth]{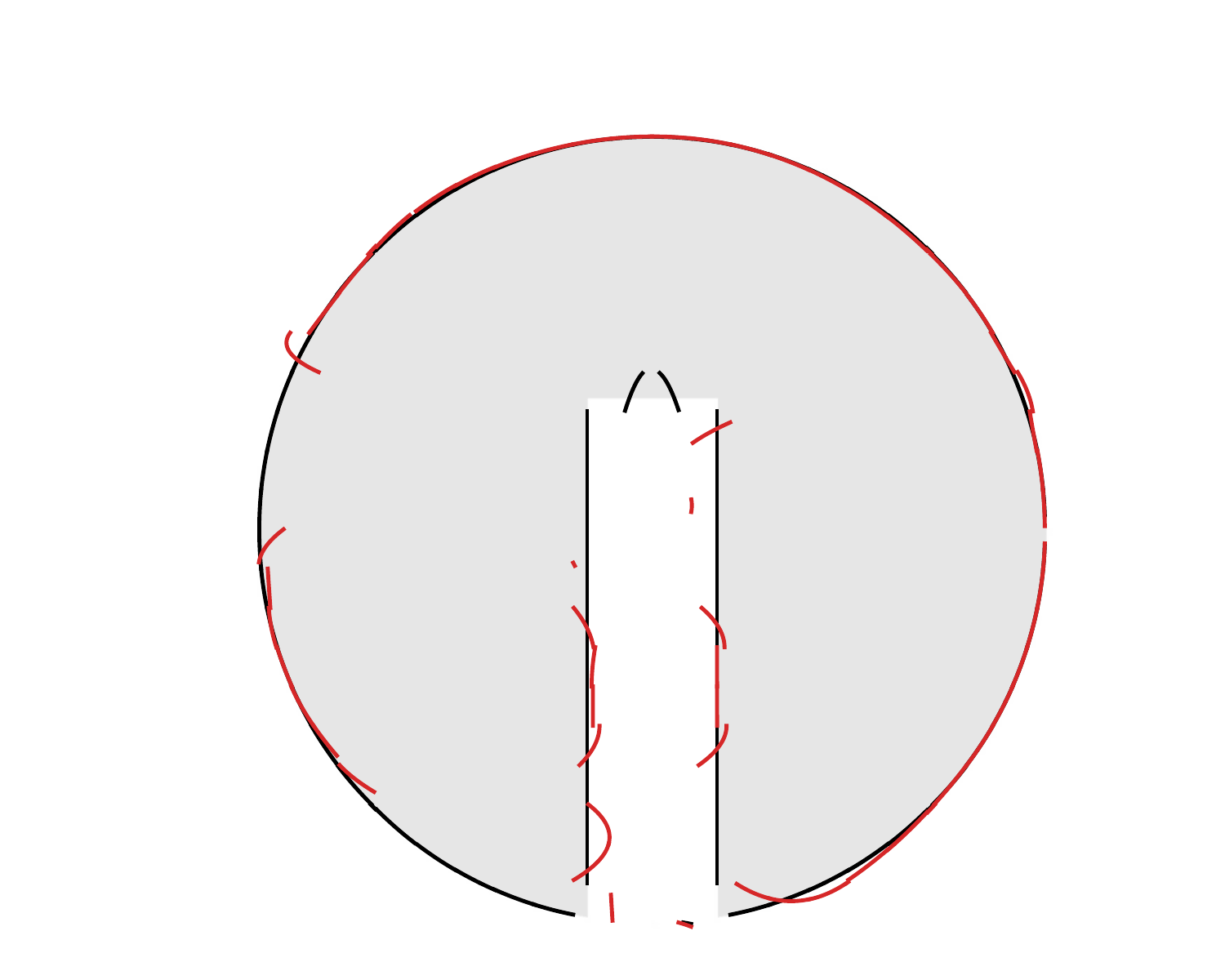}
         \caption{\corr{\perplexityinsert{quadratic-aero} and $h=1/30$.}}
     \end{subfigure}
     \hfill
     \begin{subfigure}[t]{0.32\textwidth}
         \centering
         \includegraphics[width=\textwidth]{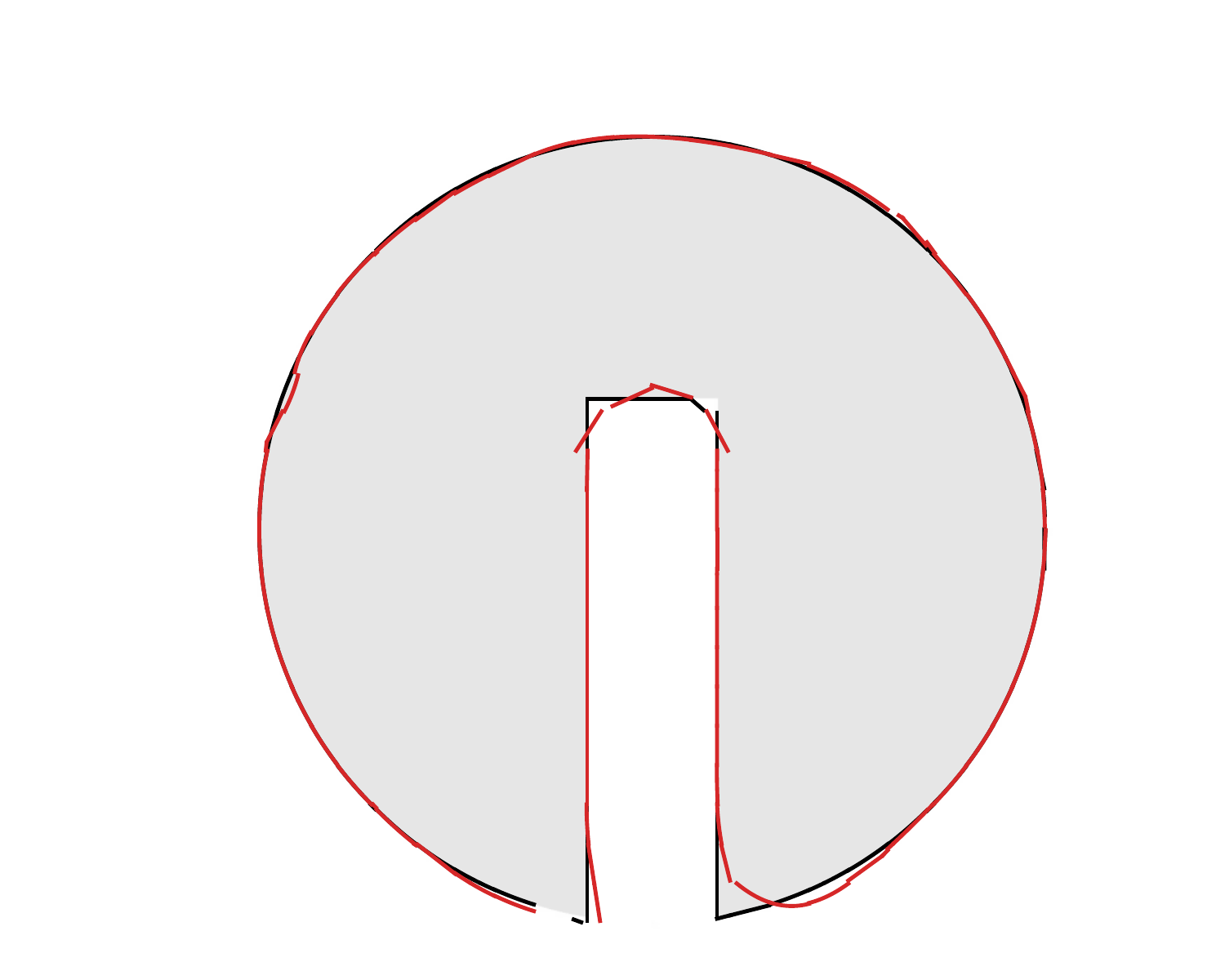}
         \caption{\corr{\perplexityinsert{aero-qelvira-vertex} and $h=1/30$.}}
     \end{subfigure}
     \vfill
          \begin{subfigure}[t]{0.32\textwidth}
         \centering
         \includegraphics[width=\textwidth]{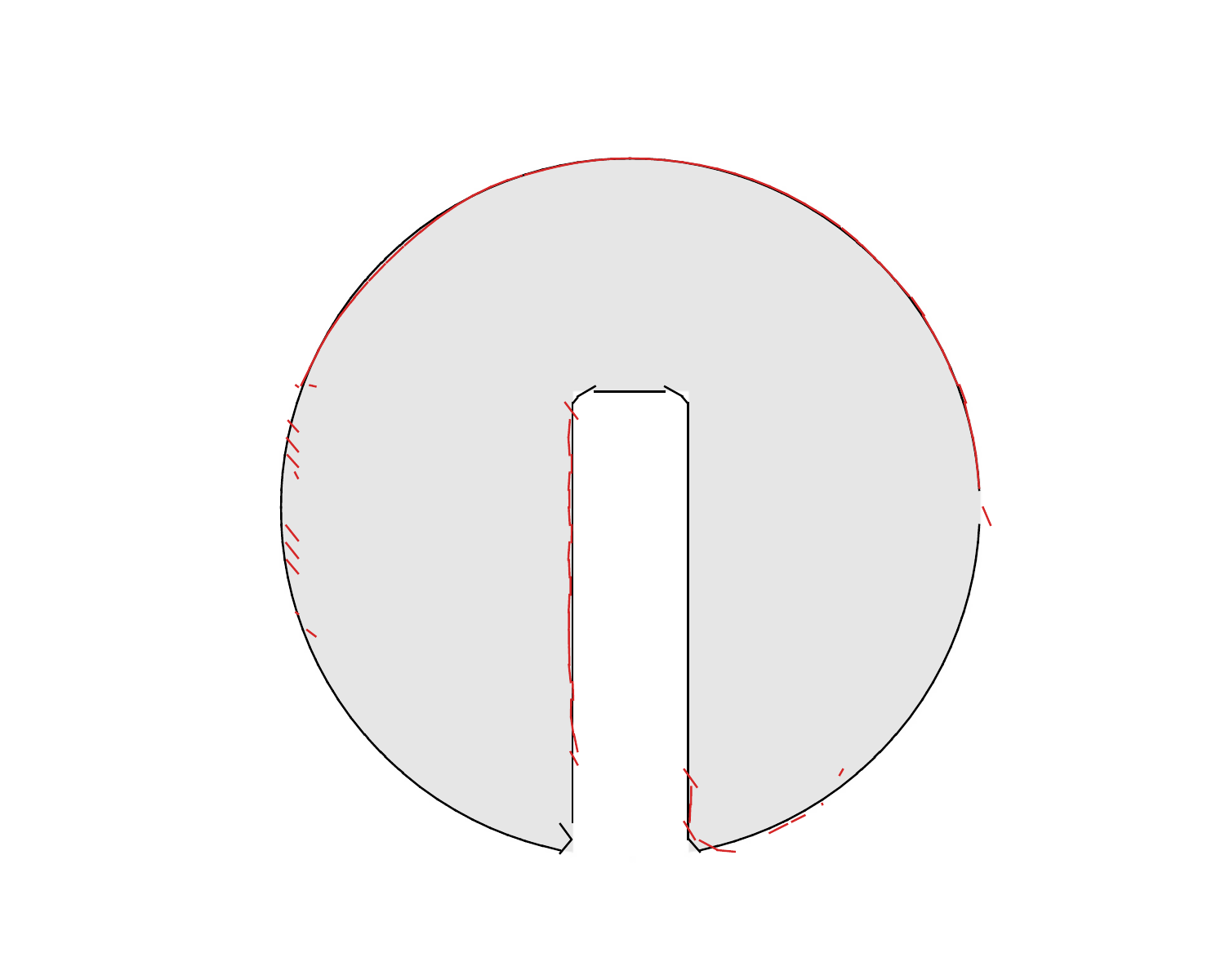}
         \caption{\corr{\perplexityinsert{elvira-w-oriented} and $h=1/60$.}}
     \end{subfigure}
     \hfill
     \begin{subfigure}[t]{0.32\textwidth}
         \centering
         \includegraphics[width=\textwidth]{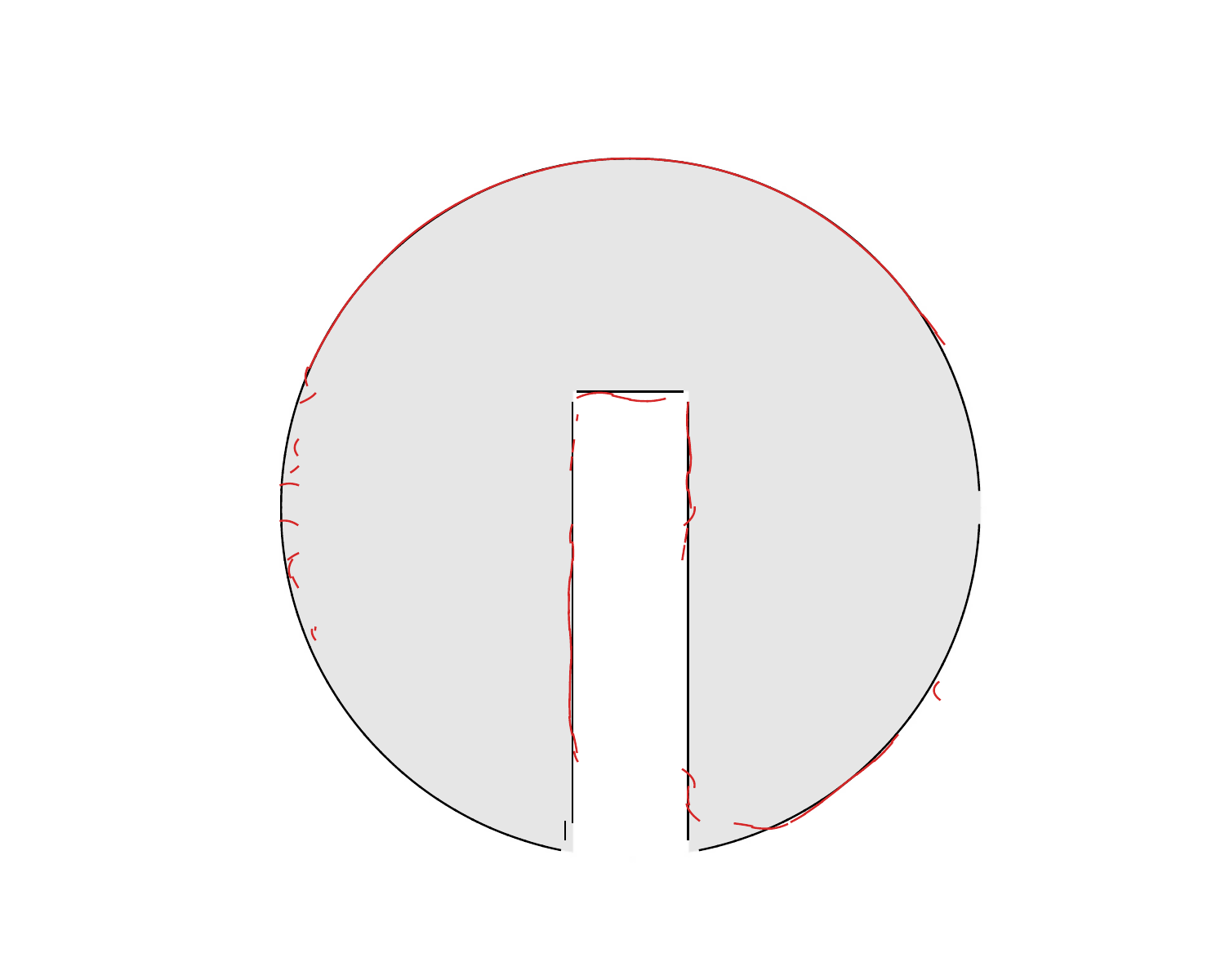}
         \caption{\corr{\perplexityinsert{quadratic-aero} and $h=1/60$.}}
     \end{subfigure}
     \hfill
     \begin{subfigure}[t]{0.32\textwidth}
         \centering
         \includegraphics[width=\textwidth]{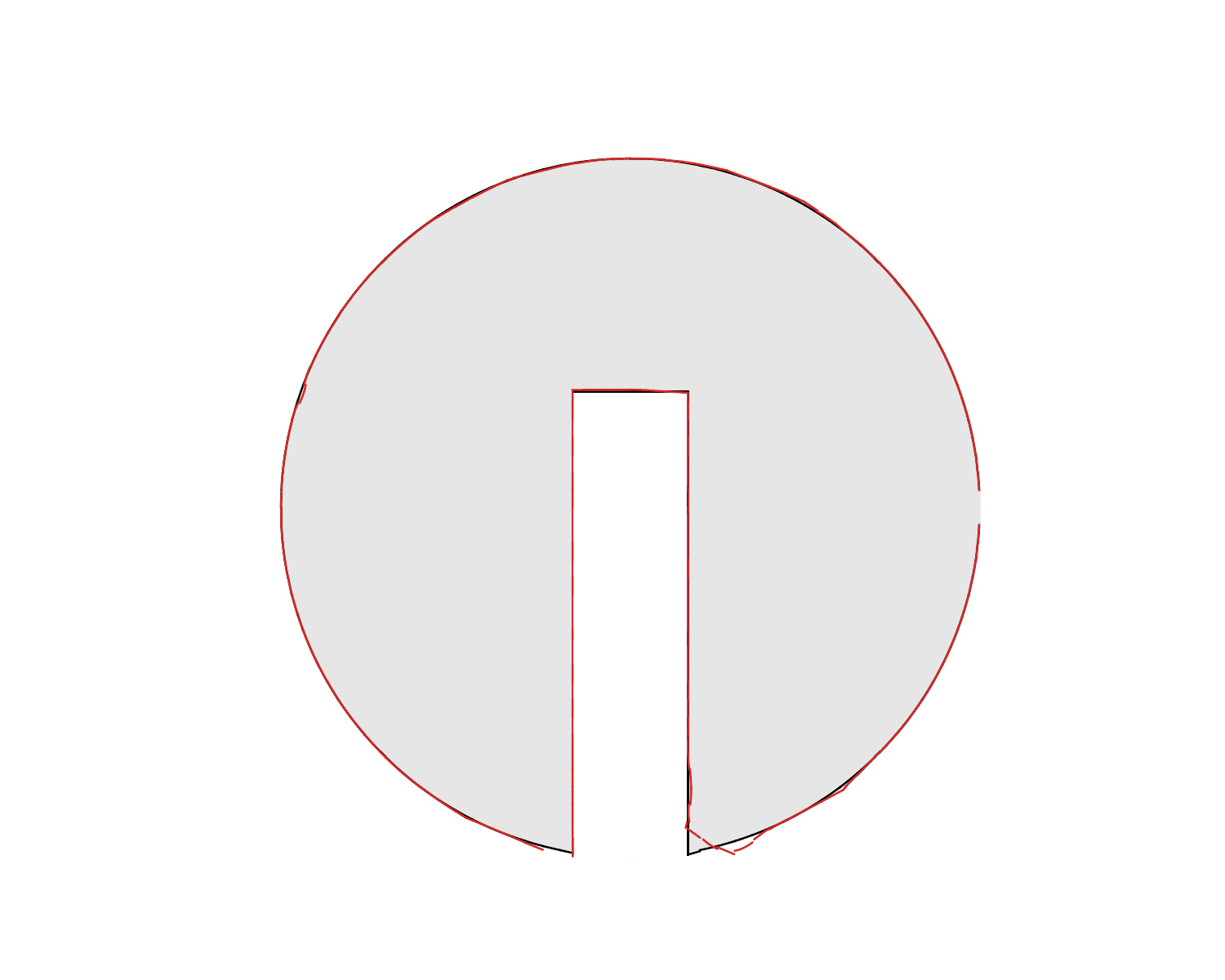}
         \caption{\corr{\perplexityinsert{aero-qelvira-vertex} and $h=1/60$.}}
     \end{subfigure}
	\caption{\corr{Reconstructions of the Zalesak notched circle at initial time (black) and final time (red).}}
	\label{fig:scheme-Zalesak}
\end{figure}

In this section we present two other Finite Volumes tests:
\begin{itemize}
	\item Smooth domain: in Figure \ref{fig:scheme-error-batata} we see that the three reconstruction methods keep the error low, limiting the diffusion. In Figure \ref{fig:scheme-batata} the reconstructions at start and after evolution show that the shape is mostly preserved.
	\item Zalesak notched circle: in Figure \ref{fig:scheme-error-batata} we see again that the presence of corners can be addressed by \perplexityinsert{aero-qelvira-vertex} and not by the other methods. This is evidenced even more by looking at the reconstructions in Figure \ref{fig:scheme-Zalesak}.
\end{itemize}

\begin{figure}
     \centering
     \begin{subfigure}[b]{\textwidth}
         \centering
         \caption{\corr{$h=1/30$}.}
         \includegraphics[width=\textwidth]{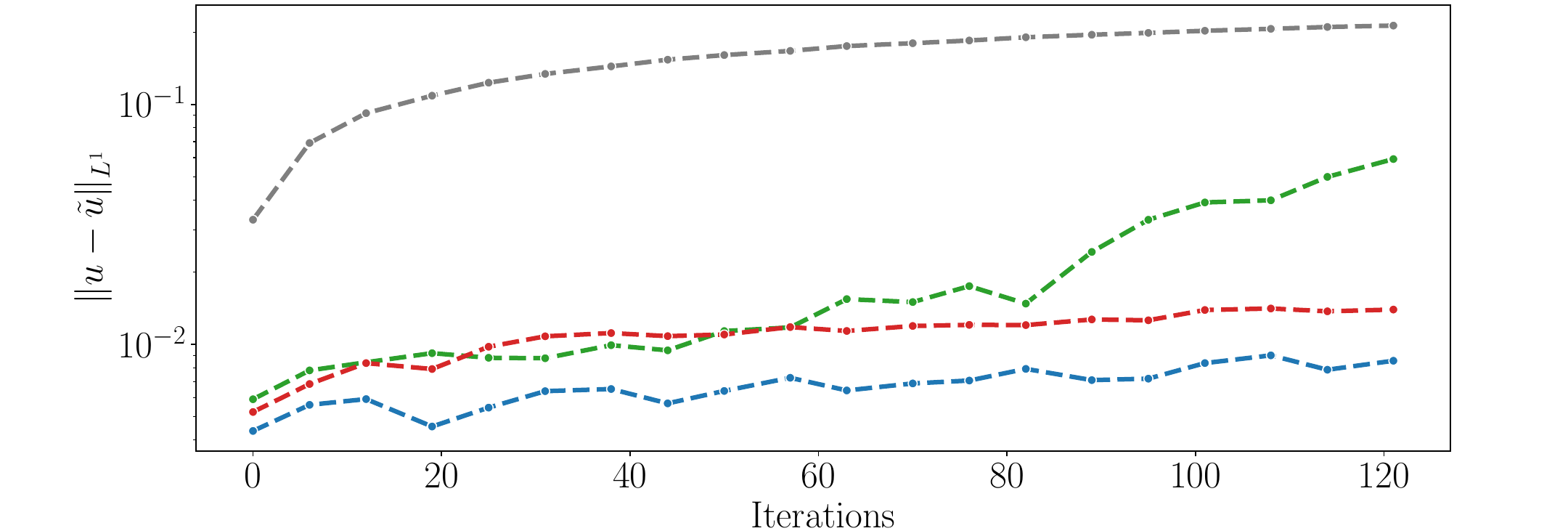}
         \label{fig:scheme-corners-30}
     \end{subfigure}
     \hfill
     \begin{subfigure}[b]{\textwidth}
         \centering
         \caption{\corr{$h=1/60$}.}
         \includegraphics[width=\textwidth]{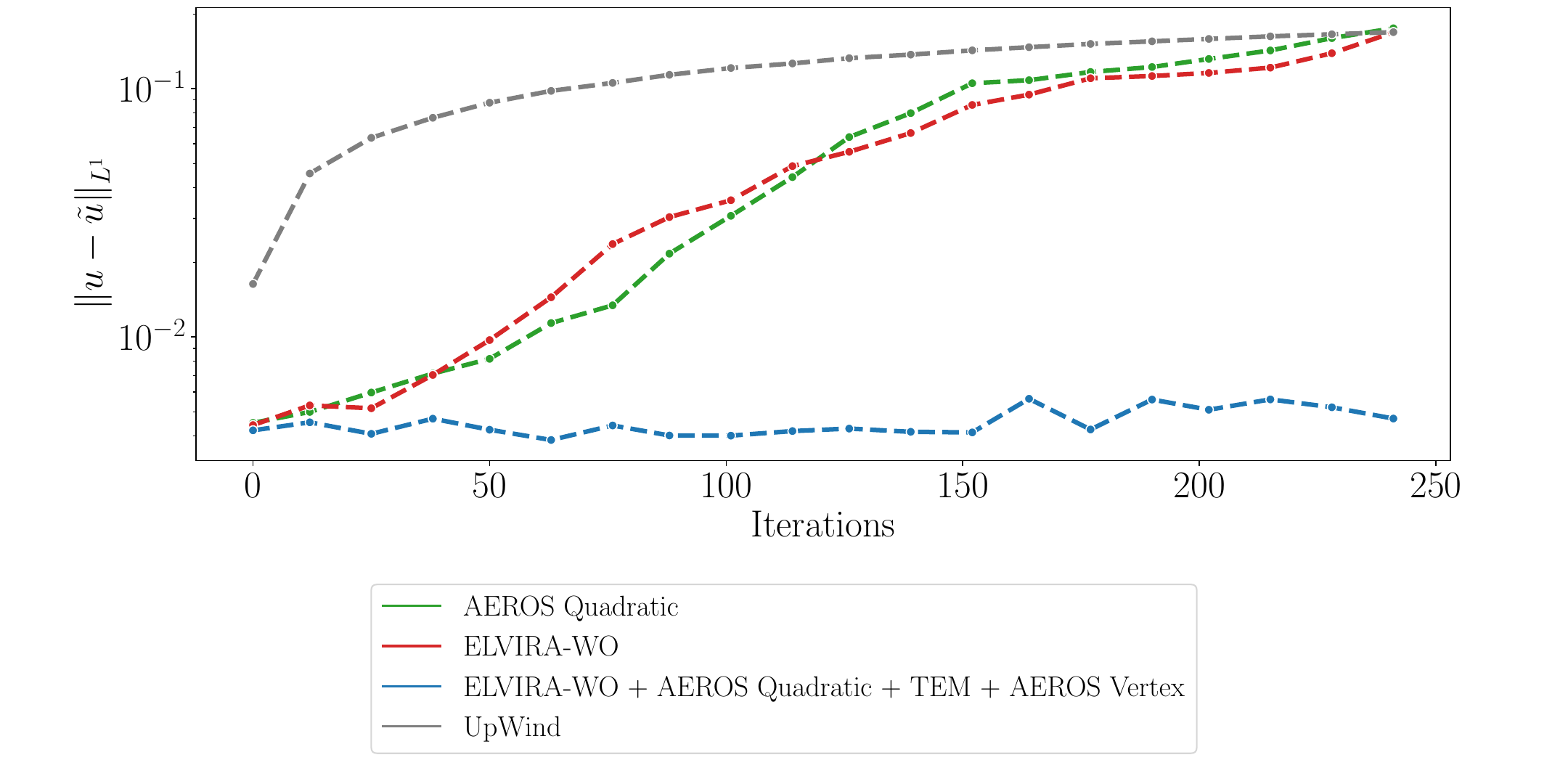}
         \label{fig:scheme-corners-60}
     \end{subfigure}
	\caption{Time evolution of the finite volume scheme $L^1$ error for the Zalesak notched circle.}
    \label{fig:scheme-error}
\end{figure}

\end{document}